\documentclass[12pt]{article}
\usepackage{amssymb}
\usepackage{tikz}
\usepackage{mathrsfs}
\def\build#1_#2^#3{\mathrel{\mathop{\kern 0pt#1}\limits_{#2}^{#3}}}\def\rde{\mathscr}
\def\Z{{\bf Z}}
\def\R{{\rde R}}\def\P{{\rde P}\!}\def\E{{\rde E}}\def\F{{\rde F}}\def\Mod{{\rde{M}\!od}}\def\M{{\rde M}}\def\U{{\rde U}}\def\G{{\rde G}}\def\H{{\rde H}}
\def\A{{\rde A}}\def\B{{\rde B}}\def\C{{\rde C}}\def\Nil{{\rde N}\! il}\def\Cl{{$\overline{\rm Cl}\ $}}\def\WE{{\bf WE}}\def\Ex{{\bf E}}
\def\Hom{{\rm Hom}}\def\Map{{\rm Map}}\def\Ker{{\rm Ker}}\def\Coker{{\rm Coker}}\def\Im{{\rm Im}}\def\fl{\longrightarrow}
\def\cqfd{\hfill\vbox{\hrule\hbox{\vrule height6pt depth0pt\hskip 6pt \vrule height6pt}\hrule\relax}}
\def\noi{\noindent}
\def\hfl#1#2{\smash{\mathop{\hbox to 12 mm{\rightarrowfill}}\limits^{\scriptstyle#1}_{\scriptstyle#2}}}
\def\hflp#1#2{\smash{\mathop{\hbox to 5 mm{\rightarrowfill}}\limits^{\scriptstyle#1}_{\scriptstyle#2}}}
\def\vfl#1#2{\llap{$\scriptstyle #1$}\left\downarrow\vbox to 6mm{}\right.\rlap{$\scriptstyle #2$}}
\def\bhfl#1#2{\smash{\mathop{\hbox to 12 mm{\leftarrowfill}}\limits^{\scriptstyle#1}_{\scriptstyle#2}}}
\def\diagram#1{\def\normalbaselines{\baselineskip=0pt\lineskip=10pt\lineskiplimit=1pt} \matrix{#1}}
\def\pv{\raise 2pt\hbox{$\bigwedge$}}\def\v{{}^\vee}\def\bv{{}^\wedge}
\begin{document}

\overfullrule=0pt
\
\vskip 64pt
\centerline{\bf Regular exact categories and algebraic K-theory}
\vskip 12pt
\centerline{Pierre Vogel\footnote{Universit\'e Paris Cit\'e and Sorbonne Universit\'e, CNRS, IMJ-PRG, F-75013 Paris, France}}
\vskip 48pt
\noi{\bf Abstract.} We introduce a new notion of regularity for rings and exact categories and we show important results in algebraic K-theory. In particular we prove 
a strong vanishing theorem for Nil groups and give an explicit class of groups, much bigger than Waldhausen's class Cl, such that every group in this class has trivial 
Whitehead groups.
\vskip 12pt
\noi{\bf Keywords:} Regularity, algebraic K-theory, functor Nil, Whitehead groups.

\noi{\bf Mathematics Subject Classification (2020):} 18E10, 19Dxx
\vskip 24pt
\noi{\bf Introduction.}
\vskip 12pt
Quillen's construction associates to any essentially small exact category $\A$ its algebraic K-theory which is an infinite loop space $K(\A)$ and this
correspondence is a functor from the category of essentially small exact categories to the category $\Omega sp_0$ of infinite loop spaces.

If $\A$ is the category $\P_A$ of finitely generated projective right modules over some ring $A$, we get a functor $A\mapsto K(A)=K(\P_A)$ and this functor can be
enriched into the non connective spectrum $\underline K$ containing also the negative part of the algebraic K-theory (see [B], [Ka]). That is $\underline K$ is a functor 
from the category of rings to the category $\Omega sp$ of $\Omega$-spectra and the natural transformation $K(A)\rightarrow \underline K(A)$ induces a homotopy equivalence 
from $K(A)$ to the connective part of $\underline K(A)$.

Denote by a left-flat bimodule a pair $(A,S)$ where $A$ is a ring and $S$ is a $A$-bimodule flat on the left. For each left-flat bimodule $(A,S)$ one has an
exact category $\Nil(A,S)$ where the objects are the pairs $(M,\theta)$ where $M$ is an object in $\P_A$ and $\theta:M\rightarrow M\otimes_A S$ is a nilpotent morphism
of right $A$-modules.

The correspondence $M\mapsto (M,0)$ induces a morphism $K(A)\rightarrow K(\Nil(A,S))$ and this morphism has a retraction coming from the functor $(M,\theta)\mapsto M$. Thus
there is a functor $Nil$ from the category of left-flat bimodules to $\Omega sp_0$ which is unique up to homotopy such that:
$$K(\Nil(A,S))\simeq K(A)\times Nil(A,S)$$ 

It is proven in [V] the following:
\vskip 12pt
\noi{\bf Theorem:} {\sl There is a functor $\underline N il$ from the category of left-flat bimodules to the category $\Omega sp$ of $\Omega$-spectra and a natural
transformation $Nil\rightarrow\underline Nil$ such that the following holds for every left-flat bimodule $(A,S)$:

$\bullet$ the map $Nil(A,S)\rightarrow\underline Nil(A,S)$ induces a homotopy equivalence from $Nil(A,S)$ to the connective part of $\underline Nil(A,S)$

$\bullet$ if $R$ is the tensor algebra of $S$, then there is a homotopy equivalence in $\Omega sp$:$$\underline K(R)\build\longrightarrow_{}^\sim \underline K(A)\times\Omega^{-1}(\underline Nil(A,S))$$

Moreover if $A$ is regular coherent on the right, the spectrum $\underline Nil(A,S)$ is contractible.}
\vskip 12pt
This functor $\underline Nil$ plays an important role in algebraic K-theory (see [V], [W1]), but essentially when the ring $A$ is regular coherent. 

If the ring $A$ is not coherent, very little is known about $\underline Nil(A,S)$.

In this paper we introduce a new notion of regularity for rings and, more generally, for exact categories. And this notion appears to be very useful in algebraic K-theory.
\vskip 12pt
Let $\M$ be a Grothendieck category and $\C$ be a class of objects in $\M$. We say that $\C$ is exact if, for every short exact sequence in $\M$:
$$0\fl X\fl Y\fl Z\fl 0$$
such that two of the modules $X,Y,Z$ are in $\C$, then the third one is also in $\C$.

We say that $\C$ is cocomplete if $\C$ is stable under filtered colimit.

Consider a ring $A$. Let $\Mod_A$ be the category of right $A$-modules and $\P_A$ be the subcategory of finitely generated projective modules in $\Mod_A$.

We say that $A$ is regular on the right, if $\Mod_A$ is the only cocomplete exact class in $\Mod_A$ containing $\P_A$. We have also
a notion of regularity on the left by considering left modules instead of right modules.

This definition of regularity seems to be very different from the classical one but we have the following:
\vskip 12pt
\noi{\bf Theorem 1:} {\sl A ring is regular coherent on the right if and only if it is regular on the right and coherent on the right.}
\vskip 12pt
\noi{\bf Remark:} This theorem remains true if right is replaced by left.
\vskip 12pt
The main result of this paper is the following vanishing theorem:
\vskip 12pt
\noi{\bf Vanishing Theorem 2:} {\sl Let $(A,S)$ be a left-flat bimodule. Suppose $A$ is regular on the right. Then the spectrum  $\underline Nil(A,S)$ is contractible.}
\vskip 12pt
An immediate consequence of this result is the following: in the paper [V], all results remain true if 'regular coherent' is replaced by 'regular'.
\vskip 12pt
Consider the smallest class of groups \Cl  satisfying the following:

$\bullet$ \Cl contains the trivial group

$\bullet$ let $H$, $G_1$ and $G_2$ be three groups in \Cl, and $\alpha:H\fl G_1$ and $\beta:H\fl G_2$ be two injective morphisms. Then the amalgamated free
product $\Gamma$ defined by the cocartesian diagram:
$$\diagram{H&\hfl{\alpha}{}&G_1\cr\vfl{\beta}{}&&\vfl{}{}\cr G_2&\hfl{}{}&\Gamma\cr}$$
belongs to \Cl

$\bullet$ let $H$ and $G$ be two groups in \Cl and $\alpha:H\fl G$ and $\beta:H\fl G$ be two injective morphisms. Then the corresponding HNN extension of $G$
belongs to \Cl

$\bullet$ if a group $G$ is the union of a filtered family of subgroups in \Cl then $G$ belongs to \Cl.
\vskip 12pt
This class of groups contains Waldhausen's class Cl, it is stable under subgroup and extension. Moreover, for every group $G$ in \Cl, the group ring
$\Z[G]$ is regular on the right (and on the left).
\vskip 12pt
\noi{\bf Theorem 3:} {\sl Let $A$ be a ring and $G$ be a group in \Cl. Suppose $A$ is regular on the right. Then the $\Omega$-spectrum $\underline Wh^A(G)$ is 
contractible.}
\vskip 12pt
Because of this result, it seems to be convenient to propose the following conjecture:
\vskip 12pt
\noi{\bf Conjecture:} {\sl Let $A$ be a ring and $G$ be a group. Suppose $A$ and $\Z[G]$ are regular on the right. Then the $\Omega$-spectrum $\underline Wh^A(G)$ is 
contractible.}
\vskip 12pt
In order to prove all these results we have to consider many exact categories. Thus we have to extend the definition of regularity to essentially small exact categories.
\vskip 12pt
If we want to extend the notion of regularity of a ring $A$ to the regularity of an essentially small exact category $\A$, we have to replace $\P_A$ by $\A$ and the
inclusion $\P_A\subset\Mod_A$ by an inclusion $\A\subset \A\v$, where $\A\v$ is an abelian category which satisfies certain properties.

In fact, for each essentially small exact category $\A$ we introduce a notion of strict cocompletion of $\A$, which is a Grothendieck category $\B$ containing $\A$
as a fully exact subcategory and satisfying certain properties. This notion is important because of the following result:
\vskip 12pt
\noi{\bf Theorem 4:} {\sl Let $\A$ be an essentially small exact category. Then $\A$ has a strict cocompletion which is unique up to equivalence. Moreover if $\A$ and $\B$
are essentially small exact categories and $\A\subset\A\v$ and $\B\subset\B\v$ are strict cocompletions of $\A$ and $\B$, then every exact functor $F:\A\fl\B$ has an
extension $F\v:\A\v\fl\B\v$, unique up to isomorphism, such that $F\v$ is right-exact and respects infinite direct sums.}
\vskip 12pt
Now we can extend the regularity condition for rings to regularity for essentially small exact categories. It suffices to replace the inclusion $\P_A\subset\Mod_A$ (for
some ring $A$) by the inclusion $\A\subset \A\v$ of an essentially small exact category $\A$ into its strict cocompletion $\A\v$.

Thus $\A$ is said to be regular if $\A\v$ is the only cocomplete exact class in $\A\v$ containing $\A$.
\vskip 12pt
A crucial point in this paper, which is needed in order to prove all these results, is the following theorem:
\vskip 12pt
Let $F:\A\fl\B$ be an exact functor between two essentially small exact categories. We say that $F$ is a domination if the following conditions are satisfied:

1) for any $Y$ in $\B$, there is an object $X\in\A$ and a deflation $F(X)\fl Y$

2) for any $X$ in $\A$, any $Y$ in $\B$ and any deflation $f:Y\fl F(X)$, there is an object $Z\in\A$, a morphism $g:F(Z)\fl Y$ and a deflation $h:Z\fl X$ such that: 
$F(h)=f\circ g$.
\vskip 12pt
\noi{\bf Domination Theorem 5:} {\sl Let $\B$ be an essentially small exact category and $\A$ be a fully exact subcategory of $\B$. Suppose $\A$ is regular and the 
inclusion $\A\subset\B$ is a domination. Then $\B$ is regular and the morphism $K_i(\A)\fl K_i(\B)$ induced by the inclusion is injective for $i=0$ and bijective for 
$i>0$. 

Moreover if $\A$ is stable under direct summand in $\B$, then $K(\A)\fl K(\B)$ is a homotopy equivalence and any object in $\B$ has a finite $\A$-resolution.}
\vskip 12pt
We have also the following result:
\vskip 12pt
\noi{\bf Resolution Theorem 6:} {\sl Let $\A$ be an essentially small exact category and $\A\v$ be a strict cocompletion of $\A$. Suppose $\A$ is regular. Then, for any
$\A$-resolution $C$ of some module in $\A\v$, the kernel of $d:C_n\fl C_{n-1}$ is, for $n$ large enough, a direct summand of a module in $\A$.}
\vskip 12pt
There is a last result in homological algebra for regular rings:
\vskip 12pt
\noi{\bf Theorem 7:} {\sl Let $A$ be a ring. Suppose $A$ is regular on the left. Then any acyclic complex of projective right $A$-modules is contractible.} 
\vskip 12pt

The paper is organized as follows:

The section 1 is devoted to the construction of cocompletions for exact categories and the proof of Theorem 4.

In section 2, theorem 5 and 6 are proven. As a consequence we get a proof of Theorem 1. We prove also Theorem 7.

In section 3, we prove many properties of regular rings and regular groups. In particular, we show that every group in $\overline{Cl}$ is regular and prove Theorem 3.

The last section is devoted to the proof of Theorem 2.
\vskip 24pt
\noi{\bf 1. Cocompletion of exact categories.}
\vskip 12pt
\noi{\bf Writing conventions:}

If $\A$ is an additive category, objects and morphisms in $\A$ will be called $\A$-modules and $\A$-morphisms.
 
Covariant functors and contravariant functors will be called functors and cofunctors.

If $\Phi:\A\fl\B$ is a functor between two categories and $f:A\fl Y$ is a morphism in $\A$, its induced morphism $\Phi(\A)\fl\Phi(\B)$ will be often denoted $f$.
Thus a morphism $f:\A\fl\B$ induces a morphism $f:\Phi(\A)\fl \Phi(\B)$.
\vskip 12pt
Recall that an exact category (in the sense of Quillen) is an additive category $\A$ together with two subcategories of $\A$: the category of inflations of $\A$ and the 
category of deflations on $\A$ (with the notations of [GR]) such that the following holds:

$\bullet$ each isomorphism in $\A$ is an inflation and a deflation.

$\bullet$ each inflation $X\fl Y$ is a monomorphism and has a cokernel $Z$ such that $Y\fl Z$ is a deflation.

$\bullet$ each deflation $Y\fl Z$ is an epimorphism and has a kernel $X$ such that $X\fl Y$ is an inflation.

$\bullet$ (cobase change) for any inflation $f:X\fl Y$ and any morphism $g:X\fl X'$ there is a cocartesian diagram:
$$\diagram{X&\hfl{f}{}&Y\cr\vfl{g}{}&&\vfl{}{}\cr X'&\hfl{}{}&Y'\cr}$$
where $X'\fl Y'$ is an inflation.

$\bullet$ (base change) for any deflation $f:X\fl Y$ and any morphism $g:Y'\fl Y$, there is a cartesian diagram:
$$\diagram{X'&\hfl{}{}&Y'\cr\vfl{}{}&&\vfl{g}{}\cr X&\hfl{f}{}&Y\cr}$$
where $X'\fl Y'$ is a deflation.

A sequence $X\fl Y\fl Z$ is called a conflation if $X\fl Y$ is an inflation with cokernel Z (or equivalently if $Y\fl Z$ is a deflation with kernel X). A short exact
sequence in $\A$ is a sequence $0\fl X\fl Y\fl Z\fl 0$ in $\A$ such that $X\fl Y\fl Z$ is a conflation in $\A$. 

If necessary, if $\A$ is an exact category, objects, morphisms, inflations, deflations, conflations  and exact sequences in $\A$ will be called $\A$-modules, 
$\A$-morphisms, $\A$-inflations, $\A$-deflations, $\A$-conflations and short $\A$-exact sequences.

An abelian category is an exact category, by taking as inflations and deflations the monomorphisms and the epimorphisms of the category.

An exact functor $\Phi:\A\fl\B$ between two exact categories is a functor that respects conflations.

An exact functor $\Phi:\A\fl\B$ is said to be fully exact if it is fully faithful and if, for any sequence $S=(X\fl Y\fl Z)$ in $\A$, $S$ is an $\A$-conflation if and
only if $\Phi(S)$ is a $\B$-conflation.
\vskip 12pt
\noi{\bf 1.1 Lemma:} {\sl Let $\A$ be an exact category and $\C$ be a class of $\A$-modules. Suppose $\C$ is stable under extension. Then there is a unique fully exact
subcategory $\B$ of $\A$ such that $\C$ is the class of $\B$-modules.}
\vskip 12pt
\noi{\bf Proof:} Since $\C$ is stable under extension, it is easy to see that $\C$ is stable under base change and cobase change. Therefore $\C$ generates a well
defined fully exact subcategory of $\A$. The result follows.\cqfd
\vskip 12pt
\noi{\bf Definition:} Let $\A$ be an exact category. We say that a functor $F:\A\fl\B$ is a cocompletion if the following conditions are satisfied:

(cc0) $\B$ is a cocomplete abelian category
 
(cc1) $F$ is fully exact

(cc2) for any $\A$-module $X$, $F(X)$ is finitely presented in $\B$

(cc3) for any non-zero $\B$-module $Y$, there is an $\A$-module $X$ and a non-zero morphism $F(X)\fl Y$

(cc4) for any $\A$-module $X$, any $\B$-epimorphism $\varphi:U\fl V$ and any $\B$-morphism $f:F(X)\fl V$, there exist an $\A$-deflation $g:Y\fl X$ and a commutative 
diagram in $\B$:
$$\diagram{F(Y)&\hfl{}{}&U\cr\vfl{g}{}&&\vfl{\varphi}{}\cr F(X)&\hfl{f}{}&V\cr}$$
\vskip 12pt
If $F$ is an inclusion, we say that $\B$ is a strict cocompletion of $\A$.
\vskip 12pt
\noi{\bf 1.2 Theorem:} {\sl Every essentially small exact category has a strict cocompletion which is unique up to equivalence. Moreover, if $F:\A\fl\B$ is a cocompletion 
of an essentially small exact category $\A$, the following holds:

(cc5) $\B$ is a Grothendieck category

(cc6) for any $\B$-module $Y$, there is a family of $\A$-modules $X_i$ indexed by a set and an epimorphism $\build\oplus_i^{}F(X_i)\fl Y$

(cc7) $\A$ is stable under extension in $\B$, that is: for any $\B$-conflation $F(X)\fl E\fl\Phi(Y)$ with $X$ and $Y$ in $\A$, there is an $\A$-conflation 
$X\fl Z\fl Y$ and a commutative diagram:
$$\diagram{F(X)&\hfl{}{}&F(Z)&\hfl{}{}&F(Y)\cr\vfl{=}{}&&\vfl{\simeq}{}&&\vfl{=}{}\cr F(X)&\hfl{}{}&E&\hfl{}{}&F(Y)\cr}$$
where $F(Z)\fl E$ is an isomorphism

(cc8) for any $\A$-morphism $f:X\fl Y$, the morphism $f:F(X)\fl F(Y)$ is an epimorphism if and only if there is an $\A$-module $Z$ such that the morphism 
$X\oplus Z\build\fl_{}^{f\oplus0}Y$ is a deflation.}
\vskip 12pt
\noi{\bf Proof:} Consider an essentially small exact category $\A$.
\vskip 12pt
\noi{\bf 1.3 Lemma:} {\sl  Let $F:\A\fl\B$ be a functor. Suppose conditions (cc0), (cc1), (cc3) and (cc4) are satisfied for $F$. Then the conditions (cc6), (cc7) and 
(cc8) are also satisfied.}
\vskip 12pt
\noi{\bf Proof:} 

\noi{\bf Condition (cc8):}

Consider an $\A$-morphism $f:X\fl Y$ inducing an epimorphism $f:F(X)\fl F(Y)$. Because of condition (cc4), there is a deflation $g:T\fl Y$ and a commutative diagram:
$$\diagram{F(T)&\hfl{}{}&F(X)\cr\vfl{g}{}&&\vfl{f}{}\cr F(Y)&\hfl{=}{}&F(Y)\cr}$$

Because of condition (cc1), this diagram comes from a commutative diagram:
$$\diagram{T&\hfl{h}{}&X\cr\vfl{g}{}&&\vfl{f}{}\cr Y&\hfl{=}{}&Y\cr}$$
and $fh=g:T\fl Y$ is a deflation.

Consider the morphism $\varphi=f\oplus fh:X\oplus T\fl Y$. This morphism is the composite of the morphism $\varphi_1=$Id$_X\oplus fh:X\oplus T\fl X\oplus Y$ and
the morphism $\varphi_2=f\oplus$Id$_Y:X\oplus Y\fl Y$. But $\varphi_1$ is a deflation because $fh$ is a deflation and $\varphi_2$ is also a deflation because it has a 
section. Then $\varphi$ is a deflation.

Let $Z$ be the kernel of Id$_X\oplus h:X\oplus T\fl X$. Then we have: $X\oplus T=X\oplus Z$ and the deflation $\varphi$ is the morphism $f\oplus 0:X\oplus Z\fl Y$.
Hence the condition (cc8) is satisfied.
\vskip 12pt
\noi{\bf Condition (cc7):}

Let $X$ and $Z$ be two $\A$-modules, $U$ be a $\B$-module and:
$$0\fl F(X)\build\fl_{}^\alpha U\build\fl_{}^\beta F(Z)\fl 0$$
be an exact sequence. Because of the condition (cc4), there is a deflation $f:T\fl Z$ and a morphism $g:F(T)\fl U$ such that: $f=\beta g$. Since $f$ is a deflation, the
morphism $X\oplus T\build\fl_{}^{0\oplus f} Z$ is also a deflation and its kernel is an $\A$-module $K$. Thus we have a commutative diagram with exact lines:
$$\diagram{0&\hfl{}{}&F(K)&\hfl{}{}&F(X)\oplus F(T)&\hfl{}{}&F(Z)&\hfl{}{}&0\cr&&\vfl{}{}&&\vfl{}{}&&\vfl{}{}&&\cr0&\hfl{}{}&F(X)&\hfl{\alpha}{}&U&\hfl{\beta}{}&F(Z)&
\hfl{}{}&0\cr}$$
Since $F$ is fully faithful, the morphism $F(K)\fl F(X)$ comes from an $\A$-morphism $\lambda:K\fl X$ and, because $K\fl X\oplus T$ is an inflation, we get, by cobase 
change, a conflation:
$$X\fl Y\fl Z$$
and a commutative diagram:
$$\diagram{F(X)&\hfl{}{}&F(Y)&\hfl{}{}&F(Z)\cr\vfl{=}{}&&\vfl{\simeq}{}&&\vfl{=}{}\cr F(X)&\hfl{}{}&U&\fl{}{}&F(Z)\cr}$$
Hence the condition (cc7) is satisfied.
\vskip 12pt
\noi{\bf Condition (cc6):}

Since $\A$ is essentially small it has a skeleton $\U$ and every $\A$-module is isomorphic to a module in $\U$.

Let $Y$ be a $\B$-module and $M$ be the set of pairs $(X,f)$ with $X\in\U$ and $f\in\Hom(F(X),Y)$. Each $(X,f)$ in $M$ induces a morphism $f:F(X)\fl Y$ and we get a 
morphism
$$\varphi:\build\oplus_{(X,f)}^{} F(X)\build\fl_{}^f Y$$
Let $Z$ be the cokernel of this morphism. Suppose $Z$ is non-zero. Then, because of the condition (cc3), there is an $\A$-module $T$ and a non-zero morphism $F(T)\fl Z$
and, because $Y\fl Z$ is epic, the condition (cc4) implies there is a deflation $T'\fl T$ and a commutative diagram:
$$\diagram{F(T')&\hfl{}{}&Y\cr\vfl{}{}&&\vfl{}{}\cr F(T)&\hfl{}{}&Z\cr}$$
Since $T'\fl T$ is a deflation, $F(T')\fl F(T)$ is epic and the composite morphism $F(T')\fl Z$ is non-zero. Since $\U$ is a skeleton of $\A$, $T'$ is isomorphic to
some $X\in \U$ and there is a pair $(X,f)\in M$ such that the composite morphism $F(X)\fl Y\fl Z$ is non-zero. But this fact contradicts the definition of $Z$. Therefore
we have: $Z=0$ and $\varphi$ is an epimorphism. Then condition (cc6) is satisfied and the lemma is proven.\cqfd 

\vskip 12pt
Let $\A$ be an essentially small exact category. Let Mod$_\A$ be the category of cofunctors from $\A$ to the category $Ab$ of abelian groups and $\A\v$
be the subcategory of left-exact cofunctors in Mod$_\A$.

We have a functor $F:\A\fl \A\v$ sending each $\A$-module $X$ to the cofunctor $X\v:Y\mapsto \Hom(Y,X)$. This functor is known to be the Gabriel-Quillen embedding. 
\vskip 12pt
\noi{\bf 1.4 Lemma:} {\sl The functor $F:\A\fl \A\v$ is a cocompletion of $\A$ satisfying all the properties of Theorem 1.2.}
\vskip 12pt
\noi{\bf Proof:} Because of Lemma 1.3, it suffices to check conditions (cc1), (cc2), (cc3), (cc4) and (cc5).

It is known that $\A\v$ is an abelian category and that $F:\A\fl \A\v$ is fully exact (see [Ga] II section 2, [Q] section 2, [Ke] Appendix A, [TT] 
Appendix A). It is also known the following property:

if $f:E\fl E'$ is an epimorphism in $\A\v$, then, for each $\A$-module $X$ and each element $x\in E'(X)$, there is a deflation $g:Y\fl X$ such that 
$g^*(x)\in f_*(E(Y))$. Then conditions (cc1) and (cc4) are satisfied.

Let $\C\subset$Mod$_\A$ be the full subcategory of cofunctors $F$ such that, for each $Z\in\A$ and each $z\in F(Z)$, there is a deflation $f:Y\fl Z$ such that
$f^*(z)=0$. The category $\C$ is a localizing subcategory of Mod$_\A$ and $\A\v$ is equivalent to Mod$\A/\C$ (see [Ke] Proposition A.2). Therefore $\A\v$ is a Grothendieck
category (see [Ga]) and property (cc5) is satisfied.

Consider a non-zero cofunctor $E$ in $\A\v$. Since $E\not=0$, there is an $\A$-module $X$ such that $E(X)\not=0$. But we have: $E(X)=\Hom(X\v,E)$ and we get a non-zero 
morphism $X\v\fl E$. Then condition (cc3) is satisfied.

Recall that an object $X$ in a cocomplete abelian category is finitely presented if and only if Hom$(X,-)$ commutes with filtered colimit.

Let $E_i, i\in I$ be a filtered system in $\A\v$ and $E$ its colimit. Thus we have compatible morphisms $\alpha_i:E_i\fl E$. Let $X$ be an $\A$-module. The morphisms 
$\alpha_i$ induce compatible morphisms $\Hom(X\v,E_i)\fl \Hom(X\v,E)$. Denote $F(X)$  the colimit of the modules $\Hom(X\v,E_i)$. Since 
$\build\lim_{{\fl\atop i}}^{}$ is exact, $F$ is left-exact and belongs to $\A\v$. Thus we get a morphism $f:F\fl E$ and compatible morphisms $\beta_i:E_i\fl F$.

It is easy to check that $f$ is the unique morphism such that $f\beta_i=\alpha_i$ for all $i\in I$. Therefore $f$ is an isomorphism and we get isomorphisms:
$$\build\lim_{{\fl\atop i}}^{}\Hom(X\v,E_i)\simeq F(X)\simeq \Hom(X\v,F)\simeq \Hom(X\v,E)\simeq \Hom(X\v,\build\lim_{{\rightarrow\atop i}}^{} E_i)$$
Hence $X\v$ is finitely presented in $\A\v$, and property (cc2) is satisfied.\cqfd
\vskip 12pt
\noi{\bf 1.5 Lemma:} {\sl The category $\A$ has a strict cocompletion satisfying all properties (cc*).}
\vskip 12pt
\noi{\bf Proof:} Let $\C$ be the disjoint union of the objects of $\A$ and the objects of $\A\v$. Then there is a unique way to put on $\C$ a structure of category such
that the inclusions $\A\subset\C$ and $\A\v\subset\C$ and the map $\pi=F\oplus$Id$:\C\fl\A\v$ are functors. Moreover $\pi$ is an equivalence of categories. Therefore
$\C$ is a strict cocompletion of $\A$ satisfying all properties (cc*).\cqfd
\vskip 12pt
Let $\Phi:\A\fl \B$ be a cocompletion of $\A$. If $E$ is a $\B$-module, we have a cofunctor $\A\fl$ Ab sending each $\A$-module $X$ to the $\Z$-module $\Hom(\Phi(X),E)$.
This cofunctor $\Psi(E)$ is left exact and belongs to $\A\v$. The correspondence $E\mapsto\Psi(E)$ is a functor $\Psi$ and we have:
$$\Psi(E)(X)=\Hom(\Phi(X),E)=\Hom(X\v,\Psi(E))$$
for any $\A$-module $X$ and any $\B$-module $E$. Note that all conditions (cc*) are satisfied by $\B$ and $\A\v$ except that $\B$ is not necessarily a Grothendieck
category. 
\vskip 12pt
\noi{\bf 1.6 Lemma:} {\sl The functor $\Psi:\B\fl\A\v$ is exact and respects direct sums.}
\vskip 12pt
\noi{\bf Proof:} Let $0\fl E\fl E'\fl E''\fl 0$ be a short exact sequence in $\B$. Let $F_0$ and $F_1$ be the kernel and the cokernel of $\Psi(E)\fl\Psi(E')$ and $F_2$ be
the cokernel of $\Psi(E')\fl\Psi(E'')$.

For each $\A$-module $X$ the composite morphism $F_0(X)\fl\Hom(\Phi(X),E)\fl\Hom(\Phi(X),E')$ is null. But the morphism $\Hom(\Phi(X),E)\fl\Hom(\Phi(X),E')$ is injective.
Therefore the morphism $F_0(X)\fl\Hom(\Phi(X),E)$ is null for each $X$ and the morphism $F_0\fl \Psi(E)$ is also null. Hence we have $F_0=0$ and $\Psi(E)\fl\Psi(E')$ is
a monomorphism.

Let $X$ be an $\A$-module. Consider the following sequence:
$$\Hom(\Phi(X),E')\fl \Hom(\Phi(X),E'')\fl F_2(X)$$
Let $x$ be an element of $F_2(X)$. This element corresponds to a morphism $f:X\v\fl F_2$. Since $\Psi(E'')\fl F_2$ is an $\A\v$-epimorphism, the property (cc4) implies 
that there exist a deflation $Y\fl X$ such that the composite morphism $Y\v\fl X\v\fl F_2$ factors through $\Psi(E'')$ by a morphism $Y\v\fl E''$ and, in the commutative
following diagram:
$$\diagram{\Hom(\Phi(X),E')&\hfl{}{}&\Hom(\Phi(X),E'')&\hfl{}{}& F_2(X)\cr\vfl{}{}&&\vfl{}{}&&\vfl{}{}\cr\Hom(\Phi(Y),E')&\hfl{}{}&\Hom(\Phi(Y),E'')&\hfl{}{}& F_2(Y)\cr}$$
the image of $x\in F_2(X)$ in $F_2(Y)$ lifts to an element $y\in\Hom(\Phi(Y),E'')$. Let $\overline E$ be the $\B$-module defined by the cartesian square:
$$\diagram{\overline E&\hfl{}{}&\Phi(Y)\cr\vfl{}{}&&\vfl{y}{}\cr E'&\hfl{}{}&E''\cr}$$
Since $E'\fl E''$ is an epimorphism $\overline E\fl \Phi(Y)$ is also an epimorphism and, because of condition (cc4), there is a deflation $Y'\fl Y$ such that the
morphism $\Phi(Y')\fl\Phi(Y)$ factors through $\overline E$. Then, up to replacing $Y$ by $Y'$, we may as well suppose that $x$ lifts to an element $y\in\Hom(\Phi(Y),E')$.

But the composite morphism $\Hom(\Phi(Y),E')\fl\Hom(\Phi(Y),E'')\fl F_2(Y)$ is null and the image of $x$ in $F_2(Y)$ is also null. On the other hand, the deflation
$Y\fl X$ induces a monomorphism $F_2(X)\fl F_2(Y)$ and we have: $x=0$. Hence $F_2(X)=0$ for each $X$ and the cokernel $F_2$ of $\Psi(E')\fl\Psi(E'')$ is null.
Therefore $\Psi(E')\fl\Psi(E'')$ is an epimorphism.  

Consider the following diagram:
$$\Hom(\Phi(X),E')\fl F_1(X)\fl \Hom(\Phi(X),E'')$$
Let $x$ be an element in $F_1(X)$ sent to $0$ in $\Hom(\Phi(X,E'')$. Since $\Psi(E')\fl F_1$ is an epimorphism, there is a deflation $Y\fl X$ and a 
commutative diagram:
$$\diagram{Y\v&\hfl{}{}&\Psi(E')\cr\vfl{}{}&&\vfl{}{}\cr X\v&\hfl{}{}&F_1\cr}$$
where $X\v\fl F_1$ corresponds to $x\in F_1(X)$. We have a commutative diagram:
$$\diagram{\Hom(\Phi(X),E')&\hfl{}{}&F_1(X)&\hfl{}{}&\Hom(\Phi(X),E'')\cr\vfl{}{}&&\vfl{}{}&&\vfl{}{}\cr\Hom(\Phi(Y),E')&\hfl{}{}&F_1(Y)&\hfl{}{}&\Hom(\Phi(Y),E'')\cr}$$
In this diagram, the image $y$ of $x$ in $F_1(Y)$ can be lifted to an element $y'\in\Hom(\Phi(Y),E')$ and, because this element goes to $0$ in $\Hom(\Phi(Y,E'')$, it comes
from an element $y\in\Hom(\Phi(Y),E)$.

On the other hand, the composite morphism $\Psi(E)\fl\Psi(E')\fl F_1$ is null and the composite morphism $\Hom(\Phi(Y),E)\fl\Hom(\Phi(Y),E')\fl F_1(Y)$ also. Therefore $y$
goes to $0$ in $F_1(Y)$ and, because $F_1(X)\fl F_1(Y)$ is a monomorphism, we have $x=0$. Hence $F_1(X)\fl\Hom(\Phi(X),E'')$ is a monomorphism and $F_1\fl\Psi(E'')$ is
also a monomorphism.

Hence, in the following sequence:
$$0\fl \Psi(E)\build\fl_{}^\alpha\Psi(E')\build\fl_{}^\beta\Psi(E'')\fl 0$$
$\alpha$ is a monomorphism, $\beta$ is an epimorphism and the morphism 
$$\Coker(\Psi(E)\build\fl_{}^\alpha\Psi(E'))\fl\Psi(E'')$$ 
is a monomorphism. Therefore the sequence is exact and $\Psi$ respects exact sequences.
\vskip 12pt
Let $E_i$, $i\in I$ be a family of $\B$-modules and $E$ be its direct sum. Let $X$ be an $\A$-module. Since $\Phi(X)$ is finitely presented, every morphism from $\Phi(X)$
in $E$ factors through a finite sum of the $X_i$'s and we have:
$$\Hom(\Phi(X),E)=\Hom(\Phi(X),\build\oplus_i^{} E_i)=\ \build\oplus_i^{}\Hom(\Phi(X),E_i)$$
Therefore $\Psi(E)$ is the direct sum of the $\Psi(E_i)$ and $\Psi$ respect direct sums.\cqfd
\vskip 12pt
\noi{\bf 1.7 Lemma:} {\sl The functor $\Psi$ is fully faithful.}
\vskip 12pt
\noi{\bf Proof:} Let $E$ and $E'$ be two $\B$-modules. The functor $\Psi$ induces a map $\Hom_\B(E,E')\fl\Hom_{\A\v}(\Psi(E),\Psi(E'))$ and we have to prove that this 
map is an isomorphism.
 
Because of the condition (cc6) there is a family $X_i$, $i\in I$ of $\A$-modules and an epimorphism $\alpha:\ \build\oplus_i^{} \Phi(X_i)\fl E$. Let $C_0$ be the direct
sum $\build\oplus_i^{} \Phi(X_i)$.
By doing the same thing with the kernel of $\alpha$, we get a family $Y_j$, $j\in J$ and an epimorphism $\beta:\ \build\oplus_j^{}\Phi(Y_j)\fl\Ker(\alpha)$. Thus we have 
an exact sequence in $\B$:
$$C_1\build\fl_{}^\beta C_0\build\fl_{}^\alpha E\fl 0$$
with $C_1=\ \build\oplus_j^{} \Phi(Y_j)$.

Because of the condition (cc2), every $\Phi(Y_j)$ is finitely presented in $\B$ and each morphism $\beta:\ \Phi(Y_j)\fl\ \build\oplus_i^{} \Phi(X_i)$ factors through a 
finite sum of the $\Phi(X_i)$ and because of the condition (cc1), the morphism $\beta$ is induced by $\A$-morphisms $\beta_{ij}:Y_j\fl X_i$ such that, for each $j\in J$, 
only finitely many of the $\beta_{ij}$ are non-zero.

Thus we have the following exact sequence in $\B$:
$$\build\oplus_j^{} \Phi(Y_j)\build\fl_{}^{\beta_{**}}\ \build\oplus_i^{} \Phi(X_i)\build\fl_{}^\alpha E\fl 0$$
and, because of Lemma 1.6, we get an exact sequence in $\A\v$:
$$\build\oplus_j^{} Y_j\v\build\fl_{}^{\beta_{**}}\ \build\oplus_i^{} X_i\v\fl\Psi(E)\fl 0$$

Thus we get the following commutative diagram with exact lines:
$$\diagram{0&\hflp{}{}&\Hom(E,E')&\hflp{\alpha^*}{}&\build\prod_i^{}\Hom(\Phi(X_i),E')&\hflp{\beta_{**}^*}{}&\build\prod_j^{}\Hom(\Phi(Y_j),E')\cr
&&\vfl{\lambda_0}{}&&\vfl{\lambda_1}{}&&\vfl{\lambda_2}{}\cr
0&\hflp{}{}&\Hom(\Psi(E),\Psi(E'))&\hflp{}{}&\build\prod_i^{}\Hom(X_i\v,\Psi(E'))&\hflp{\beta_{**}^*}{}&\build\prod_j^{}\Hom(Y_j\v,\Psi(E'))\cr}$$
where the vertical lines are induced by the functor $\Psi$. Hence $\lambda_1$ and $\lambda_2$ are isomorphisms and the map: $\Hom(E,E')\fl\Hom(\Psi(E),\Psi(E'))$ is an
isomorphism. Thus $\Psi$ is fully faithful.\cqfd
\vskip 12pt
\noi{\bf 1.8 Lemma:} {\sl The functor $\Psi$ is an equivalence of categories.}
\vskip 12pt
\noi{\bf Proof:} Let $M$ be an $\A\v$-module. As above, we find families $X_i$, $i\in I$ and $Y_j$, $j\in J$ of $\A$-modules, morphisms $\beta_{ij}:Y_j\fl X_i$ and
an exact sequence:
$$\build\oplus_j^{} Y_j\v\build\fl_{}^{\beta_{**}}\ \build\oplus_i^{} X_i\v\fl M\fl 0$$
such that, for each $j\in J$, only finitely may of the $\beta_{ij}$ are non-zero.

Let $E$ be the cokernel of the morphism:
$$\build\oplus_j^{} \Phi(Y_j)\build\fl_{}^{\beta_{**}}\ \build\oplus_i^{}\Phi(X_i)$$
Then it is easy to check that $M$ is isomorphic to $\Psi(E)$.\cqfd
\vskip 12pt
Because of this last lemma, $\B$ is a Grothendieck category and this fact completes the proof of Theorem 1.2.\cqfd
\vskip 12pt
\noi{\bf 1.9 Theorem:} {\sl Let $\A$ and $\B$ be two essentially small exact categories, $F:\A\fl\B$ be an exact functor and $\A\v$ and $\B\v$ be strict cocompletions of 
$\A$ and $\B$. Then $F$ has an extension $F\v:\A\v\fl\B\v$, unique up to isomorphism, which is a right-exact functor respecting (infinite) direct sums.}
\vskip 12pt 
\noi{\bf Proof:} Since $\A\v$ and $\B\v$ are Grothendieck categories, there are functorial $\Ker$, $\Coker$ and (infinite) direct sums in $\A\v$ and $\B\v$. 

Let $\Sigma\A$ be the following category:

An object in this category is a pair $(I,X_*)$, where $I$ is a set and $X_*$ is a family of $\A$-modules $X_i$, $i\in I$, and a morphism from $(I,X_*)$ to $(J,Y_*)$ is
a family of $\A$-morphisms $\alpha_{ji}:X_i\fl Y_j$ such that the sum $\sum_j\alpha_{ji}$ is finite for every $i\in I$.

The category $\Sigma\A$ is an additive category and the direct sum operation is a functor $S:\Sigma\A\fl\A\v$.
\vskip 12pt
\noi{\bf 1.10 Lemma:} {\sl The functor $S$ is fully faithful.}
\vskip 12pt
\noi{\bf Proof:} Let $(I,X_*)$ and $(J,Y_*)$ be two $\Sigma\A$-modules and $\alpha:S(I,X_*)\fl S(J,Y_*)$ be an $\A\v$-morphism. This morphism comes from a morphism
$\alpha':\build\oplus_{i\in I}^{}X_i\fl \build\oplus_{j\in J}^{}Y_j$ and, because each $X_*$ is finitely presented in $\A\v$, $\alpha'$ sends each $X_i$ to a finite
sum of the $Y_j$'s. Therefore we have:
$$\alpha'\in\build\prod_i^{}\build\oplus_j^{} \Hom(X_i,Y_j)$$
and $\alpha$ can be lifted in a unique way to a morphism $(I,X_*)\fl(J,Y_*)$. Hence, the functor $S$ is fully faithful.\cqfd
\vskip 12pt
For the same reason we have an additive category $\Sigma\B$ and a faithful functor $S:\Sigma\B\fl\B\v$. Moreover the functor $F:\A\fl\B$ induces a functor 
$\Sigma\A\fl\Sigma\B$ still denoted $F$. 

Let $\U$ be a skeleton of $\A$.

Let $E$ be a $\A\v$-module and $M(E)$ be the set of pairs $(X,f)$ with $X\in\U$ and $f\in\Hom(X,E)$. For each $u=(X,f)\in M(E)$ we set: $X_u=X$ and $f_u=f$. We set:
$L(E)=(M(E),X_*)\in\Sigma\A$. The family of morphisms $f_u$, $u\in M(E)$, is an $\A\v$-morphism $\varphi_E:S(L(E))\fl E$. Thus we have a functor $L:\A\v\fl \Sigma\A$ and
a morphism of functors $\varphi:SL\fl$ Id. Moreover $\varphi_E$ is an epimorphism for every $\A\v$-module $E$. Then for every $\A\v$-module $E$, we have an exact
sequence:
$$0\fl E'\fl SL(E)\build\fl_{}^{\varphi_E} E\fl 0$$
and, by setting: $L'(E)=L(E')$, we get an exact sequence:
$$SL'(E)\build\fl_{}^\beta SL(E)\build\fl_{}^{\varphi_E} E\fl 0$$
Since $S$ is fully faithful, $\beta$ comes from a $\Sigma\A$-morphism $\psi_E:L'(E)\fl L(E)$.

Thus we get two functors $L$ and $L'$ from $\A\v$ to $\Sigma\A$, two morphisms of functors: $\varphi:SL\fl$ Id and $\psi:L'\fl L$ and an exact sequence:
$$SL'\build\fl_{}^{S\psi}SL\build\fl_{}^\varphi \hbox{Id}\fl 0$$

The $\Sigma\A$-morphism $\psi:L'\fl L$ induces a $\Sigma\B$-morphism $F\psi:FL'\fl FL$. Let $\Phi:\A\v\fl\B\v$ be the cokernel of $SF\psi:SFL'\fl SFL$. Then we have an 
exact sequence:
$$SFL'\fl SFL\fl\Phi\fl0$$
and, for every $\A\v$-module $E$, we have an exact sequence in $\B\v$:
$$SFL'(E)\fl SFL(E)\fl\Phi(E)\fl0$$

Suppose that $F':\A\v\fl\B\v$ is an exact functor which respects direct sums and is an extension of $F:\A\fl \B$.

Since $F'$ is exact, the exact sequence $SL'\fl SL\fl \hbox{Id}\fl 0$ induces an exact sequence:$F'SL'\fl F'SL\fl F'\fl 0$ and, because $F'$ respects direct sums and is an 
extension of $F$ we have: $F'SL=SF'L=SFL$ and $F'SL'=SF'L'=SFL'$. Thus we get an exact sequence:
$$SFL'\fl SFL\fl F'\fl 0$$
and an isomorphism of functors $F'\fl\Phi$.

Thus the last thing to do is to show that $\Phi$ is right-exact, that it respects direct sums and that the restriction of $\Phi$ to $\A$ is isomorphic to $F$.
\vskip 12pt
\noi{\bf 1.11 Lemma:} {\sl We have an isomorphism of functors $\varepsilon:SF\build\fl_{}^\sim\Phi  S$.}
\vskip 12pt
\noi{\bf Proof:} Let $\A'$ be the fully exact subcategory of $\A$ generated by $\U$. This category is small and equivalent to $\A$.
 
Let $X=(I,X_*)$ be a $\Sigma\A'$-module and $E$ be the $\A\v$-module $SX=\sum_i X_i$. The inclusion $X_i\subset Y$ will be denoted $\varepsilon_i$. Let $J$ be the 
complement in $M(E)$ of the set $\{(X_i,\varepsilon_i)\}$. Then we have a decomposition: $L(E)=X\oplus Y$, where $Y=(J,Y_*)$ is a $\Sigma\A'$-module.

The morphism $\varphi_E:SL(E)\fl E$ restricts on a morphism $SY\fl SX$ coming from a $\Sigma\A'$-morphism $f:Y\fl X$ and we have a split exact sequence:
$$0\fl Y\build\fl_{}^{-f\oplus 1} X\oplus Y\build\fl_{}^{1\oplus f} X\fl 0$$
By applying the same method for $Y$, we find a $\Sigma\A'$-module $Z$ and a split exact sequence in $\Sigma\A$:
$$0\fl Z\fl Y\oplus Z\fl Y\fl 0$$
with: $L'(E)=L(SY)=SY\oplus SZ$. Moreover the composite morphism $SZ\oplus SY\fl SY\fl SX\oplus SY$ is the morphism $\psi_E$.

By applying the functor $F$ we get split exact sequences in $\Sigma\B$:
$$0\fl F(Y)\fl F(X)\oplus F(Y)\fl F(X)\fl 0$$
$$0\fl F(Z)\fl F(Y)\oplus F(Z)\fl F(Y)\fl 0$$
Thus we get an exact sequence:
$$0\fl F(Z)\fl F(Y)\oplus F(Z)\build\fl_{}^\beta F(X)\oplus F(Y)\build\fl_{}^\alpha F(X)\fl 0$$
such that: $FL(E)=F(X)\oplus F(Y)$, $FL'(E)=F(Y)\oplus F(Z)$, $F(\psi_E)=\beta$. Hence we get a functorial isomorphism $\varepsilon:SF(X)\build\fl_{}^\sim \Phi(E)=\Phi 
S(X)$.

Suppose now that $X=(I,X_*)$ is an $\Sigma\A$-module. Since $\U$ is a skeleton of $\A$, there are isomorphisms $u_i:X_i\fl Y_i$ with $Y_i$ in $\U$. These isomorphisms
induce an isomorphism $u:X\fl Y$ with $Y\in\Sigma\A'$ and we have a commutative diagram:
$$\diagram{SF(X)&\hfl{u}{\sim}&SF(Y)\cr\vfl{}{}&&\vfl{\varepsilon}{\sim}\cr \Phi S(X)&\hfl{u}{\sim}&\Phi S(Y)\cr}$$
where the morphism $\varepsilon_X:SF(X)\fl \Phi S(X)$ is the composite $u^{-1}\varepsilon u$. It is easy to see that this isomorphism doesn't depend on the choice of $u$
and we get a well defined isomorphism $\varepsilon: SF\build\fl_{}^\sim \Phi S$.\cqfd
\vskip 12pt 
For every $X\in\A\v$ and every $\A\v$-morphism $f:X\fl Y$, we set:
$$F\v(X)=\left\{\matrix{F(X)&\ \hbox{if}\ X\in\A\cr \Phi(X)&\hbox{otherwise}\cr}\right.$$
$$F\v(f)=\left\{\matrix{F(f):F(X)\fl F(Y)&\ \hbox{if}\ X,Y\in\A\cr \varepsilon F(f):F(X)\fl\Phi(X)&\ \hbox{if}\ X\in\A, Y\not\in\A\cr
\Phi(f)\varepsilon^{-1}:\Phi(X)\fl F(Y)&\ \hbox{if}\ X\not\in\A, Y\in\A\cr \Phi(f):\Phi(X)\fl\Phi(Y)&\ \hbox{otherwise}\cr}\right.$$
Then it is easy to prove the following:
\vskip 12pt
\noi{\bf 1.12 Lemma:} {\sl $F\v$ is a functor from $\A\v$ to $\B\v$ which is equal to $F$ on $\A$. Moreover we have the equality $SF=F\v S$ and an exact sequence:
$$SFL'\fl SFL\fl F\v\fl0$$}
\vskip 12pt

\noi{\bf 1.13 Lemma:} {\sl The functor $F\v$ respects epimorphisms.} 
\vskip 12pt 
\noi{\bf Proof:} Let $f:E\fl E'$ be an $\A\v$-epimorphism. We have: $L(E')=(I,X_*)$, with $I=M(E')$ and $i=(X_i,g_i)$, for each $i\in I$.  and an epimorphism:
$$\build\oplus_{i\in I}^{} F(X_i)\build\fl_{}^{\oplus f_i} F\v(E')$$
Because of property (cc4), for each $i\in I$, there is an $\A$-module $Y_i$ and a commutative diagram:
$$\diagram{Y_i&\hfl{h_i}{}&E\cr\vfl{\alpha_i}{}&&\vfl{f}{}\cr X_i&\hfl{g_i}{}& E'\cr}$$
where $\alpha_i$ is a deflation.

We have two $\Sigma\A$-modules: $X=L(E')=(I,X_*)$ and $Y=(I,Y_*)$. Since morphisms $\alpha_i:Y_i\fl X_i$ are deflations, each morphism $F(\alpha_i):F(Y_i)\fl F(X_i)$
is an epimorphism and $SF(\alpha)=F\v S(\alpha)$ is also an epimorphism. Therefore in the following commutative diagram:
$$\diagram{F\v SL(Y)&\hfl{}{}&F\v(E)\cr\vfl{F\v SL(\alpha)}{}&&\vfl{}{F\v(f)}\cr F\v SL(X)&\hfl{F\v(\varphi_{E'})}{}&F\v(E')\cr}$$
both morphisms $F\v SL(\alpha)$ and $F\v(\varphi_{E'})$ are epimorphisms. Hence $F\v(f)$ is also an epimorphism.\cqfd
\vskip 12pt
\noi{\bf 1.14 Lemma:} {\sl Let $X\fl Y\fl Z$ be a sequence in $\Sigma\A$ such that: 
$$SX\fl SY\fl SZ\fl 0$$ 
is exact in $\A\v$. Then
$$SF(X)\fl SF(Y)\fl SF(Z)\fl 0$$
is exact in $\B\v$.}
\vskip 12pt 
\noi{\bf Proof:} Suppose $Z=(I,Z_*)$.  Since $SY\fl SZ$ is an epimorphism, there are $\A$-conflations $A_i\build\fl_{}^{\alpha_i} B_i\build\fl_{}^{\beta_i} Z_i$ and 
commutative diagrams:
$$\diagram{A_i&\hfl{\alpha_i}{}&B_i&\hfl{\beta_i}{}&Z_i\cr\vfl{}{}&&\vfl{}{}&&\vfl{}{}\cr K&\hfl{}{}&SY&\hfl{}{}&SZ\cr}$$
where $K$ is the kernel of $SY\fl SZ$.

Moreover, because $SX\fl K$ is an epimorphism, there are deflations $\gamma_i:C_i\fl A_i$ and commutative diagrams:
$$\diagram{C_i&\hfl{\gamma_i}{}&A_i\cr\vfl{}{}&&\vfl{}{}\cr SX&\hfl{}{}&K\cr}$$
Then we get a commutative diagram in $\Sigma\A$:
$$\diagram{C&\hfl{\alpha\gamma}{}&B&\hfl{\beta}{}&Z\cr\vfl{}{}&&\vfl{}{}&&\vfl{=}{}\cr X&\hfl{}{}&Y&\hfl{}{}&Z\cr}$$
with: $A=(I,A_*)$, $B=(I,B_*)$ and $C=(I,C_*)$.

Since $F$ is exact each sequence $F(A_i)\fl F(B_i)\fl F(Z_i)$ is a $\B$-conflation and each $F(C_i)\fl F(A_i)$ is a deflation. Then the sequence:
$$F\v SC\fl F\v SB\fl F\v SZ\fl 0$$
is exact.

On the other hand, we have an epimorphism $SX\oplus SB\fl SY$ and then an epimorphism: $F\v SX\oplus F\v SB\fl F\v SY$.

We have a commutative diagram:
$$\diagram{F\v SC&\hfl{\alpha\gamma}{}&F\v SB&\hfl{\beta}{}&F\v SZ&\hfl{}{}&0\cr\vfl{}{}&&\vfl{}{}&&\vfl{=}{}&&\cr F\v SX&\hfl{}{}&F\v SY&\hfl{}{}&F\v SZ&\hfl{}{}&0\cr}$$
where the top line is exact, the composite morphism $F\v SX\fl F\v SZ$ is null and the morphism $F\v SX\oplus F\v SB\fl F\v SY$ is onto. Therefore the bottom line is exact
and the lemma is proven.\cqfd
\vskip 12pt
\noi{\bf 1.15 Lemma:} {\sl The functor $F\v$ is right-exact.}
\vskip 12pt  
\noi{\bf Proof:} Let $0\fl E\fl E'\fl E''\fl0$ be a short exact sequence in $\A\v$. We have:
$$M(E)\subset M(E')\ \ \ \hbox{and}\ \ \ M(\Ker(\varphi_E))\subset M(\Ker(\varphi_{E'}))$$

Let $I$ and $J$ be the complements of $M(E)$ and $M(\Ker(\varphi_E))$ in $M(E')$ and $M(\Ker(\varphi_{E"}))$. Thus we can set:
$$L(E)=X\hskip 24pt L'(E)=X\oplus X'\hskip 24pt L(E')=Y\hskip 24pt L'(E')=Y\oplus Y'$$
with: $X'=(I,X'_*)$ and $Y'=(J,Y'_*)$ and we get a commutative diagram in $\A\v$ with exact lines and exact columns:
$$\diagram{0&\hfl{}{}&SY&\hfl{}{}&SY\oplus SY'&\hfl{}{}&SY'&\hfl{}{}&0\cr&&\vfl{}{}&&\vfl{}{}&&\vfl{\alpha}{}&&\cr
0&\hfl{}{}&SX&\hfl{}{}&SX\oplus SX'&\hfl{}{}&SX'&\hfl{}{}&0\cr&&\vfl{}{}&&\vfl{}{}&&\vfl{\beta}{}&&\cr
0&\hfl{}{}&E&\hfl{}{}&E'&\hfl{}{}&E''&\hfl{}{}&0\cr&&\vfl{}{}&&\vfl{}{}&&\vfl{}{}&&\cr&&0&&0&&0&&\cr}\leqno{(D)}$$
In particular the sequence:
$$SY'\build\fl_{}^\alpha SX'\build\fl_{}^\beta E''\fl 0$$
is exact.
 
In the diagram $F\v(D)$, the top two lines and the two left columns are exact. Therefore we get an exact sequence:
$$F\v(E)\fl F\v(E')\fl \Coker(F\v(\alpha))\fl 0$$ 

For each $i\in I$ the morphism $X'_i\fl E''$ is isomorphic to an element of $M(E'')$ and we get morphisms $\gamma_i:X'_i\fl SL(E'')$. These morphisms induce a morphism
$\gamma:X'\fl L(E'')$ and the following diagram is commutative:
$$\diagram{SX'&\hfl{\gamma}{}&SL(E'')\cr\vfl{\beta}{}&&\vfl{}{\varphi_{E''}}\cr E''&\hfl{=}{}&E''\cr}$$

Let $K$ be the kernel of $\varphi_{E''}:SL(E'')\fl E''$. Since $SL'(E'')\fl K$ is onto, for each $j\in J$, there is a deflation $\delta_j;Z_j\fl Y'_j$ and a commutative 
diagram:
$$\diagram{Z_j&\hfl{}{}&SL'(E'')\cr\vfl{\delta_j}{}&&\vfl{}{}\cr Y'_j&\hfl{\gamma\alpha}{}&K\cr}$$
Thus we have a $\Sigma\A$-module $Z$, a morphism $\delta:Z\fl Y$ and a commutative diagram with exact lines:
$$\diagram{SZ&\hfl{\alpha\delta}{}&SX'&\hfl{\beta}{}&E''&\hfl{}{}&0\cr \vfl{}{}&&\vfl{}{}&&\vfl{}{=}&&\cr SL'(E'')&\hfl{}{}&SL(E'')&\hfl{}{}&E''&\hfl{}{}&0\cr}$$
This exact sequence implies an exact sequence in $\A\v$:
$$SZ\fl SX'\oplus SL'(E'')\fl SL(E'')\fl 0$$
and, because of Lemma 1.4, the sequence:
$$SFZ\fl SFX'\oplus SFL'(E'')\fl SFL(E'')\fl 0$$
is exact.

On the other hand the sequence:
$$SFL'(E'')\fl SFL(E'')\fl F\v(E'')\fl 0$$
is exact. Therefore in the following diagram:
$$\diagram{SFZ&\hfl{\alpha\delta}{}&SFX'&\hfl{\beta}{}&F\v(E'')&\hfl{}{}&0\cr \vfl{}{}&&\vfl{}{}&&\vfl{}{=}&&\cr SFL'(E'')&\hfl{}{}&SFL(E'')&\hfl{}{}&F\v(E'')&\hfl{}{}&0
\cr}$$
the two lines are exact. In particular the sequence:
$$SFZ\fl SFX'\fl F\v(E'')\fl 0$$
is exact and, because $Z\fl Y'$ is induced by deflations, $SFZ\fl SFY'$ is an epimorphism. Hence the sequence:
$$SFY'\build\fl_{}^{F\v(\alpha)} SFX'\fl F\v(E'')\fl 0$$
is exact and and we have:
$$F\v(E'')\simeq \Coker(F\v(\alpha))\simeq \Coker(F\v(E)\fl F\v(E'))$$
Hence $F\v$ is right-exact and the lemma is proven.\cqfd
\vskip 12pt
\noi{\bf 1.16 Lemma:} {\sl The functor $F\v$ respects direct sums.}
\vskip 12pt
\noi{\bf Proof:} Let $E_i$, $i\in I$, be a family of $\A\v$-modules. For each $i\in I$, there are two $\Sigma\A$-modules $X_i$ and $Y_i$ and an exact sequence in $\A\v$:
$$SY_i\fl SX_i\fl E_i\fl 0$$
Let $E$ the direct sum of the $E_i$. We have an exact sequence:
$$\build\oplus_{i\in I}^{} SY_i\fl \build\oplus_{i\in I}^{} SX_i\fl E\fl 0$$ 

Let $X$ be the direct sum of the $X_i$ and $Y$ be the direct sum of the $Y_i$. Then we have:
$$\build\oplus_{i\in I}^{} SX_i=SX\hskip 48pt \build\oplus_{i\in I}^{} SY_i=SY$$
and we get an exact sequence:
$$SY\fl SX\fl E$$ 
Since $F\v$ is exact, we have an exact sequence:
$$F\v SY\fl F\v SX\fl F\v(E)\fl 0$$ 
Therefore the following sequence is exact:
$$SFY\fl SFX\fl F\v(E)\fl 0$$
and we have:
$$F\v(E)=\Coker(SFY\fl SFX)=\Coker(\build\oplus_{i\in I}^{} SFY_i\fl\build\oplus_{i\in I}^{} SFX_i)$$
$$=\build\oplus_{i\in I}^{} \Coker(SFY_i\fl SFX_i)=\build\oplus_{i\in I}^{} F\v(E_i)$$ 

Hence $F\v$ respects direct sums and the lemma is proven as well as Theorem 1.9 and Theorem 4.\cqfd
\vskip 12pt
\noi{\bf Examples:} 
\vskip 12pt
Let $A$ be a ring, $\Mod_A$ be the category of right $A$-modules and $\P_A$ be the category of finitely generated projective modules in $\Mod_A$. Then $\P_A$ is an
essentially small exact category and $\Mod_A$ is a strict cocompletion of $\P_A$.
\vskip 12pt 
 
Let $A$ be a ring and $S$ be an $A$-bimodule which is flat on the left.

Let $M$ be a module in $\Mod_A$ and $\theta:M\fl M\otimes S$ be a morphism in $\Mod_A$. By iteration of $\theta$, we get morphisms $\theta^n:M\fl M\otimes S^{\otimes n}$.
 
We said that $\theta$ is nilpotent if every $x\in M$ is killed by some power $\theta^n$ of $\theta$.

Let $\Nil(A,S)$ (resp. $\Nil(A,S)\v$) be the category of pairs $(M,\theta)$ where $M$ is in $\P_A$ (resp. in $\Mod_A$) and $\theta:M\fl M\otimes S$ is a nilpotent 
morphism. Then $\Nil(A,S)$ is an essentially small exact category and $\Nil(A,S)\v$ is a strict cocompletion of $\Nil(A,S)$.
\vskip 12pt
Consider two rings $A$ and $B$ and an $(A,B)$-bimodule $S$ which is finitely generated projective on the right. We have a functor $F:\P_A\fl \P_B$ sending each 
$\P_A$-module $M$ to $M\otimes S$. This functor is exact and has an extension $F\v:\Mod_A\fl\Mod_B$ sending each module $M$ in $\Mod_A$ to the module $M\otimes S$. Then 
$F\v$ is exact if and only if $S$ is flat on the left.
\vskip 12pt
Let $\A$ be an essentially small exact category and $\A\v$ be a strict cocompletion of $\A$. An $\A$-complex is a $\Z$-graded differential $\A$-module $C$ with a 
differential of degree $-1$. An $\A$-complex $C$ has homology modules $H_*(C)$ in $\A\v$. A finite $\A$-complex is an $\A$-complex $C$ such that $C_i=0$ for $|i|$ large
enough. The category of finite $\A$-complexes is an essentially small exact category and, equipped with the inflations and the homology equivalences, it is also a
Waldhausen category denoted $\A_*$ (see [W2]).
\vskip 12pt
\noi{\bf 1.17 Theorem:} {\sl Let $\A$ be an essentially small exact category and $\A_*$ be the Waldhausen category of finite $\A$-complexes. Then the inclusion 
$\A\subset\A_*$ induces a homotopy equivalence: $K(\A)\build\fl_{}^\sim K(\A_*)$.} 
\vskip 12pt 

\noi{\bf Proof:} Let $\A\v$ be a strict cocompletion of $\A$ and $\B$ be the category of $\A\v$-modules $X$ such that there is some $\A$-module $Y$ with $X\oplus Y\in\A$. 
The category $\B$ is exact and the inclusions $\A\subset \B$ and $\B\subset\A\v$ are fully exact. Moreover $\A\v$ is a strict cocompletion of $\B$ and $\A$ is strictly 
cofinal in $\B$.

Then the inclusion $\A\subset\B$ induces a homotopy equivalence in K-theory: $K(\A)\build\fl_{}^\sim K(\B)$. 

Let $\B_*$ be the Waldhausen category of finite $\B$-complexes. 

Let $X$ and $Y$ be two $\B$-modules and $f:X\fl Y$ be a morphism which is an epimorphism in $\A\v$. Let $K$ be the kernel of $f$. 

Since $X$ and $Y$ are in $\B$, there are two $\A$-modules $A$ and $B$ such that $X\oplus A\in\A$ and $Y\oplus B\in\A$. Then $X\oplus A\oplus B$ and $Y\oplus A\oplus B$
are $\A$-modules and we have:
$$K=\Ker(X\build\fl_{}^f Y)=\Ker(X\oplus A\oplus B\build\fl_{}^{f\oplus 1\oplus 1}Y\oplus A\oplus B)$$ 
and, because of the property (cc8), there is an $\A$-module $C$ such that $K\oplus C\in\A$. Then $K$ is in $\B$ and $\B$ is stable under kernel of epimorphisms.

Hence the Gillet-Waldhausen theorem implies that $K(\B)\fl K(\B_*)$ is a homotopy equivalence.

On the other hand the inclusions $\A\subset\B$, $\A_*\subset\B_*$, $\A\subset\A_*$ and $\B\subset\B_*$ imply a commutative diagram:
$$\diagram{K(\A)&\hfl{\alpha}{}&K(\A_*)\cr\vfl{\lambda}{}&&\vfl{\lambda_*}{}\cr K(\B)&\hfl{\beta}{}&K(\B_*)\cr}$$
where $\beta$ is a homotopy equivalence.

Because $\A$ is strictly cofinal in $\B$, $\A_*$ is also strictly cofinal in $\B_*$ and $\lambda$ and $\lambda_*$ are homotopy equivalences. Therefore $\alpha$ is a
homotopy equivalence.\cqfd 
\vskip 24pt
\noi{\bf 2. Regularity.}
\vskip 12pt
Let $\A$ be an essentially small exact category and $\A\v$ be a strict cocompletion of $\A$. Let $\overline\A_*$ be the category of $\A$-complexes and 
$\A\v_*$ be the category of $\A\v$-complexes. By considering each $\A$-module as an $\A$-complex concentrated in degree $0$, we may consider $\A$ as a
subcategory of $\overline\A_*$ and also $\A\v$ as a subcategory of $\A\v_*$.

Notice that the homology of an $\A$-complex is a well-defined graded $\A\v$-module. 

These categories are exact: a conflation in one of these categories of complexes is a sequence which is a conflation in each degree.
\vskip 12pt
\noi{\bf 2.1 Theorem:} {\sl Let $\A$ be an essentially small exact category and $\A\v$ be a strict cocompletion of $\A$. Let $C$ be an $\A$-complex. Suppose $\A$ is regular
and $H_i(C)=0$ (in $\A\v$) for any large enough integer $i$. Then there is an integer $n$ such that, for any $i\geq n$, the kernel of $d:C_i\fl C_{i-1}$ in $\A\v$ is a 
direct summand of an $\A$-module.}
\vskip 12pt
\noi{\bf Proof:} Let $\C$ be the class of $\A\v$-modules $X$ such that, for each $\A$-complex $C$, every morphism from $C$ to $X$ factors through a finite $\A$-complex.
\vskip 12pt
\noi{\bf 2.2 Lemma:} {\sl The class $\C$ is cocomplete and contains $\A$.}
\vskip 12pt
\noi{\bf Proof:} Let $X_i$, $i\in I$ be a filtered system of $\A\v$-modules in $\C$ and $X$ be its colimit. Let $C$ be an $\A$-complex and $f:C\fl X$ be a morphism of
$\A\v$-complexes. This morphism is given by a morphism $g:K\fl X$ where $K$ is the cokernel in $\A\v$ of the morphism $d:C_1\fl C_0$. Since $K$ is finitely presented, the 
morphism $g$ factors through some $X_i$. Then, because $X_i$ is in $\C$, $f$ factors also through a finite complex and $X$ is in the class $\C$.

Let $X$ be an $\A$-module and $f:C\fl X$ be a morphism of $\A\v$-complexes. The morphism $f$ factors obviously through the finite $\A$-complex $X$ and $X$ belongs to 
$\C$.\cqfd
\vskip 12pt
\noi{\bf 2.3 The $k$-cone of a morphism:}
\vskip 12pt
If $C$ is a $\A\v$-complex, its underlying graded $\A\v$-module will be denoted $C^\circ$.

Let $C$ and $C'$ be two $\A\v$-complexes. For any integer $k$, denote $\Map(C,C')_k$ the set of $\A$-morphisms $C^\circ\fl C'^\circ$ of degree $k$. Thus we have:
$$\Map(C,C')_k\ \simeq\ \build\prod_{i\in\Z}^{}\Hom(C_i,C'_{i+k})$$
The graded module $\Map(C,C')$ has a differential defined by:
$$\forall f\in\Map(C,C')_k,\ \ d(f)=d\circ f-(-1)^k f\circ d$$
and $\Map(C,C')$ is a graded differential $\Z$-module. If $k$ is an integer, the elements of $\Map(C,C')_k$ are called $k$-maps or maps of degree $k$. A $k$-map $f$ with 
$d(f)=0$ is called a $k$-morphism (or a morphism of degree $k$) and a morphism in the image of $d:\Map(C,C')_{k+1}\fl\Map(C,C')_k$ is called a $k$-homotopy (or a homotopy 
of degree $k$. Thus a morphism between two $\A$-complexes is a $0$-morphism.

Let $a$ be an integer and $\varphi:C\fl C'$ be an $a$-morphism of $\A\v$-complexes. Then we have morphisms of complexes $i:C'\fl C\oplus C'$, $s:C\fl C\oplus C'$, 
$p:C\oplus C'\fl C$ and $r:C\oplus C'\fl C'$. The morphism $r$ is a retraction of the inclusion $i$ and $s$ is a section of the projection $p$. Thus we have the relations:
$$pi=0\hskip 24pt rs=0\hskip 24pt ri=1\hskip 24pt ps=1\hskip 24pt ir+sp=1$$

If $k$ is an integer, there is a unique way to modify the degree and the differential on $C\oplus C'$ in such a way that:
$$\partial^\circ i=k\hskip 24pt \partial^\circ r=-k\hskip 24pt \partial^\circ p=-1-k-a\hskip 24pt \partial^\circ s=k+a+1$$
$$d(i)=0\hskip 24pt d(p)=0\hskip 24pt d(s)=i\varphi\hskip 24pt d(r)=-(-1)^k \varphi p$$
This modified complex $C''$ will be called the $k$-cone of $\varphi$ and denoted $C_k(\varphi)$. Thus we have an exact sequence:
$$0\fl C'\build\fl_{}^i C''\build\fl_{}^p C\fl 0$$
where $i$ is a $k$-morphism, $p$ is a $(-1-k-a)$-morphism, $r$ is a retraction of $i$ and $s$ is a section of $p$. 
\vskip 12pt
\noi{\bf 2.4 Lemma:} {\sl Consider a commutative diagram of $\A\v$-complexes:
$$\diagram{C'&\hfl{h}{}&F\cr\vfl{f}{}&&\vfl{}{\varphi}\cr C&\hfl{g}{}&X\cr}$$
where $C$ and $C'$ are $\A$-complexes, $F$ is a finite $\A$-complex, $X$ is a $\A\v$-module and $f$ is a deflation with finite kernel. Then $g$ factors through 
a finite complex.}
\vskip 12pt
\noi{\bf Proof:} Let $K$ be the kernel of $f$. Since $K$ and $F$ are finite $\A$-complexes, there are integers $a<b$ such that: $K_i=F_i=0$ for any $i<a$ and any $i>b$. 
Let $H$ be the graded $\A$-module such that: $H_i=C'_i$ if $a\leq i\leq b$ and $H_i=0$ otherwise. Equipped with the differential $d=0$, $H$ is a finite $\A$-complex. The 
identities $C'_i\fl H_i$ induce a $0$-map $\alpha:C'\fl H$.

Let $\Sigma$ be the $0$-cone of the identity of $H$. We have an exact sequence:
$$0\fl H\build\fl_{}^i \Sigma\build\fl_{}^p H\fl 0$$
a retraction $r:\Sigma\fl H$ of $i$, a section $s:H\fl \Sigma$ of $p$ and:
$$\partial^\circ i=\partial^\circ r=0\hskip 24pt\partial^\circ p=-1\hskip 24pt\partial^\circ s=1$$
$$d(i)=0\hskip 24pt d(p)=0\hskip 24pt d(s)=i\hskip 24pt d(r)=-p$$

We have a morphism $\beta=d(s\alpha)=i\alpha-sd(\alpha):C'\fl \Sigma$ and then a commutative diagram:
$$\diagram{C'&\hfl{h\oplus\beta}{}&F\oplus \Sigma\cr\vfl{f}{}&&\vfl{}{\varphi\oplus0}\cr C&\hfl{g}{}&X\cr}$$
But $h\oplus\beta:C'\fl F\oplus \Sigma$ restricts to an inflation $K\fl F\oplus \Sigma$ with cokernel $F'$. Then $F'$ is a finite $\A$-complex and $g$ factors through it.
\cqfd
\vskip 12pt
\noi{\bf 2.5 Lemma:} {\sl Let $g:X\fl Y$ be an epimorphism in $\A\v$, $F$ be a finite $\A$-complex and $\varphi:F\fl Y$ be a morphism of $\A\v$-complexes. Then there is
a finite $\A$-complex $C$ and a commutative diagram:
$$\diagram{C&\hfl{\alpha}{}&F\cr\vfl{\beta}{}&&\vfl{}{\varphi}\cr X&\hfl{g}{}&Y\cr}$$
where $\alpha$ is a deflation having a section in each non-zero degree.}
\vskip 12pt
\noi{\bf Proof:} Since $X\fl Y$ is an epimorphism there is a deflation $a_0:K_0\fl F_0$ and a commutative diagram:
$$\diagram{K_0&\hfl{a_0}{}&F_0\cr\vfl{b_0}{}&&\vfl{}{\varphi_0}\cr X&\hfl{g}{}&Y\cr}$$
Set: $K_i=F_i$ for each $i\not=0$. Thus we get a graded $\A$-module $K$ and, equipped with the differential $d=0$, $K$ is a finite $\A$-complex. Moreover the identities 
$K_i=F_i$ and the maps $a_0:K_0\fl F_0$ and $0:F_0\fl K_0$ induce two $0$-maps $a:K\fl F$ and $c:F\fl K$. The morphism $b_0:K_0\fl E$ defines also a $0$-map $b:K\fl X$. 

Let $C$ be the $-1$-cone of the identity of $K$. We have an exact sequence:
$$0\fl K\build\fl_{}^i C\build\fl_{}^p K\fl 0$$
a retraction $r:C\fl K$ of $i$, a section $s:K\fl C$ of $p$ and:
$$\partial^\circ i=-1\hskip 24pt\partial^\circ r=1\hskip 24pt\partial^\circ p=\partial^\circ s=0$$
$$d(i)=0\hskip 24pt d(p)=0\hskip 24pt d(s)=i\hskip 24pt d(r)=p$$
Then we have a commutative diagram:
$$\diagram{C&\hfl{d(ar)}{}&F\cr\vfl{d(br)}{}&&\vfl{}{\varphi}\cr X&\hfl{g}{}&Y\cr}$$

Set: $\alpha=d(ar)=d(a)r+ap$. We have: $\alpha s=a$ and $\alpha s$ is a deflation in each degree. Therefore $\alpha$ is a deflation in each degree and $\alpha$ is a 
deflation. 

We have also: $\alpha sc=(d(a)r+ap)sc=apsc=ac$. Then $sc$ is a section of $\alpha$ is each non-zero degree.\cqfd
\vskip 12pt
Let $0\fl X\fl Y\fl Z\fl 0$ be an exact sequence of $\A\v$-modules.
\vskip 12pt
\noi{\bf 2.6 Lemma:} {\sl Suppose $X$ and $Y$ are in $\C$. Then $Z$ is also in $\C$.}
\vskip 12pt
\noi{\bf Proof:} Let $C$ be an $\A$-complex and $f:C\fl Z$ be a morphism of complexes. Because of property (cc4), there is an $\A$-deflation $C'_0\fl C_0$ and a 
commutative diagram
$$\diagram{C'_0&\hfl{}{}&C_0\cr\vfl{h}{}&&\vfl{f_0}{}\cr Y&\hfl{v}{}&Z\cr}$$
And, by base change, we get a cartesian square:
$$\diagram{C'_1&\hfl{}{}&C_1\cr\vfl{}{}&&\vfl{d}{}\cr C'_0&\hfl{}{}&C_0\cr}$$

Then we get an $\A$-complex $C'=(\dots\build\fl_{}^d C_2\build\fl_{}^d C'_1\fl C'_0\fl C_{-1}\build\fl_{}^d \dots)$ and a deflation $g:C'\fl C$ with finite kernel.

The morphism $h:C'_0\fl Y$ can be seen as a $0$-map from $C'$ to $Y$ and there is a unique $(-1)$-morphism $\lambda:C'\fl X$ such that: $d(h)=u\lambda$.

The $(-1)$-morphism $\lambda$ can be seen as a morphism from some suspension of $C'$ to $X$ and, because $X$ is in the class $\C$, $\lambda$ factors through a finite 
$\A$-complex $F$ via a morphism $\alpha:F\fl X$ and a $(-1)$-morphism $\mu:C'\fl F$.

Let $\Sigma$ be the $0$-cone of $\mu:C'\fl F$. We have an exact sequence:
$$0\fl F\build\fl_{}^i \Sigma\build\fl_{}^p C'\fl 0$$
a retraction $r:\Sigma\fl F$ of $i$, a section $s:C'\fl \Sigma$ of $p$ and:
$$\partial^\circ i=\partial^\circ r=\partial^\circ p=\partial^\circ s=0$$
$$d(i)=0\hskip 24pt d(p)=0\hskip 24pt d(s)=i\mu\hskip 24pt d(r)=-\mu p$$

Then we get a commutative diagram in $\A\v_*$ with exact lines:
$$\diagram{0&\hfl{}{}&F&\hfl{i}{}&\Sigma&\hfl{p}{}&C'&\hfl{}{}&0\cr&&\vfl{\alpha}{}&&\vfl{f'}{}&&\vfl{fg}{}&&\cr 0&\hfl{}{}&E&\hfl{u}{}&E'&\hfl{v}{}&E''&\hfl{}{}&0\cr}$$
with: $f'=hp+u\alpha r$.

But $Y$ is in $\C$ and $f':\Sigma\fl Y$ factors through a finite $\A$-complex $F'$ via a morphism $\varphi:\Sigma\fl F'$. 

We have a commutative diagram of morphisms:
$$\diagram{(F\build\fl_{}^i \Sigma)&\hfl{}{}&(F\build\fl_{}^{\varphi i} F')\cr\vfl{}{}&&\vfl{}{}\cr(0\fl C')&\hfl{}{}&(0\fl Z)\cr}$$
and, by applying the $0$-cone functor $C_0()$, we get a commutative diagram of complexes:
$$\diagram{C_0(F\build\fl_{}^i \Sigma)&\hfl{}{}&C_0(F\build\fl_{}^{\varphi i} F')\cr\vfl{}{}&&\vfl{}{}\cr C_0(0\fl C')&\hfl{}{}&C_0(0\fl Z)\cr}$$
But we have: $C_0(0\fl C')=C'$, $C_0(0\fl Z)=Z$ and $C_0(F\build\fl_{}^i \Sigma)\fl C_0(0\fl C')$ is a deflation with kernel $C_0(F\build\fl_{}^{\rm Id}F)$. 
This kernel is a finite acyclic $\A$-complex and we have commutative diagrams:
$$\diagram{C''&\hfl{}{}&F''\cr\vfl{g'}{}&&\vfl{}{}\cr C'&\hfl{fg}{}&Z\cr}\hskip 48pt \diagram{C''&\hfl{}{}&F''\cr\vfl{gg'}{}&&\vfl{}{}\cr C&\hfl{f}{}&Z\cr}$$
where $C''$ is the complex $C_0(F\build\fl_{}^i \Sigma)$, $F''$ is the finite $\A$-complex $C_0(F\build\fl_{}^{\varphi i} F')$ and $g':C''\fl C'$ is a deflation
with a finite acyclic kernel.

Since $g$ and $g'$ have finite kernels, the kernel of $gg'$ is a finite $\A$-complex. Then, because of Lemma 2.4, $f:C\fl Z$ factors through a finite complex.\cqfd 
\vskip 12pt
\noi{\bf 2.7 Lemma:} {\sl Suppose $X$ and $Z$ are in $\C$. Then $Y$ is also in $\C$.}
\vskip 12pt
\noi{\bf Proof:} Let $C$ be an $\A$-complex and $f:C\fl Y$ be a morphism of $\A\v$-complexes. Because $Z$ is in $\C$ there is a commutative diagram:
$$\diagram{C&\hfl{h}{}&F\cr\vfl{f'}{}&&\vfl{}{}\cr Y&\hfl{v}{}&Z\cr}$$
where $F$ is a finite $\A$-complex.

Because of Lemma 2.5, there is a finite $\A$-complex $K$ and a deflation $\alpha:K\fl F$ such that the composite $K\build\fl_{}^\alpha F\fl Z$ factors
through $Y$. Let $C'$ be the kernel of the deflation $C\oplus K\build\fl_{}^{h\oplus\alpha} F$. We have a commutative diagram with exact lines:
$$\diagram{0&\hfl{}{}&C'&\hfl{}{}&C\oplus K&\hfl{}{}&F&\hfl{}{}&0\cr&&\vfl{}{}&&\vfl{}{}&&\vfl{}{}&&\cr 0&\hfl{}{}&X&\hfl{u}{}&Y&\hfl{v}{}&Z&\hfl{}{}&0\cr}$$

Since $X$ is in $\C$, the morphism $C'\fl X$ factors through a finite complex $F'$ and, by cobase change, we get the following commutative diagram with exact lines:
$$\diagram{0&\hfl{}{}&C'&\hfl{}{}&C\oplus K&\hfl{}{}&F&\hfl{}{}&0\cr&&\vfl{}{}&&\vfl{}{}&&\vfl{}{}&&\cr 
0&\hfl{}{}&F'&\hfl{}{}&F''&\hfl{}{}&F&\hfl{}{}&0\cr&&\vfl{}{}&&\vfl{}{}&&\vfl{}{}&&\cr
0&\hfl{}{}&X&\hfl{u}{}&Y&\hfl{v}{}&Z&\hfl{}{}&0\cr}$$

Since $F$ and $F'$ are finite $\A$-complexes, $F''$ is also a finite $\A$-complex and the morphism $C\oplus K\fl Y$ factors through a finite complex. The result follows.
\cqfd
\vskip 12pt
\noi{\bf 2.8 Lemma:} {\sl Suppose $Y$ and $Z$ are in $\C$. Then $X$ is also in $\C$.}
\vskip 12pt
\noi{\bf Proof:} Let $C$ be an $\A$-complex and $f:C\fl X$ be a morphism of complexes. Let $C'$ be the quotient of $C$ by its $(-1)$-skeleton. The morphism $f$ factors
through $C'$ by a morphism $g:C'\fl X$. Since $Y$ is in $\C$ there is a finite $\A$-complex $F'$ and a commutative diagram:
$$\diagram{C'&\hfl{h}{}&F'\cr\vfl{g}{}&&\vfl{\varphi}{}\cr X&\hfl{u}{}&Y\cr}$$

Let $\Sigma$ be the $0$-cone of $h$. We have an exact sequence:
$$0\fl F'\build\fl_{}^i \Sigma\build\fl_{}^p C'\fl 0$$
a retraction $r:\Sigma\fl F'$ of $i$, a section $s:C'\fl \Sigma$ of $p$ and:
$$\partial^\circ i=\partial^\circ r=0\hskip 24pt\partial^\circ p=-1\hskip 24pt\partial^\circ s=1$$
$$d(i)=d(p)=0\hskip 24pt d(s)=ih\hskip 24pt d(r)=-hp$$

We have a commutative diagram:
$$\diagram{C'&\hfl{h}{}&F'&\hfl{i}{}&\Sigma\cr\vfl{g}{}&&\vfl{\varphi}{}&&\vfl{\psi}{}\cr X&\hfl{u}{}&Y&\hfl{v}{}&Z\cr}$$
where: $\psi=v\varphi r$ is a morphism. Moreover $ih=d(s)$ is a homotopy.

Since $Z$ is in $\C$ the morphism $\psi$ factors through a finite $\A$-complex $F''$ and we have a commutative diagram of complexes:
$$\diagram{C'&\hfl{h}{}&F'&\hfl{h'}{}&F''\cr\vfl{g}{}&&\vfl{\varphi}{}&&\vfl{\varphi'}{}\cr X&\hfl{u}{}&Y&\hfl{v}{}&Z\cr}$$
and a $1$-map $s':C'\fl F''$ such that: $h'h=d(s')$.

Because of Lemma 2.5, there exist a finite $\A$-complex $K$ and a commutative diagram:
$$\diagram{K&\hfl{\alpha}{}&F''\cr\vfl{\beta}{}&&\vfl{\varphi'}{}\cr Y&\hfl{v}{}&Z\cr}$$
where $\alpha$ is a deflation which have a section $\sigma$ in each non-zero degree. By setting $\sigma=0$ in degree $0$, $\sigma$ is a $0$-map from $F''$ to $K$.

We have a commutative diagram:
$$\diagram{C'&\hfl{h}{}&F'\oplus K&\hfl{h'\oplus\alpha}{}&F''\cr\vfl{g}{}&&\vfl{\varphi\oplus\beta}{}&&\vfl{\varphi'}{}\cr X&\hfl{u}{}&Y&\hfl{v}{}&Z\cr}$$

Since $C'_{-1}=0$ and $s'$ is a $1$-map, we have: $\alpha\sigma s'=s'$ and $\varphi\sigma s'=0$. Hence we have:
$$(h'\oplus\alpha)(h-d(\sigma s'))=h'h-\alpha d(\sigma s')=h'h-d(\alpha\sigma s')=h'h-d(s')=0$$
and we get a commutative diagram:
$$\diagram{C'&\hfl{h-d(\sigma s')}{}&F'\oplus K&\hfl{h'\oplus\alpha}{}&F''\cr\vfl{g}{}&&\vfl{\varphi\oplus\beta}{}&&\vfl{\varphi'}{}\cr X&\hfl{u}{}&Y&\hfl{v}{}&Z\cr}$$
with: $(h'\oplus\alpha)(h-d(\sigma s'))=0$.

Since $\alpha$ is a deflation, the morphism $h'\oplus\alpha$ is also a deflation with a finite kernel $F$. Hence the morphism $g:C'\fl X$ factors through the finite 
$\A$-complex $F$ and $f$ factors also through $F$. The result follows.\cqfd
\vskip 12pt
Then the class $\C$ is exact and cocomplete and contains $\A$. Since $\A$ is regular the class $\C$ is the class of all $\A\v$-modules.

Let $a$ be an integer and $C$ be an $\A$-complex such that $H_i(C)=0$ for all $i\geq a$. Let $X$ be the cokernel of $d:C_{a+1}\fl C_a$ and $C'$ be the quotient of $C$
by its $(a-1)$-skeleton. The equality $X=\Coker(d:C_{a+1}\fl C_a)$ induces a morphism $f:C'\fl X$ and, because $X$ is in the class $\C$, $f$ factors through a finite 
complex $F$. Up to killing the $(a-1)$-skeleton of $F$ we may as well suppose that $F_i=0$ for all $i<a$. We may also suppose that $F$ is non-zero. Let $b$ be the largest 
integer such that $F_b\not=0$ and $n$ be an integer with $n>b$.

For each integer $p<n$ denote $\U_p$ the class of pairs $(K,f)$ where $K$ is an $\A$-complex and $f:C\fl K$ is a morphism of complexes such that $f:C_i\fl K_i$ is an 
isomorphism for all $i<p$ and $K_i=0$ for all $i>n+1$.
\vskip 12pt
\noi{\bf 2.9 Lemma:} {\sl The class $\U_a$ is non-empty.}
\vskip 12pt
\noi{\bf Proof:} The morphism $F\fl X$ is a morphism: 
$$\Coker(F_{a+1}\fl F_a)\fl X\build\fl_{}^\sim \Coker(C_{a+1}\fl C_a)\build\fl_{}^\sim\Im(C_a\fl C_{a-1})$$ 
Thus we have an $\A$-complex:
$$K=(\dots\fl F_{a+2}\fl F_{a+1}\fl F_a\fl C_{a-1}\fl C_{a-2}\fl\dots$$
and a morphism $C\fl K$ which belongs to $\U_a$.\cqfd
\vskip 12pt
\noi{\bf 2.10 Lemma:} {\sl Let $p$ be an integer with $a<p\leq n$. Suppose $\U_{p-1}$ is non-empty then $\U_p$ is non-empty.}
\vskip 12pt
\noi{\bf Proof:} Since $a<p$, we have $H_{p-1}(C)=0$. Let $(K,f)$ be an object in $\U_{p-1}$. We have a morphism from $\Coker(K_{p+1}\fl K_p)$ to:
$$\Ker(K_{p-1}\fl K_{p-2})\simeq \Ker(C_{p-1}\fl C_{p-2})\simeq \Coker(C_{p+1}\fl C_p)$$

Let $X$ be the $\A\v$-module $\Coker(C_{p+1}\fl C_p)$. Then we have two morphisms 
$$f:\Coker(C_{p+1}\fl C_p)\fl \Coker(K_{p+1}\fl K_p)\ \hbox{and}\ \varphi:\Coker(K_{p+1}\fl K_p)\fl X$$
such that $\varphi\circ f$ is an isomorphism.

Let $Y$ be the kernel of $\varphi$. We have a decomposition: $ \Coker(K_{p+1}\fl K_p)=X\oplus Y$ and then an epimorphism $g:K_p\fl X\oplus Y$.

Since $K$ is finitely presented, $Y$ is finitely generated and there is an $\A$-module $A$ and an epimorphism $\alpha:A\fl Y$. Since $g$ is an epimorphism, there is a 
deflation $B\fl A$ and a commutative diagram:
$$\diagram{B&\hfl{}{}&A\cr\vfl{\beta}{}&&\vfl{\alpha}{}\cr K_p&\hfl{g}{}&X\oplus Y\cr}$$

Hence, up to replacing $K_{p+1}\build\fl_{}^d K_p$ by $K_{p+1}\oplus B\build\fl_{}^{d\oplus\beta} K_p$, we may as well suppose that $H_p(C)\fl H_p(K)$ is an 
isomorphism and the morphism $C_p\oplus K_{p+1}\fl K_p$ is an epimorphism with kernel $H$. Because of property (cc8), there is an $\A$-module $U$ such that $H\oplus U$ is
an $\A$-module. Then, up to adding to $K_{p+1}$ the module $U$, we may as well suppose that the kernel $H$ is an $\A$-module.

Thus we get a commutative diagram:
$$\diagram{C_{p+1}&\hfl{}{}&H&\hfl{}{}&K_{p+1}\cr\vfl{d}{}&&\vfl{}{}&&\vfl{d}{}\cr C_p&\hfl{=}{}&C_p&\hfl{}{}&K_p\cr}$$
where the right square is cartesian.

Thus we get a morphism $g:C\fl K'$ where $K'$ is the complex:
$$K'=(\dots\fl K_{p+2}\fl H\fl C_p\fl C_{p-1}\fl\dots)$$
and $(K',g)$ belongs to the class $\U_p$.\cqfd
\vskip 12pt
Because of this lemma we see that the classes $\U_a,\U_{a+1},\dots,\U_n$ are non-empty. Let $(K,f)$ be an object of $\U_n$. We have a commutative diagram with 
an exact top line:
$$\diagram{C_{n+2}&\hfl{}{}&C_{n+1}&\hfl{}{}&C_n&\hfl{}{}&C_{n-1}\cr\vfl{}{}&&\vfl{}{}&&\vfl{\sim}{}&&\vfl{\sim}{}\cr 
0&\hfl{}{}&K_{n+1}&\hfl{}{}&K_n&\hfl{}{}&K_{n-1}\cr}$$
and we have morphisms:
$$\Coker(C_{n+2}\fl C_{n+1})\build\fl_{}^\alpha K_{n+1}\build\fl_{}^\beta\Ker(K_n\fl K_{n-1})$$ 

For each integer $k$, denote $Z_k$ the kernel of $d:C_k\fl C_{k-1}$. Thus we have:
$$\Coker(C_{n+2}\fl C_{n+1})\simeq Z_n\hskip 24pt \Ker(K_n\fl K_{n-1})\simeq \Ker(C_n\fl C_{n-1})=Z_n$$
and the composite morphism $\beta\alpha:\Coker(C_{n+2}\fl C_{n+1})\fl\Ker(K_n\fl K_{n-1})$ is an isomorphism. Therefore $Z_n$ is a direct summand of the $\A$-module 
$K_n$ and theorems 2.1 and 6 are proven. \cqfd
\vskip 12pt
Now we are able to give a proof of Theorem 1:
\vskip 12pt
\noi{\bf Proof of Theorem 1:}
Let $A$ be a ring. We have to prove that $A$ is regular coherent on the right if and only if it is regular on the right and coherent on the right. But a regular coherent
ring is coherent. Then we may assume that $A$ is coherent on the right and we have to prove that $A$ is regular coherent on the right if and only if it is regular on the 
right. 

Let $\A$ be the category of finitely generated projective right $A$-modules. The category $\A\v$ of right $A$-modules is a strict cocompletion of $\A$.
Because $A$ is assumed to be coherent on the right, the category $\B$ of finitely presented  right $A$-modules is abelian.

Suppose $A$ is regular coherent on the right. Let $\C$ be a cocomplete exact class of right $A$-modules which contains $\A$. Since $\C$ is cocomplete, $\C$ contains
every projective modules in $\A\v$. Let $M$ be a finitely presented right $A$-module. Since $A$ is regular coherent on the right, $M$ has a projective resolution with
finite length. Thus we have an exact sequence in $\A\v$:
$$0\fl C_n\fl C_{n-1}\fl\dots\fl C_1\fl C_0\fl M\fl0$$
where each $C_i$ is projective. Thus we have an acyclic $\A\v$-complex $C$ (with $C_{-1}=M$). Let $Z_i$ be the kernel of $d:C_i\fl C_{i-1}$. Thus we have: $Z_n=0$, 
$Z_{-1}=M$. We have also exact sequences:
$$0\fl Z_i\fl C_i\fl Z_{i-1}\fl 0$$
Since $\C$ is exact each $Z_i$ is in $\C$ and then $M$ belongs to the class $\C$. Therefore $\C$ contains every finitely presented right $A$-module and, because $\C$ is
cocomplete, every right $A$-module belongs to $\C$. Hence $\A$ is regular and $A$ is regular on the right.

Suppose $A$ is regular on the right. Then $\A$ is regular. Let $M$ be a finitely presented right $A$-module. Since $\B$ is abelian, $M$ has a $\A$-resolution $C$.

Because of Theorem 6, there is an integer $n$ such that the kernel $Z_n$ of $d:C_n\fl C_{n-1}$ is a direct summand of some $\A$-module. Then $Z_n$ is a direct summand
of a finitely generated projective right $A$-module and $M$ has a finite $\A$-resolution:
$$0\fl Z_n\fl C_n\fl C_{n-1}\fl\dots\fl C_1\fl C_0\fl M\fl0$$
Hence $A$ is regular coherent on the right and this completes the proof of Theorem 1.\cqfd
\vskip 12pt 
\noi{\bf 2.11 Lemma:} {\sl Let $\B$ be an essentially small exact category and $\A$ be a fully exact subcategory of $\B$. Suppose $\A\subset\B$ is a domination. Then
any strict cocompletion of $\B$ is a strict cocompletion of $\A$.}
\vskip 12pt 
\noi{\bf Proof:} Since $\A\subset\B$ is a domination we have the following properties:

1) for any $\B$-module $Y$, there is an $\A$-module $X$ and a $\B$-deflation $X\fl Y$

2) for any $X$ in $\A$, any $Y$ in $\B$ and any $\B$-deflation $f:Y\fl X$, there is an $\A$-module $Z$ and a morphism $g:Z\fl Y$ such that $f\circ g:Z\fl X$ is an 
$\A$-deflation.

Let $\B\v$ be a strict cocompletion of $\B$. Consider the inclusion $\A\subset \B\v$. Since $\A$ is fully exact in $\B$, conditions (cc0), (cc1) and (cc2) for this 
inclusion are easy to check. Since $\A\subset\B$ is a domination the condition (cc3) is also satisfied. Consider the condition (cc4):

Let $X$ be an $\A$-module, $\varphi:U\fl V$ be a $\B\v$-epimorphism and $f:X\fl V$ be a morphism. Since $\B\v$ is a strict cocompletion of $\B$, there is a $\B$-module $Y$
and a commutative diagram:
$$\diagram{Y&\hfl{}{}&U\cr\vfl{g}{}&&\vfl{\varphi}{}\cr X&\hfl{f}{}&V\cr}$$
where $g$ is a $\B$-deflation and, since $\A\subset\B$ is a domination, we have a commutative diagram:
$$\diagram{Z&\hfl{}{}&Y&\hfl{}{}&U\cr\vfl{h}{}&&\vfl{g}{}&&\vfl{\varphi}{}\cr X&\hfl{=}{}&X&\hfl{f}{}&V\cr}$$
where $Z$ is an $\A$-module and $h$ an $\A$-deflation.

Hence condition (cc4) is satisfied and $\B\v$ is a strict cocompletion of $\A$.\cqfd
\vskip 12pt

\noi{\bf 2.12 Theorem:} {\sl Let $\B$ be an essentially small exact category and $\A$ be a fully exact subcategory of $\B$. Suppose $\A$ is regular and $\A\subset\B$ is a 
domination. Then $\B$ is regular and the morphism $K_i(\A)\fl K_i(\B)$ induced by the inclusion $\A\subset \B$ is an isomorphism for $i>0$ and a monomorphism for $i=0$.}
\vskip 12pt
\noi{\bf Proof:} Let $\B\v$ be a strict cocompletion of $\B$. Because of Lemma 2.11, $\B\v$ is a strict cocompletion of $\A$.

Let $\C$ be a cocomplete exact class in $\B\v$ containing $\B$. Then $\C$ contains $\A$ and, because $\A$ is regular, $\C$ is the class of all $\B\v$-modules. Therefore
$\B$ is regular.

Let $\A'$ be the class of direct summand in $\B$ of $\A$-modules. This class is a fully exact subcategory of $\B$. Since $\A$ and $\B$ are stable under extension in 
$\B\v$, $\A'$ is also stable under extension in $\B$. Moreover $\A$ is cofinal in $\A'$ and Quillen's Cofinality Theorem [Q] implies that the morphism 
$K_i(\A)\fl K_i(\A')$ induced in K-theory by the inclusion $\A\subset \A'$ is an isomorphism for $i>0$ and a monomorphism for $i=0$.

Let $0\fl X\fl Y\fl Z\fl 0$ be an exact sequence in $\B$. Suppose $X$ and $Z$ are in $\A'$. Then there are $\A$-modules $A$ and $B$ with $X\oplus A$ and $Z\oplus B$ in
$\A$ and we have an exact sequence:
$$0\fl X\oplus A\fl Y\oplus A\oplus B\fl Z\oplus B\fl 0$$
and $Y$ is in $\A'$. Hence $\A'$ is stable under extension in $\B$.

Suppose $Y$ and $Z$ are in $\A'$. Then there are $\A$-modules $A$ and $B$ with $Y\oplus A$ and $Z\oplus B$ in $\A$ and we have an exact sequence:
$$0\fl X\fl Y\oplus A\oplus B\fl Z\oplus A\oplus B\fl 0$$
and $X$ is the kernel of an $\A$-morphism which is an epimorphism in $\B\v$. Hence $\A'$ is stable under kernel of epimorphisms in $\A$.

Let $X$ be a $\B$-module. Since $\A\subset\B$ is a domination, there is a long exact sequence in $\B\v$:
$$\dots\fl C_2\fl C_1\fl C_0\fl X\fl 0$$
where $C_i$ is in $\A$ for each $i\geq 0$ and each kernel $Z_i$ of $d:C_i\fl C_{i-1}$ is in $\B$. 

Because of Theorem 2.1 there is an integer $n>0$ such that the kernel $Z_n$ of $d:C_n\fl C_{n-1}$ is a direct summand of an $\A$-module $K$ and $X$ has a finite 
$\A'$-resolution. Hence Quillen's Resolution Theorem [Q] implies that the inclusion $\A'\subset \B$ induces an isomorphism in K-theory. The result follows.\cqfd
\vskip 12pt
\noi{\bf 2.13 Corollary:} {\sl Let $\A$ be an essentially small exact category and $\A\v$ be a strict cocompletion of $\A$. Suppose $\A$ is stable under direct summand in
$\A\v$ and $\A$ is regular. Then any $\A\v$-module having an $\A$-resolution has an $\A$-resolution of finite length.}
\vskip 12pt
\noi{\bf Proof:} Let $\B$ be the category of  $\A\v$-module $X$ such that there is an exact sequence in $\A\v$:
$$\dots\fl C_n\fl C_{n-1}\fl\dots\fl C_1\fl C_0\fl X\fl 0$$
where each $\C_i$ is an $\A$-module. We can verified that $\B$ is an essentially small exact subcategory of $\A\v$ and that the inclusion $\A\subset\B$ is a domination. 

Hence Theorem 2.12 implies the result.\cqfd
\vskip 12pt
\noi{\bf 2.14 Remark:} Because of this result the notion of regularity presented here is stronger than many other similar notion appearing in the literature, at least 
for rings (see [AMM], [Ge], [Gl], [Na]).
\vskip 12pt
\noi{\bf 2.15 Theorem:} {\sl Let $\B$ be an essentially small exact category and $\A$ be a fully exact subcategory of $\B$. Suppose $\A$ is regular and $\A\subset\B$ is a 
domination.Then a $\B$-module has a finite $\A$-resolution if and only if its class in $K_0(\B)$ belongs to the image of $K_0(\A)\fl K_0(\B)$.

Moreover, for any $\B$-module $X$ with $[X]\in K_0(\A)$ and any $\A$-resolution of $X$:

$$\dots\fl C_n\fl C_{n-1}\fl\dots\fl C_1\fl C_0\fl X\fl 0$$
there is an integer $p>0$ with the following property:

For any integer $n\geq p$, there is an $\A$-module $U$ such that Ker$(C_n\fl C_{n-1})\oplus U$ belongs to $\A$.}
\vskip 12pt
\noi{\bf Proof:} Let $\sim$ be the following relation on $\B$:

If $X$ and $Y$ are $\B$-module, we said: $X\sim Y$ if there is a $\B$-module $Z$ and two $\B$-inflations $X\fl Z$ and $Y\fl Z$ with cokernels in $\A$.
\vskip 12pt
\noi{\bf 2.16 Lemma:} {\sl The relation is an equivalence relation and the direct sum operation induces on the quotient set $M=\B/\!\sim$ a commutative monoid structure.}
\vskip 12pt
\noi{\bf Proof:} It is clear that $\sim$ is symmetric and reflexive. Moreover two isomorphic $\B$-modules are equivalent. Let $X$, $Y$ and $Z$ be three $\B$-modules.
Suppose: $X\sim Y$ and $Y\sim Z$. Then there exist two $\B$-modules $U$ and $V$ and $\B$-inflations $\alpha:X\fl U$, $\beta:Y\fl U$, $\gamma:Y\fl V$ and $\delta:Z\fl V$
with cokernels in $\A$. Thus we get a commutative diagram:
$$\diagram{&&&&Z\cr&&&&\vfl{\delta}{}\cr &&Y&\hfl{\gamma}{}&V\cr&&\vfl{\beta}{}&&\vfl{g}{}\cr X&\hfl{\alpha}{}&U&\hfl{f}{}&W\cr}$$
where $W$ is obtained by cobase change.

We have $\B$-conflations:
$$\Coker(\alpha)\fl\Coker(f\alpha)\fl\Coker(f)\ \ \hbox{and}\ \ \Coker(\delta)\fl\Coker(g\delta)\fl\Coker(g)$$
But we have: $\Coker(f)\simeq\Coker(\gamma)$ and $\Coker(g)\simeq\Coker(\beta)$ and these cokernels are in $\A$.

On the other hand, $\B\v$ is a strict cocompletion of $\A$ and $\A$ is stable under extensions in $\B\v$ and then in $\B$. Hence $\Coker(f\alpha)$ and $\Coker(g\delta)$ are 
in $\A$. Thus we have: $X\sim Z$ and the relation is transitive.
Therefore $\sim$ is an equivalence relation and we have a quotient set $M$. The direct sum operation is compatible with this relation and induces on $M$ a commutative 
monoid structure.\cqfd
\vskip 12pt
Let $K_0$ be cokernel of the morphism $K_0(\A)\fl K_0(\B)$ induced by the inclusion. Every $\B$-module $X$ has a class $<\! X\! >$ in $M$ and a class $[X]$ in $K_0$.  
\vskip 12pt 
\noi{\bf 2.17 Lemma:} {\sl There are unique morphisms of monoids $\varepsilon:M\fl K_0$ and $\theta:M\fl M$ such that: 
$$\varepsilon(<\! X\! >)=[X]\ \ \hbox{and}\ \ \theta(<\! X\! >)=<\! Y\! >$$ 
for any $X\in\B$ and any conflation $Y\fl A\fl X$ with $A\in \A$.

Moreover, for any $x\in M$, we have: $\varepsilon(\theta x)=-\varepsilon(x)$.}
\vskip12pt
\noi{\bf Proof:} Let $X$ and $Y$ be two $\B$-modules. Suppose that $X\sim Y$. Then there is a $\B$-module $Z$ and two conflations: $X\fl Z\fl A$ and $Y\fl Z\fl B$ with
$A$ and $B$ in $\A$.

Thus we have in $K_0$
$$[X]=[Z]-[A]=[Z]=[Y]+[B]=[Y]$$
And the map $\B\fl K_0$ factors through $M$ via a morphism of monoids $\varepsilon:M\fl K_0$ and we have $\varepsilon(<\! X\! >)=[X]$ for any $X\in\B$. 

Let $U\build\fl_{}^u A\build\fl_{}^\alpha X$ and $V\build\fl_{}^vB\build\fl_{}^\beta Z$ be two $\B$-conflations with $A$ and $B$ in $\A$. Suppose: $X\sim Y$. Then 
there is a $\B$-module $Z$ and two inflations $f:X\fl Z$ and $g:Y\fl Z$ with cokernels in $\A$. Let $h:C_0\fl Z$ be a deflation with $C_0\in\A$. Set: 
$C=C_0\oplus A\oplus B$. We have a deflation $\gamma=h\oplus f\alpha\oplus g\beta:C\fl Z$ with kernel $W$ and a commutative diagram:
$$\diagram{U&\hfl{\lambda}{}&W&\bhfl{\mu}{}&V\cr\vfl{u}{}&&\vfl{w}{}&&\vfl{v}{}\cr A&\hfl{\varphi}{}&C&\bhfl{\psi}{}&B\cr\vfl{\alpha}{}&&\vfl{\gamma}{}&&\vfl{\beta}{}\cr 
X&\hfl{f}{}&Z&\bhfl{g}{}&Y\cr}$$
where each column is a conflation, morphisms $\lambda$, $\mu$, $f$ and $g$ are $\B$-inflations and morphisms $\varphi$ and $\psi$ are $\A$-inflations. Therefore we get 
$\B$-conflations:
$$\Coker(\lambda)\fl\Coker(\varphi)\fl\Coker(f)\hskip 12pt\hbox{and}\hskip 12pt \Coker(\mu)\fl\Coker(\psi)\fl\Coker(g)$$
Since $\varphi$ and $\psi$ are $\A$-inflations, modules $\Coker(\varphi)$ and $\Coker(\psi)$ are in $\A$ and $\Coker(\lambda)$ and $\Coker(\mu)$ are kernels  of
epimorphisms in $\A$. Because of the condition (cc8) there are two $\A$-modules $P$ and $Q$ such that $\Coker(\lambda)\oplus P$ and $\Coker(\mu)\oplus Q$ are in $\A$.

Hence up to adding $P\oplus Q$ to $C$, we may as well suppose that $\Coker(\lambda)$ and $\Coker(\mu)$ are $\A$-modules and we have: $U\sim V$. Therefore the class 
$<\! U\! >$ depends only on $<\! X\! >$ and we have a well defined map $\theta:M\fl M$ such that, for each $\B$-conflation $Y\fl A\fl X$ with $A\in\A$, we have:
$<\! Y\! >=\theta(<\! X\! >)$. It is easy to see that $\theta$ is compatible with direct sum and $\theta$ is a monoid homomorphism.

Moreover, for any $\B$-module $X$ and any conflation $Y\fl A\fl X$ with $A\in \A$, we have:
$$\varepsilon(\theta(<\! X\! >))=\varepsilon(<\! Y\! >)=[Y]=-[X]=-\varepsilon(<\! X\! >)$$
and, for any $x\in M$, we have: $\varepsilon(\theta x)=-\varepsilon(x)$.\cqfd 
\vskip 12pt
Let $M_0$ be the set of elements $x\in M$ such that $x+\theta(x)=0$. This set is a sub-monoid of $M$ and also an abelian group.
\vskip 12pt
\noi{\bf 2.18 Lemma:} {\sl There is a unique monoid homomorphism $\lambda:M\fl M_0$ which is the identity on $M_0$ such that:

$\bullet$ for each $x\in M$ we have: $\lambda\theta(x)=-\lambda(x)$

$\bullet$ for each $x\in M$, there is an integer $n\geq 0$ such that $\theta^p(x)=(-1)^p\lambda(x)$ for any $p\geq n$.}
\vskip 12pt
\noi{\bf Proof:} Consider the subcategory $\A'$ of direct summand of $\A$-modules defined in the proof of Theorem 2.12. Let $X$ be a $\A'$-module. Then there is an 
$\B$-module $Y$ with $X\oplus Y\in \A$. Let: $x=<\! X\! >$ and $y=<\! Y\! >$. Because of the conflations $X\fl X\oplus Y\fl Y$ and $0\fl X\oplus Y\fl X\oplus Y$, we have: 
$\theta(y)=x$ and $\theta(x+y)=0$ and then:
$$0=\theta(x)+\theta(y)=\theta(x)+x$$
and $x$ belongs to $M_0$.

Let $X$ be a $\B$-module and $x$ be its class in $M$. We construct by induction a sequence of $\B$-conflations $X_{n+1}\fl A_n\fl X_n$, for $n\geq0$, such that: $X_0=X$ 
and $A_n\in\A$ for all $n$. Then we have for all $n\geq0$: $<\! X_n\! >=\theta^n(x)$.

Let $A$ be the $\A$-complex:
$$A=(\dots\fl A_2\fl A_1\fl A_0\fl 0\fl\dots)$$
We have: $H_i(A)=0$ for all $i>0$ and Theorem 2.1 implies that there is an integer $n>0$ such that the kernel $X_n$ of $d:A_{n-1}\fl A_{n-2}$ is a direct summand of some
$\A$-module. Then there is an $\B$-module $Y$ such that $X_n\oplus Y$ is in $\A$.

Since $X_n\oplus Y$ is an $\A$-module, we have: $<\! X_n\! >+<\! Y\! >=<\! X_n\oplus Y\! >=0$ and, because of the conflation $Y\fl X_n\oplus Y\fl X_n$, we have:
$\theta(<\! X_n>\! )=<\! Y\! >$. Hence we get: $\theta^n(x)+\theta(\theta^n(x))=0$ and $\theta^n(x)$ belongs to $M_0$.

Since $\theta^n(x)$ belongs to $M_0$, for any $p\geq0$, we have: $\theta^{n+p}(x)=(-1)^p\theta^n(x)$ and there is a unique $y\in M_0$ such that: $\theta^p(x)=(-1)^py$
for any $p\geq n$.

Then there is a unique map $\lambda:M\fl M_0$ such that: $\theta^p(x)=(-1)^p\lambda(x)$ for any $p$ large enough and it is easy to see that $\lambda$ is a monoid 
homomorphism.\cqfd
\vskip 12pt
\noi{\bf 2.19 Lemma:} {\sl There is a unique isomorphism of groups $\mu:K_0\build\fl_{}^\sim M_0$ such that: $\mu\varepsilon=\lambda$.}
\vskip 12pt
\noi{\bf Proof:} Consider a $\B$-conflation: $X\fl Y\fl Z$. Let $A\fl X$ and $B\fl Y$ two $\B$-deflations with $A$ and $B$ in $\A$. Then we have a commutative diagram
where lines and columns are conflations: 
$$\diagram{U&\hfl{}{}&V&\hfl{}{}&W\cr\vfl{}{}&&\vfl{}{}&&\vfl{}{}\cr A&\hfl{}{}&A\oplus B&\hfl{}{}&B\cr \vfl{}{}&&\vfl{}{}&&\vfl{}{}\cr X&\hfl{}{}&Y&\hfl{}{}&Z\cr}$$
Hence we can construct a sequence of $\B$-conflations $K_n=(X_n\fl Y_n\fl Z_n)$ with $K_0=(X\fl Y\fl Z)$ and $K_1=(U\fl V\fl W)$, a sequence of $\A$-conflations 
$H_n=(A_n\fl B_n\fl C_n)$ with $H_0=(A\fl A\oplus B\fl B)$ and a sequence of morphisms: $K_{n+1}\fl H_n\fl K_n$. 

Hence, for each $n\geq0$, we have a commutative diagram where lines and columns are conflations:
$$\diagram{X_{n+1}&\hfl{}{}&Y_{n+1}&\hfl{}{}&Z_{n+1}\cr\vfl{}{}&&\vfl{}{}&&\vfl{}{}\cr A_n&\hfl{}{}&B_n&\hfl{}{}&C_n\cr \vfl{}{}&&\vfl{}{}&&\vfl{}{}\cr 
X_n&\hfl{}{}&Y_n&\hfl{}{}&Z_n\cr}$$

Hence, for $n$ large enough, we have: $<\! X_n\! >=\theta^n(<\! X\! >)$, $<\! Y_n\! >=\theta^n(<\! Y\! >)$ and  $<\! Z_n\! >=\theta^n(<\! Z\! >)$. Then, for $n$ large 
enough, there are $\B$-modules $U,V,W$ such that: $X_n\oplus U$, $Y_n\oplus V$ and $Z_n\oplus X$ are in $\A$ and we have an $\A$-conflation:
$$X_n\oplus U\oplus V\fl Y_n\oplus U\oplus V\oplus W\fl Z_n\oplus W$$

Let $x$, $y$, $z$, $u$, $v$ and $w$ be the classes of $X$, $Y$, $Z$, $U$, $V$ and $W$ in $M$. Thus we have:
$$(-1)^n\lambda(x)+u+v=(-1)^n\lambda(y)+u+v+w=(-1)^n\lambda(z)+w=0$$
and then:
$$(-1)^n\lambda(x)+\lambda(u)+\lambda(v)=(-1)^n\lambda(y)+\lambda(u)+\lambda(v)+\lambda(w)=(-1)^n\lambda(z)+\lambda(w)=0$$
$$\Longrightarrow\ \ \lambda(y)=\lambda(x)+\lambda(z)$$

Hence, for any $A$-module $X$ we have: $\lambda(<\! X\! >)=0$ and for any $\B$-conflation $X\fl Y\fl Z$ we have: $\lambda(<\! Y\! >)=\lambda(<\! X\! >)+\lambda(<\! Z\! >)$.
Then there is a unique morphism $\mu:K_0\fl M_0$ such that: $\lambda=\mu\varepsilon$.

Since $\lambda:M\fl M_0$ is the identity on $M_0$, $\mu$ is onto. Let $u$ be an element of the kernel of $\mu$ and $x$ be an element of $M$ with $\varepsilon(x)=u$.
For any integer $n>0$, we have:
$$u=\varepsilon(x)=(-1)^n\varepsilon(\theta^n x)=\varepsilon(\lambda x)$$
and there is $y\in M_0$ such that: $u=\varepsilon(y)$. Then we have:
$$0=\mu(u)=\mu\varepsilon(y)=\lambda(y)=y$$
and $\mu$ is an isomorphism.\cqfd
\vskip 12pt
Now we are able to finish the proof of Theorem  2.15.

Let $\B'$ the subcategory of $\B\v$ defined by: a $\B\v$-module $X$ is in $\B'$ if an only if there is a $\B$-module $Y$ with; $X\oplus Y\in\B$. This category is
fully exact in $\B\v$ and $\B$ is strictly cofinal in $\B'$. Then $K_0(\B)\fl K_0(\B')$ is an isomorphism and each $\B'$-module has a class in $K_0$.

If $k\geq0$ is an integer, we denote $P(k)$ the following property:

$\bullet$ if a $\B'$-module $X$ has an $\A$-resolution of length $k$ then $[X]=0$ in $K_0$.

The property $P(0)$ is clearly true.

Let $k>0$ be an integer and suppose $P(k-1)$ is true. Let $X$ be a $\B'$-module having an $\A$-resolution of length $k$:
$$0\fl C_k\fl C_{k-1}\fl\dots\fl C_1\fl C_0\fl X\fl 0$$
The $\B'$-morphism $C_0\fl X$ is an epimorphism and, because of property (cc8), its kernel $K$ is in $\B'$. Thus we have an exact sequence:
$$0\fl C_k\fl C_{k-1}\fl\dots\fl C_1\fl K\fl0$$
and, because of the property $P(k-1)$, we have $[K]=0$ in $K_0$. Then, because of the exact sequence $0\fl K\fl C_0\fl X\fl 0$, we have: $[X]=[C_0]-[K]=0$ and the property
is true. 

Therefore properties $P(k)$ are all true and, for each $\B$-module $X$ having a finite $\A$-resolution, we have $[X]=0$ in $K_0$.

Let $X$ be a $\B$-module with $[X]=0$ in $K_0$. Then we have: $0=\mu([X])=\lambda(<\! X\! >)$ and there is an integer $p$ such that $\theta^p(<\! X\! >)=0$.

Let: $\dots\fl C_2\fl C_1\fl C_0\fl X\fl 0$ be an $\A$-resolution of $X$. For any $i\geq 0$ we set:
$$Z_i=\left\{\matrix{\hbox{Ker}(C_0\fl X)&\hbox{if}\ i=0\cr \hbox{Ker}(C_i\fl C_{i-1})&\hbox{if}\ i>0\cr}\right.$$
Each $Z_i$ belongs to $\A\v=\B\v$ and we have exact sequences:
$$0\fl Z_0\fl C_0\fl X\fl 0$$
$$0\fl Z_i\fl C_i\fl Z_{i-1}\fl 0$$
in $\A\v$. Then, because of property (cc8), there are $\A$-modules $A_i$ such that each module $Z_i\oplus A_i$ belongs to $\B$ and we have $\B$-conflations:
$$Z_0\oplus A_0\fl C_0\oplus A_0\fl X$$
$$Z_i\oplus A_i\fl C_i\oplus A_i\oplus A_{i-1}\fl Z_{i-1}\oplus A_{i-1}$$
Then, for any $i\geq0$, we have: $<\! Z_i\oplus A_i\! >=\theta^{i+1}(<\! X\! >$. Hence, for $n\geq p$, we have: $<\! Z_n\oplus A_n\! >=0$ and there are two $\A$-modules
$U$ and $V$ and an exact sequence:
$$0\fl Z_n\oplus A_n\fl U\fl V\fl 0$$
Because of property (cc8), there is an $\A$-module $W$ such that $Z_n\oplus A_n\oplus W$ belongs to $\A$. Hence Theorem 2.15 and also Domination Theorem 5 are 
proven.\cqfd
\vskip 12pt
\noi{\bf Proof of Theorem 7:} Let $A$ be a left-regular ring, $\A$ be the exact category of finitely generated projective left $A$-module and $\A\v$ be 
the category of all left $A$-modules. Then $\A$ is a essentially small regular exact category and $\A\v$ is a strict cocompletion of $\A$. Moreover $\A$ is regular.

Let $C$ be an acyclic complex of projective right $A$-modules. We have to prove that $C$ is contractible.

Let $\C$ be the class of left $A$-modules $X$ such that $C\build\otimes_A^{}X$ is acyclic.

Consider a short exact sequence of left $A$-modules: $0\fl X\fl Y\fl Z\fl 0$. Since $C$ is a complex of projective modules, we have an short exact sequence of complexes:
$$0\fl C\build\otimes_A^{}X\fl C\build\otimes_A^{}Y\fl C\build\otimes_A^{}Z\fl 0$$
and a long exact sequence $(S)$ in homology:
$$\dots\fl H_{i+1}(C\build\otimes_A^{}Z)\build\fl_{}^\partial H_i(C\build\otimes_A^{}X)\fl H_i(C\build\otimes_A^{}Y)\fl H_i(C\build\otimes_A^{}Z)
\build\fl_{}^\partial \dots$$
If two of the modules $X,Y,Z$ are in $\C$, two out of three of the modules in $(S)$ are zero and it is the same for the other ones. Therefore, the third module is in $\C$ 
and the class $\C$ is exact.

Suppose $X$ is the colimit of a filtered system of modules $X_i\in\C$. Then we have:
$$H_*(C\build\otimes_A^{}\ \build\lim_\longrightarrow^{} X_i)\simeq H_*(\build\lim_\longrightarrow^{} C\build\otimes_A^{}X_i)\simeq\ 
\build\lim_\longrightarrow^{} H_*(C\build\otimes_A^{}X_i)\simeq0$$
and $\C$ is cocomplete.

Moreover, for every flat module $X\in\A\v$, the functor $-\build\otimes_A^{}X$ is exact and $C\build\otimes_A^{}X$ is acyclic. Then $\C$ contains all flat modules and
$\C$ contains $\A$. Since $\A$ is regular, $\C$ is the class of all modules in $\A\v$.

Let $X$ be a left $A$-module. For every right $A$-module $M$, we set: $\widehat M=M\build\otimes_A^{}X$. Denote also $Z_n$ the modules of $n$-cycles of $C$.

Since $C$ is acyclic, we have an exact sequence:
$$C_{n+2}\fl C_{n+1}\fl Z_n\fl 0$$
which implies the following exact sequence:
$$\widehat C_{n+2}\fl \widehat C_{n+1}\fl \widehat Z_n\fl 0$$
Since $X$ is in $\C$, $\widehat C$ is acyclic and we have:
$$\widehat Z_n\simeq \Coker(\widehat C_{n+2}\fl \widehat C_{n+1})\simeq \Im(\widehat C_{n+1}\fl \widehat C_n)\simeq \Ker(\widehat C_n\fl \widehat C_{n-1})$$
Thus we have an exact sequence:
$$0\fl \widehat Z_{n+1}\fl \widehat C_{n+1}\fl \widehat Z_n\fl 0$$
On the other hand, we have an exact sequence:
$$0\fl Z_{n+1}\fl C_{n+1}\fl Z_n\fl 0$$
which implies an exact sequence:
$$0\fl \hbox{Tor}_1^A(Z_n,X)\fl \widehat Z_{n+1}\fl \widehat C_{n+1}\fl \widehat Z_n\fl 0$$
Therefore Tor$_1^A(Z_n,X)$ is the trivial module and, because that's true for every left $A$-module $X$, $Z_n$ is a flat module. 

Hence all modules $Z_n$ are flat and Theorem 8.6 of Neeman [Ne] implies the result.\cqfd
\vskip 24pt
\noi{\bf 3 Regular rings and regular groups.}
\vskip 12pt
We say that a ring $A$ is regular on the right (or right-regular) if the category of finitely generated projective right $A$-modules is regular. If the category of
finitely generated left $A$-module is regular, we say that $A$ is regular on the left (or left-regular).

Let $\A$ be an essentially small exact category and $\A\v$ be a strict cocompletion of $\A$. Let $\C_0$ be the smallest cocomplete exact class in $\A\v$ containing 
$\A$. We say that a $\A\v$-module is $\A$-regular if it belongs to $\C_0$. Thus $\A$ is regular if and only if each $\A\v$-module is $\A$-regular.

If $\A$ is a category of finitely generated projective modules over a ring $A$, the $\A$-regularity property will be also called $A$-regularity.
\vskip 12pt
\noi{\bf 3.1 Lemma:} {\sl Let $\alpha:A\fl B$ be a ring homomorphism and $M$ be a right $A$-module. Suppose $B$ is flat, as a left $A$-module and $M$ is $A$-regular. Then
the right $B$-module $M\build\otimes_A^{}B$ is $B$-regular.}
\vskip 12pt
\noi{\bf Proof:} Let $\C$ be the class of right $A$-modules $M$ such that $M\build\otimes_A^{}B$ is $B$-regular. Since $B$ is flat, as a left $A$-module, the class $\C$ is
cocomplete and exact. Moreover $\C$ contains each finitely generated projective $A$-module and $\C$ contains each $A$-regular modules. In particular 
$M\build\otimes_A^{}B$ is $B$-regular.\cqfd
\vskip 12pt
\noi{\bf 3.2 Proposition:} {\sl Let $A_i$, $i\in I$, be a filtered system of rings and $A$ be its colimit. Suppose each $A_i$ is regular on the right and, for each 
morphism $A_i\fl A_j$, $A_j$ is flat as a left $A_i$-module. Then $A$ is regular on the right.}
\vskip 12pt
\noi{\bf Proof:} Let $i$ be an element of $I$. Since, for every $i\fl j$, $A_j$ is flat as a left $A_i$-module, $A$ is also flat as a left $A_i$-module.

Let $M$ be a finitely presented right $A$-module. We have an exact sequence:
$$A^p\build\fl_{}^f A^q\fl M\fl 0$$
and $f$ is represented by a matrix with entries in $A$. This matrix can be lift to a matrix with entries in some $A_i$. Thus we have an exact sequence:
$$A_i^p\build\fl_{}^{f_i} A_i^q\fl M_i\fl 0$$
and then an isomorphism: $M\simeq M_i\build\otimes_{A_i}^{} A$. 

Since $A_i$ is regular on the right, $M_i$ is $A_i$-regular and, because of Lemma 3.1, $M$ is $A$-regular. Therefore, the class of $A$-regular modules contains each 
finitely presented module and each filtered colimit of finitely presented modules. Then every right $A$-module is $A$-regular and $A$ is regular on the right.\cqfd
\vskip 12pt
\noi{\bf 3.3 Lemma:} {\sl Let $\alpha:A\fl B$ be a ring homomorphism. Suppose that $\alpha$ induces a split monomorphism of right $A$-modules, that $B$ is right-regular 
and that $B$ is regular as a right $A$-module. Then $A$ is right-regular.} 
\vskip 12pt
\noi{\bf Proof:} We may as well suppose that $\alpha$ is an inclusion and there is a right $A$-module $K$ such that: $B=A\oplus K$ (as right $A$-modules).

Let $\M_A$ be the category of right $A$-modules and $\M_B$ be the category or right $B$-modules. For each module $M$ in $\M_B$, we may consider $M$ as a right $A$-module.
Thus we get a functor $F:\M_B\fl \M_A$ which is clearly exact and commutes with filtered colimit.

Let $\C$ be the class of modules $M\in\M_B$ such that $F(M)$ is regular. This class is exact and cocomplete. Since $F(B)$ is regular $\C$ contains $B$. Then $\C$ contains
all finitely generated modules in $\M_B$ and, because $B$ is right-regular, $\C$ is the class of all right $B$-modules.

Let $E$ be a right $A$-module. The module $M\build\otimes_A^{} B$ belongs to $\C$ and $F(M\build\otimes_A^{} B)=M\oplus M\build\otimes_A^{} K$ is regular. Thus $M$ is a
direct summand of a regular module and $M$ is regular. Hence $A$ is right-regular.\cqfd
\vskip 12pt
\noi{\bf 3.4 Proposition:} {\sl Let $A$, $B$ and $C$ be three rings and $\alpha:C\fl A$ and $\beta:C\fl B$ be injective ring homomorphisms. Let $R$ be the ring defined
by the cocartesian diagram:
$$\diagram{C&\hfl{\alpha}{}&A\cr\vfl{\beta}{}&&\vfl{}{}\cr B&\hfl{}{}&R\cr}$$
Suppose $A$, $B$ and $C$ are right-regular and $C$-bimodules $\Coker(\alpha)$ and $\Coker(\beta)$ are flat on the left. Then $R$ is right-regular.}
\vskip 12pt
\noi{\bf Proof:} We may as well suppose that $\alpha$ and $\beta$ are inclusions. The tensor product over $C$ of two modules $M$ and $N$ will be simply denoted $M.N$.

Consider the sequence:
$$0\fl R\build\otimes_C^{} R\build\fl_{}^f R\build\otimes_A^{} R\oplus R\build\otimes_B^{} R\build\fl_{}^g R\fl 0\leqno{(S)}$$
where $f$ and $g$ are defined by:
$$f(u\otimes v)=f_A(u\otimes v)\oplus f_B(u\otimes v)=(u\otimes v)\oplus (-u\otimes v)\hskip 24pt g((u\otimes v)\oplus0)=g(0\oplus (u\otimes v))=uv$$
for all $(u,v)\in R\times R$.

We verify that $g$ sends bijectively the cokernel of $f$ to $R$. Then in order to show that $(S)$ is exact it suffices to prove that $f$ is one to one.

We have two left-flat $C$-bimodules $A'=A/C$ and $B'=B/C$.

The ring $R$ has a filtration by $C$-bimodules:
$$C=R_0\subset R_1\subset R_2\subset R_3\subset\dots$$
where $R_1$ is generated by images of $A$ and $B$ and $R_p$ is the bimodule $(R_1)^p$.

Set: $S_p=R_p/R_{p-1}$ (with $P_{-1}=0$). Each $S_p$ is a $C$-bimodule and the direct sum of these $S_p$ is a graded ring. Moreover, since the tensor product is written 
as a simple product, we have:
$$S_0=C$$
$$S_1=A'\oplus B'$$
$$S_2=A'.B'\oplus B'.A'$$
$$S_3=A'.B'.A'\oplus B'.A'.B'$$
and so one.

Let $u\in\Ker(f)$ and $p$ be an integer such that $u$ belongs to $R.R_p$. If $p>0$ then the class of $u$ in $R.S_p$ is on the form: $u=v\oplus v'$ with
$v\in R.A'.B'.A'\dots$ and $v'\in R.B'.A'.B'\dots$. By testing the condition $f_A(u)=0$ modulo $R.R_{p-1}$, we get: $v'=0$ and by doing the same with $f_B$, we get $v=0$.
Therefore, $u$ is in $R.R_{p-1}$. Hence, $u$ is in $R.C$ and we have: $u=x\otimes 1$. Therefore we have $x=0$. Thus $f$ is injective and the sequence $(S)$ is exact.

Let $M$ be a right $R$-module. Because of Lemma 3.1, modules $M\build\otimes_A^{}R$, $M\build\otimes_B^{}R$ and $M\build\otimes_C^{}R$ are  $R$-regular.

By tensoring $(S)$ by $M$ we get an exact sequence:
$$0\fl M\build\otimes_C^{}R\fl M\build\otimes_A^{}R\oplus  M\build\otimes_B^{}R\fl M\fl 0$$
and $M$ is $R$-regular.

Because that's true for every right $R$-module $M$, the ring $R$ is right-regular.\cqfd
\vskip 12pt
\noi{\bf 3.5 Proposition:} {\sl Let $A$ and $C$ be two rings and $\alpha$ and $\beta$ be two injective ring homomorphisms from $C$ to $A$. Let $R$ be the ring generated
by $A$ and a formal invertible element $t$ with the only relations:
$$\forall x\in C,\ \ t\alpha(x)=\beta(x)t$$
Suppose $A$ and $C$ are right-regular and $C$-bimodules $\Coker(\alpha)$ and $\Coker(\beta)$ are flat on the left. Then $R$ is right-regular.}
\vskip 12pt
\noi{\bf Proof:} We may as well suppose that morphisms $\alpha:C\fl A$ and $A\fl R$ are inclusions.
Denote $R_\beta$ the $(R,C)$-bimodule $R$ where the action is defined by: $x.y.c=xy\beta(c)$ for every $(x,y,c)\in R\times R\times C$ and by ${}_\alpha R$ the
$(C,R)$-bimodule $R$ where the action is defined by: $c.x.y=\alpha(c)xy$ for every $(c,x,y)\in C\times R\times R$. Consider the sequence:
$$0\fl R_\beta\build\otimes_C^{} {}_\alpha R\build\fl_{}^f R\build\otimes_A^{} R\build\fl_{}^g R\fl 0\leqno{(S)}$$
where $g$ is the product map and $f$ is defined by:
$$f(u\otimes v)=ut\otimes v-u\otimes tv$$

As above, in order to prove that $(S)$ is exact is suffices to prove that $f$ is injective.

We have in $R$ the following $C$-bimodules:
$${}_\alpha A_\alpha=A\hskip 24pt {}_\alpha A_\beta=At^{-1}\hskip 24pt {}_\beta A_\alpha=tA\hskip 24pt {}_\beta A_\beta=tAt^{-1}$$

We have a filtration of $R$ by $C$-bimodules: $C=R_0\subset R_1\subset R_2\subset\dots$ where $R_1$ is the sum of the ${}_i A_j$, for $i,j$ in $\{\alpha,\beta\}$ and:
$R_p=(R_1)^r$ for all $p>0$.

As above we denote $S$ the graded algebra associated with this filtration. Let $A'$ and $B'$ be the cokernel of $\alpha$ and $\beta$. We have the following bimodules:
$${}_\alpha E_\alpha=A'\hskip 24pt {}_\alpha E_\beta=At^{-1}\hskip 24pt {}_\beta E_\alpha=tA\hskip 24pt {}_\beta E_\beta=tB't^{-1}$$
Thus we have:
$$S_0=C\hskip 24pt S_1=\ \build\oplus_{i,j}^{}{}_i E_j$$
where the direct sum is over all $i,j\in\{\alpha,\beta\}$. Hence, we have, for $p>0$:
$$S_p=\oplus\bigl({}_{i_1}E_{j_1}.{}_{i_2}E_{j_2}.{}_{i_3}E_{j_3}.\cdots .{}_{i_p}E_{j_p}\bigr)$$
where the direct sum is over all sequences $(i_1,i_2,\dots,i_p)$ and $(j_1,j_2,\dots,j_p)$ in $\{\alpha,\beta\}$ such that: $j_k\not=i_{k+1}$ for all $k<p$.

For each $i\in\{\alpha,\beta\}$, denote $S_p(i)$ the direct sum of modules beginning with ${}_i E_\alpha$ or ${}_i E_\beta$. Thus  we have: $S_p=S_p(\alpha)\oplus
S_p(\beta)$ (if $p>0$).

Let $u$ be an element in $R.R$ such that $f(u)=0$. Suppose $u$ is in some $R.R_p$ with $p>0$. The class of $u$ in $R.S_p=R.S_p(\alpha)\oplus R.S_p(\beta)$ is on the form: 
$v\oplus w$ with $v\in R.S_p(\alpha)$ and $w\in R.S_p(\beta)$. Then $f(u)$ is in $R.R_{p+1}$ and its class in $R.S_{p+1}$ is on the form $g(w)$ where $g$ is the morphism
from $R\otimes S_p(\beta)$ to $R\otimes S_{p+1}(\beta)$ sending each $x\otimes y$ to $x\otimes ty$. This morphism is clearly injective and the class of $u$ in $R.S_p$
is an element $v\in R.S_p(\alpha)$. Then $f(u)$ belongs to $R\build\otimes_A^{} R_p$ and its class in $R\build\otimes_A^{} S_p$ is on the form $g'(v)$ where $g'$ is the
morphism $x\otimes y\mapsto -x\otimes ty$ which is injective from $R.S_p(\alpha)$ to $R.S_p(\beta)$. Therefore we have $v=0$ and $u$ belongs to $R.R_{p-1}$.

Hence, if $f(u)=0$, then $u$ is in $R.C$ and we have: $u=x\otimes 1$ for some $x\in R$. Therefore, we have:  $0=f(u)=xt\otimes 1-x\otimes t$ and then: $x=0$.

Thus the sequence $(S)$ is exact and for each right $R$-module $M$, we have an exact sequence:
$$0\fl M\build\otimes_R^{}R_\beta\build\otimes_C^{} {}_\alpha R\fl M\build\otimes_A^{} R\fl_{}^g R\fl M\fl 0$$

On the other hand $R$ is left-flat over $C$ and over $A$. Then, because of Lemma 3.1, $R$-modules $M\build\otimes_R^{}R_\beta\build\otimes_C^{} {}_\alpha R$ and
$M\build\otimes_A^{} R$ are $R$-regular. Therefore, $M$ is regular and, because this is true for each right $R$-module $M$, the ring $R$ is right-regular.\cqfd
\vskip 12pt
\noi{\bf 3.6 Proposition:} {\sl Let $C$ be a right-regular ring and $E$ be a left-flat $C$-bimodule. Then the tensor ring 
$C[E]=C\oplus E\oplus E^{\otimes 2}\otimes\dots$ is right-regular.}
\vskip 12pt
\noi{\bf Proof:} Let $R$ be the tensor ring $C[E]$. We have a filtration $C=R_0\subset R_1\subset R_2\subset\dots$ of $R$ by $C$-bimodules where: $R_1=C\oplus E$ and
$R_p=(R_1)^p$ for each $p>0$. Then, for each $p>0$ we have: $R_p/R_{p-1}\simeq E^{\otimes p}$. Therefore, $R$ is left-flat over $C$.

Consider the following sequence of $C$-bimodules:
$$0\fl R\otimes E\otimes R\build\fl_{}^f R\otimes R\build\fl_{}^g R\fl 0\leqno{(S)}$$
where $g$ is the product and $f$ is the morphism: $x\otimes e\otimes y\mapsto xe\otimes y -x\otimes ey$.

As above in order to prove that this sequence is exact it's enough to prove that $f$ is injective. 

Let $g$ be the morphism $x\otimes e\otimes y\mapsto xe\otimes y$ and $h$ be the morphism $x\otimes e\otimes y\mapsto x\otimes ey$. Then we have: $f=g-h$ and $h$ is 
injective. 

Let $u\in R\otimes E\otimes R$ be an element of the kernel of $f$. We have; $u=u_0+u_1+\dots+u_p$ with $u_k\in R\otimes E\otimes E^{\otimes k}$ and the condition 
$f(u)=0$ is:
$$\forall k,\ \ h(u_{k+1})=g(u_k)$$
Thus it is easy to see that each $u_k$ is null and $f$ is injective.

Hence, the sequence $(S)$ is exact and for each right $R$-module $M$ we have an exact sequence:
$$0\fl M\otimes E\otimes R\build\fl_{}^f M\otimes R\build\fl_{}^g M\fl 0$$
Thus, because of Lemma 3.1, $M$ is regular. Then, because that's true for all right $R$-modules, $R$ is right-regular.\cqfd
\vskip 12pt
Let $G$ be a group. Since the group ring $\Z[G]$ is isomorphic to its opposite ring, $\Z[G]$ is right-regular if and only if it is left-regular. In this case we
simply say that $G$ is regular.
\vskip 12pt
\noi{\bf 3.7 Lemma:} {\sl Let $G$ be a group such that $\Z$, equipped with the trivial $G$-action is a regular $G$-module. Then $G$ is regular. Moreover, if $A$ is
a right-regular ring, the group ring $A[G]$ is right-regular.}
\vskip 12pt
\noi{\bf Proof:} Let $M$ and $N$ be two $G$-modules. We denote $M\otimes N$ the $\Z$-module $M\build\otimes_\Z^{} N$ equipped with the following diagonal right 
$G$-action:
$$\forall x\in M,\ \forall y\in N,\ \forall g\in G,\ \ (x\otimes y)g=xg\otimes yg$$

If $M$ is free over $\Z$, denote $\C(M)$ the class of $G$-modules $N$ such that $M\otimes N$ is $G$-regular. This class is clearly cocomplete and exact. 

Let $M^\circ$ be the module $M$ equipped with the trivial $G$-action. We have an isomorphism $M\otimes\Z[G]\build\fl_{}^\sim M^\circ\otimes \Z[G]$ sending, for each 
$g\in G$, $x\otimes g$ to $xg^{-1}\otimes g$. Then $M\otimes\Z[G]$ is a free $G$-module and $\C(M)$ contains $\Z[G]$. Therefore, $\C(M)$ contains every projective 
$G$-module and then every regular right $G$-module.

Since $\Z$ is regular, $\Z$ belongs to $\C(M)$ and $M\otimes\Z\simeq M$ is regular. Then every $G$-module which is free over $\Z$ is regular and it is easy to
see that $G$ is regular.
\vskip 12pt
Let $\C$ be the class of right $A$-modules $M$ such that $M\otimes\Z[G]$ is $A[G]$-regular. This class is cocomplete and exact and contains $A$. Then, because $A$ is 
right-regular, for every $A$-module $M$, $M\otimes\Z[G]$ is $A[G]$-regular and $M\otimes F$ is $A[G]$-regular for every free $G$-module $F$.

If $M$ is a right $A[G]$-module and $N$ a right $G$-module, we denote $M\otimes N$ the right $A$-module $M\build\otimes_\Z^{} N$ equipped with the diagonal 
right $G$-action.

Let $N$ be a right $G$-module. Let $\alpha_1:F_i\fl N$ and $\alpha_2:F_2\fl N$ be two epimorphisms of right $G$-modules where $F_1$ and $F_2$ are free. Denote $K_1$ and
$K_2$ the kernels of $\alpha_1$ and $\alpha_2$ and $H$ be the kernel of $\alpha_1\oplus \alpha_2:F_1\oplus F_2\fl N$. It is not difficult to see that $H$ is isomorphic
to $K_1\oplus F_2$ and to $K_2\oplus F_1$. Since $M\otimes F_1$ and $M\otimes F_2$ are regular, $M\otimes K_1$ is regular if and only if $M\otimes K_2$ is regular.
Hence, the regularity of $M\otimes K_1$ only depends on $N$.

Let $\C'(M)$ be the class of all right $G$-modules $N$ such that there is an exact sequence of $G$-modules $0\fl K\fl F\fl N\fl 0$ such that  $F$ is free and 
$M\otimes K$ is regular.

Let $N_i$ be a filtered system of $G$-modules in $\C'(M)$ and $N$ its colimit. Then there is a system of free $G$-modules $F_i$ with basis $B_i$ and compatible 
epimorphisms $F_i\fl N_i$ with kernel $K_i$ such that for each $i\rightarrow j$ the corresponding morphism $F_i\fl F_j$ is a monomorphism sending $B_i$ into $B_j$.
Let $F$ and $K$ be the colimits of the $F_i$'s and the $K_i$'s.
For each $i$, $M\otimes K_i$ is regular  and, because the class of regular right $A[G]$-modules is cocomplete the module:
$$\build\lim_\fl^{} M\otimes K_i\simeq M\otimes \build\lim_\fl^{} F_i=M\otimes K$$
is regular. On the other hand we have an exact sequence:
$$0\fl K\fl F\fl N\fl 0$$
where $F$ is free. Therefore, $M\otimes K$ is regular and $N$ belongs to $\C'(M)$. Then the class $\C'(M)$ is cocomplete.

Let $0\fl N_1\fl N_2\fl N_3\fl 0$ be an exact sequence of right $G$-modules. Then there is a commutative diagram of right $G$-modules with exact lines:
$$\diagram{0&\hfl{}{}&F_1&\hfl{}{}&F_2&\hfl{}{}&F_3&\hfl{}{}&0\cr &&\vfl{}{}&&\vfl{}{}&&\vfl{}{}&&\cr0&\hfl{}{}&N_1&\hfl{}{}&N_2&\hfl{}{}&N_3&\hfl{}{}&0\cr}$$
where $F_1$, $F_2$ and $F_3$ are free, and vertical morphisms are epimorphisms with kernel $K_1$, $K_2$ and $K_3$.

Suppose two of the modules $N_1,N_2,N_3$ are in $\C'(M)$. Then two of the modules $M\otimes K_1$, $M\otimes K_2$ and $M\otimes K_3$ are regular. Therefore, the third
module $M\otimes K_i$ is regular and $K_i$ belongs to $\C'(M)$. Hence, $\C'(M)$ is exact and, because $G$ is regular $\C'(M)$ is the class of all $G$ modules.

Let $K$ be the kernel of the augmentation morphism $\Z[G]\fl \Z$. Since $\Z$ belongs to $\C'(M)$, $M\otimes K$ is regular. But we have an exact sequence:
$$0\fl M\otimes K\fl M\otimes\Z[G]\fl M\otimes\Z\fl 0$$
where $M\otimes K$ and $M\otimes\Z[G]$ are regular. Then $M\otimes\Z\simeq M$ is regular. and $A[G]$ is right-regular.\cqfd
 \vskip 12pt
\noi{\bf 3.8 Proposition:} {\sl We have the following properties:

1) A group with finite homological dimension is regular

2) If a group $G$ is regular, every subgroup of $G$ is regular

3) Let $G_1$, $G_2$ and $H$ be three regular groups and $\alpha_1:H\fl G_1$ and $\alpha_2:H\fl G_2$ be two injective groups homomorphisms. Then the corresponding 
amalgamated  free product is regular

4) Let $H$ and $G$ be two regular groups and $\alpha:H\fl G$ and $\beta:H\fl G$ be two injective group homomorphisms. Then the corresponding HNN extension is regular.

5) Let $G_i$ be a filtered system of regular groups where each morphism $G_i\fl G_j$ is injective. Then the colimit of this system is regular

6) If a group $G$ is an extension of two regular groups, then $G$ is regular.}
\vskip 12pt
\noi{\bf Proof:} Let $G$ be a group with finite homological dimension. Then there is an integer $n>2$ such that $H_n(G,M)=0$ for all $G$-module $M$. Let $C$ be a
projective resolution of $\Z$. Then we have an exact sequence:
$$\dots\fl C_n\fl C_{n-1}\fl\dots\fl C_1\fl C_2\fl \Z\fl 0$$
where each $C_i$ a projective $G$-module. Let $Z$ be the kernel of $d:C_{n-2}\fl C_{n-3}$. We have the exact sequence:
$$\dots\fl C_{n+1}\fl C_n\fl C_{n-1}\fl Z\fl 0$$

Let $M$ be a $G$-module. We have:
$$0=H_n(G,M)\simeq H_n(C,M)\simeq\hbox{Tor}_1^{\Z[G]}(Z,M)$$
and Tor$_1^{\Z[G]}(Z,M)=0$ for every $G$-module $M$. Hence, $Z$ is flat and then regular. Thus we have an exact sequence:
$$0\fl Z\fl C_{n-2}\fl C_{n-3}\fl\dots\fl C_1\fl C_0\fl \Z\fl 0$$
where $Z$ and the $C_i$'s are all regular. Hence, $Z$ is regular and, because of Lemma 3.7, $G$ is regular

Properties 2), 3), 4) and 5) are direct consequences of Lemmas 3.3, 3.4, 3.5 and 3.2.

Let $G$ be a group and $H$ be a normal subgroup of $G$. Let $\Gamma$ be the quotient $G/H$ and suppose that $H$ and $\Gamma$ are regular.

Let $\A_0\v$, $\A_1\v$ and $\A_2\v$ be the categories of right modules over $G$, $H$ and $\Gamma$ respectively.

For $i=0,1,2$, the category of finitely generated projective modules in $\A_i\v$ will be denoted $\A_i$.

Let $\C_0$ be the class of $G$-regular modules in $\A_0\v$.

Let $F_1:\A_1\v\fl \A_0\v$ be the functor: $M\mapsto M\build\otimes_H^{}\Z[G]$ and $F_2:\A_2\v\fl \A_0\v$ be the functor sending each $\A_2\v$-module $M$ to $M$ equipped
with the induced $G$-action.

It is easy to see that $F_1$ and $F_2$ are exact functors commuting with filtered colimit.

Let $\C_1$ be the class of $\A_1\v$-modules $M$ such that $F_1(M)\in\C_0$. This class is exact and cocomplete and contains $\A_1$. Since $H$ is regular, we have: 
$\C_1=\A_1\v$ and $\C_1$ contains $\Z$ equipped with the trivial $H$-action. Hence, $\C_0$ contains the module $\Z[\Gamma]$.

Let $\C_2$ be the class of $\A_2\v$-modules $M$ such that $F_2(M)\in\C_0$. This class is exact and cocomplete and contains $\Z[\Gamma]$. Since $\Gamma$ is regular, $\C_2$
contains $\A_2$ and we have: $\C_2=\A_2\v$. In particular $\C_2$ contains $\Z$ equipped with the trivial $\Gamma$-action. 

Therefore, $\C_0$ contains $\Z$ equipped with the trivial $G$-action and, because of Lemma 3.7, $G$ is regular.\cqfd

\vskip 12pt
Let $\G$ be the category of groups and injective morphisms of groups and \Cl be the smallest class of groups satisfying the following properties:

$\bullet$ \Cl contains the trivial group

$\bullet$ \Cl is stable under taking filtered colimit (in $\G$) 
 
$\bullet$ \Cl is stable under taking amalgamated free product and HNN extensions (in $\G$)

This class is defined in the same way that Waldhausen's class Cl but without any condition of regular coherency. Therefore, Cl is contained in the class \Cl.
\vskip 12pt
\noi{\bf 3.9 Lemma:} {\sl The class \Cl is stable under taking subgroup and extension. Moreover each group in \Cl is regular.}
\vskip 12pt
\noi{\bf Proof:} Because of Lemma 3.8, the class \Cl is contained in the class of regular groups and each group in \Cl is regular.

Let $\C(0)$ be the class reduced to the trivial group and for each ordinal $a>0$, we denote $\C(a)$ the class of groups obtained by amalgamated free 
product, HNN extension or filtered colimit with groups that which are each in some $\C(b)$ with $b<a$. Then the class \Cl is the union of all these classes $\C(a)$.

If $G$ is a group in \Cl, the smallest ordinal $a$ such that $G$ is in $\C(a)$ is denoted $|G|$.

The class \Cl is stable under taking subgroup. Waldhausen's proof (Proposition 19.3 in [W1]) applies exactly in our situation.

For each ordinal $a$ we note $P(a)$ the following property:

$\bullet$ For each exact sequence of groups $1\fl H\fl\Gamma\build\fl_{}^f G\fl 1$ with $H\in$\Cl and $G$ in $\C(a)$, the group $\Gamma$ belongs to \Cl.

The property $P(0)$ is obviously true.

Let $a>0$ be an ordinal such that $P(b)$ is true for all $b<a$  and $1\fl H\fl\Gamma\build\fl_{}^f G\fl 1$ be an exact sequence of groups with $H\in$\Cl and $G$ in 
$\C(a)$. We have three possibilities $G$ is a filtered colimit, an amalgamated free product or an HNN extension.

In the first case, there is a filtered system of subgroups $G_i$ of $G$ such that: $|G_i|<a$ and $G=\build\lim_\rightarrow^{} G_i$. Then $f^{-1}(G)$ is the colimit
of the filtered system of subgroups $f^{-1}(G_i)$. By induction each group $f^{-1}(G_i)$ is in \Cl and its colimit $\Gamma$ belongs to \Cl.

In the second case there is two subgroups $G_1$ and $G_2$ of $G$ which intersect in a subgroup $H$ such that $|H|<a$, $|G_1|<a$, $|G_2|<a$ and $G$ is the amalgamated
free product $G_1\build*_H^{} G_2$. Set: $G'_1=f^{-1}(G_1)$, $G'_2=f^{-1}(G_2$ and $H'=f^{-1}(H)$. By induction $H'$, $G'_1$ and $G'_2$ are in \Cl and $\Gamma$ is the 
amalgamated free product $G'_1\build*_{H'}^{} G'_2$. Then $\Gamma$ belongs to the class \Cl.

In the third case there is a subgroup $G_1$ of $G$, a subgroup $H$ of $G_1$ and an element $t\in G$ such that $|G|<a$, $|H|<a$ and $\Gamma$ is the HNN extension of $G_1$
via the inclusion $H\subset G_1$ and the morphism $\beta:H\fl G_1$ sending each $x\in H$ to $txt^{-1}$. In particular we have: $tHt^{-1}\subset G_1$.

Let $\theta$ be an element in $f^{-1}(t)$ and set: $G'_1=f^{-1}(G_1)$ and $H'=f^{-1}(H)$. It is easy to see that $\theta H'\theta^{-1}\subset G'_1$ and $H$ is the HNN
extension of $G'_1$ via the inclusion $H'\subset G'_1$ and the morphism $x\mapsto \theta x\theta^{-1}$ from $H'$ to $G'_1$. Hence, $\Gamma$ belongs to \Cl. 

Then properties $P(a)$ are all true and \Cl is stable under extension.\cqfd
\vskip 12pt
\noi{\bf Proof of Theorem 3:}
\vskip 12pt
Let $A$ be a ring and $G$ be a group. We have a map of $\Omega$-spectra $H(G,\underline K(A))\fl \underline K(A[G])$ called the assembly map where $H(BG,\underline K(A))$
is the homology spectrum with coefficients in $\underline K(A))$.  Thus we have, for each integer $n$:
$$\pi_n(H(BG,\underline K(A)))=H_n(BG,\underline K(A))$$
The loop space of the homotopy fiber of the assembly map is called the Whitehead spectrum of $G$ with coefficients in $A$. Then we get a homotopy fibration of 
$\Omega$-spectra: $$H(BG,\underline K(A))\fl \underline K(A[G])\fl \underline Wh^A(G)$$

It is proven in [V] the following: if a group $\Gamma$ is the amalgamated free product $G_1\build*_H^{} G_2$ via injective morphisms $\alpha_1:H\fl G_1$ and 
$\alpha_2:H\fl G_2$, there is a left-flat $A[H]\times A[H]$-bimodule $S$, a homotopy cartesian square of spectra:
$$\diagram{\underline K(A[H])&\hfl{\alpha_1}{}&\underline K(A[G_1])\cr\vfl{\alpha_2}{}&&\vfl{}{}\cr \underline K(A[G_2])&\hfl{}{}&\underline K(A[\Gamma])'\cr}$$
and a homotopy equivalence: $\underline K(A[\Gamma])\simeq \underline K(A[\Gamma])'\oplus \Omega^{-1}\underline Nil(A[H]\times A[H],S)$.

But we have also a homotopy cartesian square of spectra:
$$\diagram{H(\pi,\underline K(A))&\hfl{\alpha_1}{}&H(G_1,\underline K(A))\cr\vfl{\alpha_2}{}&&\vfl{}{}\cr 
H(G_2,\underline K(A))&\hfl{}{}&H(\Gamma,\underline K(A))\cr}$$

The relative part of these two homotopy cartesian squares of spectrum is the following  homotopy cartesian square of spectra:
$$\diagram{\underline Wh^A(H)&\hfl{\alpha_1}{}&\underline Wh^A(G_1)\cr\vfl{\alpha_2}{}&&\vfl{}{}\cr\underline Wh^A(G_2)&\hfl{}{}&\underline Wh^A(\Gamma)'\cr}$$
and a homotopy equivalence: $\underline Wh^A(\Gamma)\simeq \underline Wh^A(\Gamma)'\oplus \Omega^{-1}\underline Nil(A[H]\times A[H],S)$.

Then, if $H$ is regular and $A$ is right-regular, Lemma 3.7 implies that $A[H]$ and $A[H]\times A[H]$ are right-regular and Theorem 2 implies that 
$\underline Nil(A[H]\times A[H],S)$ is contractible. Hence, we get a homotopy cartesian square:
$$\diagram{\underline Wh^A(H)&\hfl{\alpha_1}{}&\underline Wh^A(G_1)\cr\vfl{\alpha_2}{}&&\vfl{}{}\cr\underline Wh^A(G_2)&\hfl{}{}&\underline Wh^A(\Gamma)\cr}$$

If $\Gamma$ is the HNN extension of a group $G$ via two injective groups homomorphisms $\alpha: H\fl G$ and $\beta:H\fl G$, with exactly the same method we get a homotopy
fibration of spectra:
$$\underline Wh^A(H)\build\fl_{}^{\alpha-\beta}\underline Wh^A(G)\fl \underline Wh^A(\Gamma)'$$
and a homotopy equivalence: $\underline Wh^A(\Gamma)\simeq \underline Wh^A(\Gamma)'\oplus \Omega^{-1}\underline Nil(A[\pi]\times A[\pi],S)$ where $S$ is some 
left-flat $A[H]\times A[H]$-bimodule. Therefore, if $A$ is right-regular and $H$ is regular we have a homotopy fibration:
$$\underline Wh^A(H)\build\fl_{}^{\alpha-\beta}\underline Wh^A(G)\fl \underline Wh^A(\Gamma)$$
It is then easy to prove by induction that, for every $G\in$\Cl, the spectrum $\underline Wh^A(G)$ is contractible and Theorem 3 is proven.\cqfd
\vskip 24pt
\noi{\bf 4. The Vanishing Theorem.}
\vskip 12pt
Let $\Ex$ be the class of essentially small exact categories and $\WE$ be the class of pairs $(\A,w)$ where $\A$ is an essentially small exact category and $w$ is a 
class of $\A$-morphisms with the following properties:

$\bullet$ $w$ contains the isomorphisms

$\bullet$ if $f$ and $g$ are two composable $\A$-morphisms and two of the morphisms $f,g,g\circ f$ are in w, so is the third

$\bullet$ Suppose we have a commutative diagram:
$$\diagram{A&\hfl{}{}&B&\hfl{}{}&C\cr\vfl{f}{}&&\vfl{g}{}&&\vfl{h}{}\cr A'&\hfl{}{}&B'&\hfl{}{}&C'\cr}$$
where lines are conflations. Then if two of the morphisms $f,g,h$ are in $w$ so is the third.

For $(\A,w)$ in $\WE$, $w$ will be called the class of equivalences and $\A$ equipped with the inflations and the equivalences is a Waldhausen category satisfying  
saturation and extension axioms (see [W2]). 

Such a pair will be called a Waldhausen exact category. 

Denote also $\WE^h$ the class of Waldhausen exact categories having a cylinder functor satisfying the cylinder axiom that is: there is a functor $T$ from the category of 
morphisms of the Waldhausen exact category $\A$  to $\A$ and, for every $\A$-morphism $f:A\fl B$ a natural sequence: $A\oplus B\build\fl_{}^\alpha T(f)\build\fl_{}^\beta
B$ where $\alpha$ is an inflation and $\beta$ is an equivalence such that $\beta\alpha$ is the identity on $B$ and the morphism $f$ on $A$.

An exact functor $F:(\A,w)\fl(\B,w)$ between two Waldhausen exact categories is an exact functor: $\A\fl \B$ that respects equivalences.
A subcategory of a Waldhausen exact category $\A$ is a subcategory of the category $\A$ where the Waldhausen structure is induced by the Waldhausen structure on $\A$.

Notice that the class of exact functors between categories in $\Ex$, $\WE$ or $\WE^h$ is not a set. Thus $\Ex$, $\WE$ and $\WE^h$ equipped with exact functors are not 
really categories but these classes, equipped with isomorphism classes of exact functors are categories (still denoted $\Ex$, $\WE$ and $\WE^h$) and by taking 
isomorphisms as equivalences, we get an inclusion of categories: $\Ex\subset\WE$ compatible with the K-theory functor. Moreover the inclusion $\A\subset \A_*$ induces an 
functor $\varepsilon:\Ex\fl \WE^h$ also compatible with K-theory (because of Theorem 1.17).
\vskip 12pt
Throughout this section, $A$ will be a right-regular ring and $S$ a left-flat $A$-bimodule. The category of finitely generated projective right $A$-modules will be denoted 
$\A$. The category of all right $A$-modules is a strict cocompletion of $\A$ and will be denoted $\A\v$. 

Let $M$ be an $\A\v$-module and $\theta:M\fl M\otimes S$ be an $\A\v$-morphism. We say that $\theta$ is nilpotent if, for every $x\in M$ there is an integer $n>0$ such 
that $x$ is killed by the morphism $\theta^n:M\fl X\otimes S^{\otimes n}$. 

Let $\Nil(A,S)\v$ be the class of pairs $(M,\theta)$ where $M$ is an $\A\v$-module and $\theta:M\fl M\otimes S$ is a nilpotent morphism. A morphism between two such pairs 
$(M,\theta)$ and $(N,\theta)$ is a morphism $\varphi:M\fl N$ commuting with $\theta$. Then with these morphisms $\Nil(A,S)\v$ is a Grothendieck category and the 
correspondence $(M,\theta)\mapsto M$ is an exact functor $F:\Nil(A,S)\v\fl\A\v$.  

If $n>0$ is an integer, we have the following full subcategories of $\Nil(A,S)\v$:

$\bullet$ $\Nil(A,S,n)\v$ is the category of $X\in\Nil(A,S)\v$ such that: $\theta^n X=0$

$\bullet$ $\Nil(A,S)$ is the category of $X\in\Nil(A,S)\v$ such that $F(X)\in\A$

$\bullet$ $\Nil(A,S,n)$ is the category $\Nil(A,S)\cap \Nil(A,S)\v$

Notice that $\Nil(A,S)$ is a fully exact subcategory of $\Nil(A,S)\v$ and $\Nil(A,S,n)$ is a fully exact subcategory of $\Nil(A,S,n)\v$. Moreover $F$ induces exact
functors $F:\Nil(A,S)\fl\A$ and $F:\Nil(A,S,n)\fl\A$.

The functor $F:(M,\theta)\mapsto M$ has a section: $M\mapsto (M,0)$. Then there are well-defined infinite loop spaces $Nil(A,S)$, $Nil(A,S,n)$ and homotopy equivalences:
$$K(\Nil(A,S))\simeq K(A)\times Nil(A,S)\hskip 24pt K(\Nil(A,S,n))\simeq K(A)\times Nil(A,S,n)$$

The $i$-th homotopy group of $Nil(A,S)$ and $Nil(A,S,n)$  will be denoted $Nil_i(A,S)$ and $Nil_i(A,S,n)$.
\vskip 12pt

In order to prove the Vanishing Theorem about $\Nil(A,S)$, we will need to consider more complicated categories.

Consider a bigraded differential algebra $\Lambda$ with a differential of bidegree $(-1,0)$. If an element $x\in\Lambda$ has bidegree $\partial^\circ(x)=(a,b)$ then 
$a=\partial_1^\circ(x)$ will be called the degree of $x$ and $b=\partial_2^\circ(x)$ its weight. 

Let $\R$ be the category of such bigraded differential algebras for which every homogeneous element has a non-negative weight.

Let $\Lambda$ be an algebra in $\R$. There is an involution on $\Lambda$ sending each element $x\in\Lambda$ of degree $a$ to $\overline x=(-1)^a x$. Therefore, 
for all $x,y$ in $\Lambda$ we have:
$$d(xy)=d(x)y+\overline x d(y)\ \ \ \hbox{and}\ \ \ \overline{d(x)}=-d(\overline x)$$

If $p$ and $q$ are two integers the set of $x\in\Lambda$ such that: $\partial^\circ(x)=(p,q)$ will be denoted $\Lambda[p,q]$ and the two-sided ideal of $\Lambda$
generated by homogeneous elements in $\Lambda$ of weight $\geq q$ will be denoted $J_q(\Lambda)$.

Let $X$ be an $\A\v$-complex. A $\Lambda$-action on $X$ is a family of $A$-linear maps $\Lambda[p,q]\build\otimes_\Z^{} X\fl X\build\otimes_A^{}S^{\otimes q}$ of 
degree $p$ denoted $a\otimes x\mapsto ax$ with the following properties:
$$\forall x\in X,\ \ \ 1.x=x$$
$$\forall p,q,r,s\geq0,\ \forall a\in\Lambda[p,q],\ \ \forall b\in\Lambda[r,s],\ \ \forall x\in X,\ \ \ a(bx)=(ab)x\in X\otimes S^{\otimes(q+s)}$$ 
$$\forall a\in\Lambda,\ \ \forall x\in X,\ \ \ d(ax)=d(a)x+\overline a d(x))$$
Such a $\Lambda$-action on $X$ is said to be nilpotent if, for every $x\in X$, there is an integer $n$ such that $J_n(\Lambda)x=0$.

The category of finite $\A$-complexes equipped with a nilpotent $\Lambda$-action will be denoted $\Nil(\Lambda)$ and the category of $\A\v$-complexes equipped with a 
nilpotent $\Lambda$-action will be denoted $\Nil(\Lambda)\v$.

By forgetting the $\Lambda$-action, we get a forgetful functor
$F:\Nil(\Lambda)\fl\A_*$ and we have a unique structure of exact category on $\Nil(\Lambda)$ such that $(X\fl Y\fl Z)$ is a $\Nil(\Lambda)$-conflation if and only if
$F(X)\fl F(Y)\fl F(Z)$ is an $\A_*$-conflation. Hence, $\Nil(\Lambda)$ is an essentially small exact category and the homology equivalences induce on $\Nil(\Lambda)$ a 
structure of Waldhausen exact category. Thus $\Nil(\Lambda)$ belongs to $\WE$. 

\vskip 12pt
An algebra $\Lambda\in\R$ is said to be free if there is a homogeneous set $U\subset\Lambda$ (called a basis of $\Lambda$) such that $\Lambda$ is freely generated by $U$ 
(as a ring). We denote $\R_0$ the category of algebras $\Lambda\in\R$ having a basis $U$ such that:

$\bullet$ for all $x\in U$, we have: $\partial_2^\circ(x)>0$

$\bullet$ for any integer $q$, the set of elements of $U$ of weight $q$ is finite.
\vskip 12pt
Let $n>0$ be an integer. We have in $\Nil(\Lambda)\v$ a subcategory $\Nil(\Lambda,n)\v$ defined by:
$$X\in\Nil(\Lambda,n)\v\ \Longleftrightarrow\ \ J_n(\Lambda)X=0$$
We set also:
$$\Nil(\Lambda,n)=\Nil(\Lambda,n)\v\cap\Nil(\Lambda)$$

The category $\Nil(\Lambda,n)$ is a fully exact subcategory of $\Nil(\Lambda)$ and we have a filtration of $\Nil(\Lambda)$:
$$0=\Nil(\Lambda,0)\subset\Nil(\Lambda,1)\subset\dots$$

Let $X$ be a $\Nil(\Lambda)$-module. A filtration: $0=X_0\subset X_1\subset X_2\subset\dots\subset X_n=X$ is said to be admissible if each quotient $X_i/X_{i-1}$ belongs 
to $\Nil(\Lambda,1)$ and a $\Nil(\Lambda)$-module is said to be admissible if it has an admissible filtration. Admissible $\Nil(\Lambda)$-modules are the objects of a 
subcategory $\Nil(\Lambda)'$ of $\Nil(\Lambda)$. We have also subcategories $\Nil(\Lambda,n)'$ of $\Nil(\Lambda)$ defined by:
$$\Nil(\Lambda,n)'=\Nil(\Lambda,n)\cap\Nil(\Lambda)'$$

Categories $\Nil(\Lambda)'$ and $\Nil(\Lambda,n)'$ are Waldhausen exact categories and, as exact categories, they are fully exact in $\Nil(\Lambda)$.
\vskip 12pt
Since the inclusion $\Z\subset\Lambda$ has a retraction $\varepsilon:\Lambda\fl\Z$ sending each element of a basis of $\Lambda$ to zero, there are infinite loop spaces 
$Nil(\Lambda)$, $Nil(\Lambda)'$, $Nil(\Lambda,n)$, $Nil(\Lambda,n)'$ (for $n>0$) and homotopy equivalences:
$$K(\Nil(\Lambda))\simeq K(A)\times Nil(\Lambda)\hskip 24pt K(\Nil(\Lambda)')\simeq K(A)\times Nil(\Lambda)'$$
$$K(\Nil(\Lambda,n))\simeq K(A)\times Nil(\Lambda,n)\hskip 24pt K(\Nil(\Lambda,n)')\simeq K(A)\times Nil(\Lambda,n)'$$

\vskip 12pt
\noi{\bf 4.1 Lemma:} {\sl Rings $\Z$ and $\Z[\theta]$ with $\partial^\circ(\theta)=(0,1)$ are algebras in $\R_0$.

Moreover $\Nil(\Z)$ is equivalent to the category $\A_*$ of finite $\A$-complexes, $\Nil(\Z[\theta])$ is equivalent to the category $\Nil(A,S)_*$ of finite 
$\Nil(A,S)$-complexes and $\Nil(\Z[\theta],n)$ is equivalent to the category $\Nil(A,S,n)_*$ of finite $\Nil(A,S,n)$-complexes.}
\vskip 12pt
\noi{\bf Proof:} Rings $\Z$ and $\Z[\theta]$ are obviously in $\R_0$ (with basis $\emptyset$ and $\{\theta\}$) and we have: $d(\theta)=0$. If $\Lambda=\Z$, a finite 
$(\Lambda,\A,S)$-complex is just a finite $\A$-complex $C$ and such a complex is always nilpotent. Then we have: $\Nil(\Z)\simeq\A_*$.

Suppose $\Lambda=\Z[\theta]$ with $\partial^\circ(\theta)=(0,1)$. Then a finite $(\Lambda,\A,S)$-complex is just a finite $\A$-complex equipped with a morphism
$\theta:C\fl C\otimes S$ of degree $0$. If this action is nilpotent we have:
$$\forall x\in C, \ \exists q>0,\ \hbox{such that}\ \ \theta^q(x)=0$$ 
Therefore, the action is nilpotent if and only if the action of $\theta$ is nilpotent in each degree and a $\Nil(\Lambda))$-module is a finite $\Nil(A,S)$-complex.

The same holds with $\Nil(\Lambda,n)$-modules.\cqfd
\vskip 12pt

Let $X$ be a $\Nil(\Lambda)\v$-module and $X_*=(0=X_0\subset X_1\subset X_2\subset\dots)$ be a filtration of $X$ by $\Nil(\Lambda)\v$-submodules. Then any element in $X$
belongs to some $X_k$. This filtration is said to be special if, for each $i>0$ and each homogeneous element $\theta\in\Lambda$ of weight $q$, we have: 
$$\theta X_i\subset X_{i-q}\otimes S^{\otimes q}$$
(with $X_i=0$ for $i\leq 0$).

The category of $X\in\Nil(\Lambda)$-modules equipped with a special admissible filtration will be denoted $\F(\Lambda)$.

If $n>0$ is an integer, the category of $\Nil(\Lambda)$-modules equipped with a special admissible filtration of length $n$ will be denoted $\F_n(\Lambda)$. The category
$\F_n(\Lambda)$ is an exact category. For any  $\F_n(\Lambda)$-module $X$, its filtration will be denoted $0=F_0(X)\subset F_1(X)\subset\dots\subset F_n(X)$. Then
for any $i\leq n$, the correspondence $X\mapsto F_i(X)$ is an exact functor from $\F_n(\Lambda$ to $\Nil(\Lambda,i)'$. We have also exact functors $X\mapsto \Phi_i(X)$ 
sending each $\F_n(\Lambda)$-module $X$ to $F_i(X)$ equipped with the filtration $0=F_0(X)\subset\dots\subset F_i(X)$. The functor $\Phi_i$ is an exact functor from 
$\F_n(\Lambda)$ to $\F_i(\Lambda)$. 

Consider a sequence of $\A\v$-complexes $X_i$, with $i\leq n$ and $X_i=0$ for any $i\leq 0$, stabilization morphisms $\gamma:X_{i-1}\fl X_i$ and a $\Lambda$-action on 
$\build\oplus_i^{}X_i$ such that:

$\bullet$ for any $\theta\in\Lambda$ of weight $p$, and any $i\leq n$ we have: $\theta X_i\subset X_{i-p}\otimes S^{\otimes p}$ and the following diagram is commutative:
$$\diagram{X_{i-1}&\hfl{\theta}{}&X_{i-1-p}\otimes S^{\otimes p}\cr \vfl{\gamma}{}&&\vfl{\gamma}{}\cr X_i&\hfl{\theta}{}&X_{i-p}\otimes S^{\otimes p}\cr}$$

The category of such data will be denoted $\F_n(\Lambda)\v$. The correspondence $X\mapsto X_i$ is a functor $F_i:\F_n(\Lambda)\v\fl \Nil(\Lambda,i)\v$ and the
 correspondence sending $X$ to the sequence $X_j$, $j\leq i$ is a functor $\Phi_i:\F_n(\Lambda)\v\fl \F_i(\Lambda)\v$.

It can be verified that $\F_n(\Lambda)\v$ is a strict cocompletion of $\F_n(\Lambda)$ and that functors $F_i:\F_n(\Lambda)\v\fl \Nil(\Lambda,i)\v$ and 
$\Phi_i:\F_n(\Lambda)\v\fl\F_i(\Lambda)\v$ are induced by functors $F_i:\F_n(\Lambda)\fl \Nil(\Lambda,i)'$ and $\Phi_i:\F_n(\Lambda)\fl\F_i(\Lambda)$.

The subcategory of such data in $\F_n(\Lambda)\v$ such that each stabilization morphism $\gamma:X_{i-1}\fl X_i$ is an inclusion will be denoted $\F_n(\Lambda)\bv$\ . This
category is the category of $\Nil(\Lambda,n)\v$-modules equipped with a special filtration of length $n$.

Let $X$ be a $\Nil(\Lambda,n)\v$-module. For any integer $p\geq0$, the set of elements $x\in X$ such that $J_p(\Lambda)x=0$ will be denoted $T_p(X)=X_p$. Thus we get a 
special filtration of $X$: $0=X_0\subset X_1\subset X_2\subset\dots X_n=X$ and we have a functor $T:\Nil(\Lambda,n)\v\fl \F(\Lambda,n)\bv$ witch is right 
adjoint to $F_n:\F_n(\Lambda)\bv\fl \Nil(\Lambda,n)\v$. Notice that, for any $\Nil(\Lambda,n)\v$-module $X$, we have: $T_p(X)=F_p(T(X))$.

\vskip 12pt
\noi{\bf 4.2 Lemma:} {\sl Let $\Lambda$ be an algebra in $\R_0$. Then categories $\Nil(\Lambda)$ and $\F_n(\Lambda)$ have $k$-cone functors and cylinder functor satisfying
the cylinder axiom and these functors respect subcategories $\Nil(\Lambda)'$, $\Nil(\Lambda,n)$ and $\Nil(\Lambda,n)'$ of $\Nil(\Lambda)$. Hence each of these categories 
belongs to $\WE^h$.}
\vskip 12pt
\noi{\bf Proof:}
Let $X$ and $X'$ be two $\Nil(\Lambda)$-modules and $X^\circ$ and $X'^\circ$ be their underlying graded $\A$-modules. For each integer $k$, we denote $\Map(X,X')_k$
the set of $A$-linear morphisms $f:X^\circ\fl X'^\circ$ of degree $k$ such that, for any $\theta\in\Lambda$ of degree $a$ and any $x\in X$ we have:
$$f(\theta x)=(-1)^{ak}\theta f(x)$$

These sets $\Map(X,X')_k$ define  a graded $\Z$-module $\Map(X,X')$ and this graded module has a differential:
$$\forall k\in\Z, \forall f\in\Map(X,X')_k,\ \ d(f)=d\circ f-(-1)^k f\circ d$$

If $k$ is an integer, the elements of $\Map(X,X')_k$ are called $k$-maps (or maps of degree $k$). A $k$-map $f$ with $d(f)=0$ is called a $k$-morphism (or a morphism of 
degree $k$). Thus a morphism between two $\Nil(\Lambda)$-modules is a $0$-morphism.

Let $a$ be an integer and $\varphi:X\fl X'$ be an $a$-morphism of $\A\v$-complexes. Then we have morphisms of complexes $i:X'\fl X\oplus X'$, $s:X\fl X\oplus X'$, 
$p:X\oplus X'\fl X$ and $r:X\oplus X'\fl X'$. We have also the following relations:
$$pi=0\hskip 24pt rs=0\hskip 24pt ri=1\hskip 24pt ps=1\hskip 24pt ir+sp=1$$

If $k$ is an integer, there is a unique way to modify the degree, the $\Lambda$-action and the differential on $X\oplus X'$ in such a way that
$i$, $r$, $p$ and $s$ are $k$-map, $(-k)$-map, $(-1-k-a)$-map and $(1+k+a)$-map respectively and:
$$d(i)=0\hskip 24pt d(p)=0\hskip 24pt d(s)=i\varphi\hskip 24pt d(r)=-(-1)^k \varphi p$$

This modified complex $X''$ will be called the $k$-cone of $\varphi$ and denoted $C_k(\varphi)$. Thus we have an exact sequence:
$$0\fl X'\build\fl_{}^i X''\build\fl_{}^p X\fl 0$$
where $i$ is a $k$-morphism, $p$ is a $(-1-k-a)$-morphism, $r$ is a $-k$-map and a retraction of $i$ and $s$ is a $(k+a+1)$-map and a section of $p$. 

If $X$ and $X'$ are equipped with admissible special filtrations and $\varphi$ respects these filtrations, all these constructions respect the filtrations and $k$-cones
functors are well-defined in $\F_n(\Lambda)$.

If $f:X\fl X'$ is a morphism (of degree $0$), we have two exact sequences:
$$0\fl X'\build\fl_{}^i C_0(f)\build\fl_{}^p X\fl 0$$
$$0\fl X\build\fl_{}^j C_0(p)\build\fl_{}^q C_0(f)\fl 0$$
Let $\sigma$ be the standard section of $q$. Thus we have: $d(\sigma)=jp$ and $d(\sigma i)=0$. Then $\sigma i:X'\fl C_0(p)$ is a morphism (of degree $0$) and we have an 
exact sequence:
$$0\fl X'\build\fl_{}^{\sigma i}C_0(p)\fl C_0(\hbox{Id}_X)\fl 0$$
Therefore, $T(f)=C_0(p)$ is a cylinder functor which satisfies the cylinder axiom.
 
Moreover all these constructions respect subcategories $\Nil(\Lambda)'$, $\Nil(\Lambda,n)$ and $\Nil(\Lambda,n)'$ of $\Nil(\Lambda)$.\cqfd
\vskip 12pt
Forgetful functors $\F_n(\Lambda)\v\fl\F_n(\Z)\v$ and $\Nil(\Lambda,n)\v\fl \Nil(\Z,n)\v$ will be denoted $X\mapsto X^\circ$.
\vskip 12pt
\noi{\bf 4.3 Lemma:} {\sl Let $\Lambda$ be an algebra in $\R_0$ and $n>0$ be an integer. Let $X,Y,Y'$ be three $\F_n(\Lambda)\bv$-modules, $f:X\fl Y$ and $\beta:Y\fl Y'$ 
be two $\F_n(\Lambda)\bv$-morphisms, $X'$ be an $\F_n(\Z)\bv$-module and
$$\diagram{X^\circ&\hfl{\alpha}{}&X'\cr\vfl{f}{}&&\vfl{f'}{}\cr Y^\circ&\hfl{\beta}{}&Y'^\circ\cr}\leqno{(D)}$$
be a commutative diagram in $\F_n(\Z)\bv$ such that:

1) for any $k<n$, $f:F_k(X)\fl F_k(Y)$ is a surjective homology equivalence 

2) $\beta:F_n(Y)\fl F_n(Y')$ is monic and $\theta F_n(Y')$ is contained in $F_{n-p}(Y)\otimes S^{\otimes p}$ for any $\theta\in\Lambda$ of weight $p>0$

3) $\alpha:X^\circ\fl X'$ is monic and its cokernel is projective and bounded from below.

Then there is a $\Lambda$-action on $X'$ compatible with $\alpha$ and $f'$ such that $\theta F_n(X')$ is contained in $F_{n-p}(X)\otimes S^{\otimes p}$ for any 
$\theta\in\Lambda$ of weight $p>0$.}
\vskip 12pt
\noi{\bf Proof:} For any complex $C$ the graded module of cycles of $C$ will be denoted $C^\bullet$.

Suppose $\alpha$ induces an isomorphism $F_{n-1}(X)\build\fl_{}^\sim F_{n-1}(X')$. Then we have an exact sequence:
$$0\fl F_n(X)\fl F_n(X')\fl K\fl 0$$
where $K$ is a projective $A$-complex. Let $s:K\fl F_n(X')$ be a section of $F_n(X')\fl K$. The map $s$ belongs to $\Map(K,F_n(X'))_0$ and $d(s)$ is a 
$(-1)$-morphism from $K$ to $F_n(X)$. 

Let $U$ be a basis of $\Lambda$. We put on $U$ a total order such that we have for any $x,y\in U$ with; $\partial^\circ(x)=(a,b)$ and $\partial^{\circ}(y)=(a',b')$:
$$x<y\Longrightarrow b<b'\ \ \hbox{or}\ \ \bigl(b=b'\ \ \hbox{and}\ \  a\leq a'\bigr)$$
For any $x\in U$ the subalgebra of $\Lambda$ generated by the elements $y$ in $U$ with $y<x$ will be denoted $\Lambda(<x)$.

Let $\theta$ be an element of $U$ of bidegree $(a,b)$ such that the action of any $x<\theta$ is defined on $X'$. Then $\Lambda(<\theta)$ acts on $X'$. In order to extend 
the $\Lambda(<\theta)$-action to a $\Lambda(\leq \theta)$-action, we just have to define a map $\theta s$ in $\Map(K,F_{n-b}(X))_a\otimes S^{\otimes b}$.

Because of property 1), we have an exact sequence:
$$F_{n-b}(X)\build\fl_{}^{d\oplus f} F_{n-b}(X)^\bullet\oplus F_{n-b}(Y)\build\fl_{}^{-f\oplus d}F_{n-b}(Y)^\bullet\fl 0$$
and, because $S$ is left-flat and $K$ is projective, we have an exact sequence:
$$\Map(K,F_{n-b}(X))_a\otimes S^{\otimes b}\build\fl_{}^{d\oplus f} \Map(K,F_{n-b}(X))^\bullet_{a-1}\otimes S^{\otimes b}\oplus\Map(K,F_{n-b}(Y))_a\otimes
 S^{\otimes b}$$
$$\build\fl_{}^{-f\oplus d}\Map(K,F_{n-b}(Y))^\bullet_{a-1}\otimes S^{\otimes b}\fl 0$$
By induction, $d(\theta)s$ is well defined in $\Map(K,F_{n-b}(X))_{a-1}\otimes S^{\otimes b}$ and $d(\theta)s+\overline\theta d(s)$ is well-defined in 
$\Map(K,F_{n-b}(X))'_{a-1}\otimes S^{\otimes b}$. Moreover $d(\theta)s+\overline\theta d(s)\oplus\theta f' s$ belongs to the kernel of $-f\oplus d$.

Therefore, there is a map $\theta s$ in $\Map(K,F_{n-b}(X))_a\otimes S^{\otimes b}$ such that: $d(\theta s)=d(\theta)s+\overline\theta d(s)$ and $f(\theta x)=\theta f's$.

Thus we are able to construct the desired $\Lambda$-action  by induction and the lemma is proven if $\alpha:F_{n-1}(X)\fl  F_{n-1}(X')$ is an isomorphism.

Consider the general case. Let $p<n$ be an integer such that we have constructed a compatible $\Lambda$-action on $F_{p-1}(X')$.

We set: $X'_1=\Phi_p(X')$, $X_1=\Phi_p(X)+\Phi_{p-1}(X')\subset X'_1$, $Y_1=\Phi_p(Y)$ and $Y'_1=\Phi_p(Y')$. We can apply the lemma with the diagram: 
$$\diagram{X_1^\circ&\hfl{}{}&X'_1\cr\vfl{}{}&&\vfl{}{}\cr Y_1^\circ&\hfl{}{}&{Y'_1}^\circ\cr}$$
in $\F_p(\Z)\bv$ and we get the desired $\Lambda$-action on $F_p(X')$. The lemma follows by induction.\cqfd
\vskip 12pt
\noi{\bf 4.4 Lemma:} {\sl Let $\Lambda$ be an algebra in $\R_0$, $n>0$ be an integer, $X$ be an $\F_n(\Lambda)\bv$-module, $C$ be a finite $\A$-complex and 
$f:C\fl F_n(X)^\circ$ be a morphism of complexes.

Then there is an $\F_n(\Lambda)$-module $Y$ such that: $F_n(Y)^\circ=F_{n-1}(Y)^\circ\oplus C$ and an $\F_n(\Lambda)\bv$-morphism $g:Y\fl X$ such that: $g=f$ on $C$.}
\vskip 12pt
\noi{\bf Proof:} There is an $\F_n(\Z)\v$-module $Z$ and a morphism $f:Z\fl X^\circ$ with the following properties, for any integer $i\leq n$:

1) $f:F_i(Z)\fl F_i(X)^\circ$ is a surjective homology equivalence

2) $F_i(Z)/F_{i-1}(Z)$ is bounded from below and free over $A$.

Because of Lemma 4.3 applied to the diagram:
$$\diagram{\Phi_{i-1}(Z)^\circ&\hfl{}{}&\Phi_i(Z)\cr\vfl{}{}&&\vfl{}{}\cr\Phi_{i-1}(X)^\circ&\hfl{}{}&\Phi_i(X)^\circ\cr}$$
we can construct, by induction, a $\Lambda$-action on $Z$ compatible with $f$. Then $Z$ belongs to $\F_n(\Lambda)\bv$ and $f$ is an $\F_n(\Lambda)\bv$-morphism.

Let $Z'$ be the $\F_n(\Lambda)\bv$-module defined by: $F_i(Z')=F_i(Z)$ for $i<n$ and $F_n(Z')=F_{n-1}(Z)$. We have also an $\F_n(\Z)\bv$-module $Z''$ defined by:
$F_i(Z'')=F_i(Z)$ for $i<n$ and $F_n(Z'')=F_{n-1}(Z)\oplus C$. By applying Lemma 4.3 to the diagram:
$$\diagram{Z'^\circ&\hfl{}{}&Z''\cr\vfl{}{}&&\vfl{}{}\cr X^\circ&\hfl{=}{}&X^\circ\cr}$$
we get a $\Lambda$-action on $Z''$ compatible with the inclusion $Z'\subset Z''$ and the morphism $Z''\fl X$.

Let $\widehat B$ be a basis of $F_{n-1}(Z)$ compatible with the degree and the inclusions $F_{i-1}(Z)\subset F_i(Z)$. If $B$ is a subset of $\widehat B$ the submodule of 
$F_n(Z)$ generated by $B$ will be denoted $A[B]$.

Since $C$ is a finite complex, there is a finite set $B\subset \widehat B$ such that: $d(B)\subset A[B]$ and $\theta(C)$ and $\theta(B)$ are contained in $A[B]\otimes 
S^{\otimes p}$ for any $\theta\in\Lambda$ of weight $p>0$.

Then the $\F_n(\Lambda)$-module $Y$ generated by $B$ and $C$ has the desired property and the lemma is proven.\cqfd
\vskip 12pt
\noi{\bf 4.5 Lemma:} {\sl Let $\Lambda$ be an algebra in $\R_0$ and $n>0$ be an integer. Then $\Nil(\Lambda)\v$ and $\Nil(\Lambda,n)\v$ are strict cocompletions of 
$\Nil(\Lambda)'$ and $\Nil(\Lambda,n)'$.
Moreover the inclusions $\Nil(\Lambda)'\subset \Nil(\Lambda)$ and $\Nil(\Lambda,n)'\subset \Nil(\Lambda,n)$ are dominations and $\Nil(\Lambda)'$ and $\Nil(\Lambda,n)'$ 
are regular exact categories.}
\vskip 12pt
\noi{\bf Proof:} we set: $\Nil=\Nil(\Lambda)$, $\Nil'=\Nil(\Lambda)'$ and $\Nil\v=\Nil(\Lambda)\v$ or: $\Nil=\Nil(\Lambda,n)$, $\Nil'=\Nil(\Lambda,n)'$ and 
$\Nil\v=\Nil(\Lambda,n)\v$.

It is easy to check conditions (cc0), (cc1) and (cc2) for the inclusion $\Nil'\subset\Nil\v$.

Let $X$ be a non-zero $\Nil\v$-module. Let $p$ be the smallest integer such that $T_p(X)\not=0$ and $i$ be the smallest integer such that $T_p(X)_i\not=0$. Let
$M$ be an $\A$-module and $f:M\fl T_p(X)_i$ be a non zero morphism. Consider $M$ as an $\A$-complex concentrated in degree $i$ and equipped with the trivial 
$\Lambda$-action. This module belongs to $\Nil'$ and we have a non zero morphism $M\fl X$ with $M\in\Nil'$. Hence, condition (cc3) is satisfied.

Let $X$ be a $\Nil'$-module, $U\fl V$ be an $\Nil\v$-epimorphism and $f:X\fl V$ be a $\Nil\v$-morphism. Let $W$ be the $\Nil\v$-module defined by the cartesian square:
$$\diagram{W&\hfl{}{}&U\cr\vfl{g}{}&&\vfl{}{}\cr X&\hfl{f}{}&V\cr}$$
Then the morphism $g:W\fl X$ is an epimorphism.

Suppose $\Nil$ is the category $\Nil(\Lambda)$. Then, for $n$ large enough, $X$ is in $\Nil(\Lambda,n)'$ and $T_n(W)\fl X$ is an epimorphism. Then, up to replacing $W$ by 
$T_n(W)$, we may as well suppose, in any case, that $X$ is a $\Nil(\Lambda,n)'$-module and that $W$ is a $\Nil(\Lambda,n)\v$-module.

Since $X$ is admissible, there is a filtration:
$$0=X_0\subset X_1\subset\dots\subset X_m=X$$
such that each quotient $X_i/X_{i-1}$ is in $\Nil(\Lambda,1)$. We set: $X_{ij}=T_j(X_i)$, $W_i=T(g^{-1}(X_i))$ and $W_{ij}=F_i(W_i)=T_j(g^{-1}(X_i))$.

For any $i\leq m$, there is a finite complex $K_i$ and a morphism $h_i:K_i\fl F_n(W_i)$ such that the composite morphism $K_i\build\fl_{}^{h_i}F_n(W_i)\build\fl_{}^g X_i
\fl X_i/X_{i-1}$ is onto. The kernel of this morphism will be denoted $K'_i$.

Because of Lemma 4.4, there are $\F_n(\Lambda)$-modules $Y_i$ and morphisms $u_i:Y_i\fl W_i$ such that $F_n(Y_i)^\circ=F_{n-1}(Y_i)^\circ\oplus K_i$ and $u_i=h_i$ on
$K_i$.

For any $i\leq m$, we set: $Z_i=\build\oplus_{j\leq i}^{} Y_i$. Morphisms $u_i$ induce morphisms $v_i:Z_i\fl W_i$.

Let $i\leq m$. Suppose $g\circ v_{i-1}:F_n(Z_{i-1})\fl F_n(X_{i-1})$ is a $\Nil(\Lambda,n)'$-deflation. Let $Y'_i$ be the $\F_n(\Lambda)$-module defined by: 
$F_i(Y'_i)=F_i(Y_i)$ for $i<n$ and $F_n(Y'_i)=F_{n-1}(Y_i)\oplus K'_i$. Then we have:
$$\Ker(F_n(Z_{i-1}\oplus Y_i)\rightarrow F_n(X_i))=\Ker(F_n(Z_{i-1}\oplus Y'_i)\rightarrow F_n(X_{i-1}))$$
and, since $F_n(Z_{i-1})\fl F_n(X_{i-1})$ is a $\Nil(\Lambda,n)'$-deflation, we have an exact sequence:
$$0\!\fl\! \Ker(F_n(Z_{i-1})\rightarrow F_n(X_{i-1}))\!\fl\! \Ker(F_n(Z_{i-1}\oplus Y'_i)\rightarrow F_n(X_{i-1}))\!\fl\! F_n(Y'_i)\!\fl\! 0$$
But $\Ker(F_n(Z_{i-1})\rightarrow F_n(X_{i-1}))$ and $F_n(Y'_i)$ are admissible. Then the kernel of $F_n(Z_{i-1}\oplus Y'_i)\rightarrow F_n(X_{i-1})$ is also 
admissible and $F_n(Z_i)\fl F_n(X_i)$ is a $\Nil(\Lambda,n)'$-deflation.

Thus we are able to prove by induction that each morphism $F_n(Z_i)\fl F_n(X_i)$ is a $\Nil(\Lambda,n)'$-deflation  and we get a commutative diagram:
$$\diagram{F_n(Z_m)&\hfl{v_m}{}&W\cr\vfl{gv_m}{}&&\vfl{}{}\cr X&\hfl{=}{}&X\cr}$$
where $gf_m$ is a $\Nil(\Lambda,n)'$-deflation.

Thus the condition (cc4) is satisfied and $\Nil\v$ is a strict cocompletion of $\Nil'$.

Let $X$ be a $\Nil$-module. Because of Lemma 4.4, there is a $\Nil'$-module $Y$ and a $\Nil$-deflation $Y\fl X$.

Let $X$ be a $\Nil'$-module, $g:Y\fl Z$ be a $\Nil$-deflation and $f:X\fl Z$ be a $\Nil$-morphism. Then $Y$ and $Z$ are $\Nil\v$-modules and $g$ is an epimorphism. Since
$\Nil'\subset\Nil$ is a strict cocompletion, the condition (cc4) is satisfied and there is a commutative diagram in $\Nil$:
$$\diagram{X'&\hfl{}{}&Y\cr\vfl{h}{}&&\vfl{g}{}\cr X&\hfl{f}{}&Z\cr}$$
where $h$ is a $\Nil'$-deflation. Hence, $\Nil'\subset \Nil$ is a domination.

Let $\C$ be a cocomplete exact class in $\Nil\v$ containing $\Nil'$.

Let $X$ be a $\Nil(\Lambda,1)\v$-module. Suppose the differential on $X$ is zero. Then, because $\A$ is regular, $X$ belongs to $\C$. Otherwise, we have an exact sequence:
$$0\fl Z\fl X\fl X/Z\fl 0$$
where $Z$ is the complex of cycles in $C$. Since $Z$ and $X/Z$ have zero differentials, we have: $Z\in\C$ and $C/Z\in \C$. Thus $X$ belongs to $\C$.

Let $X$ be a $\Nil\v$-module contained in $\Nil(\Lambda,p)\v$, with $p>1$. We have an exact sequence:
$$0\fl Y\fl X\fl Z\fl 0$$
where $Y$ is in $\Nil(\Lambda,p-1)\v$ and $Z$ is in $\Nil(\Lambda,1)\v$. Thus, by induction, we see that $\C$ contains each class $\Nil(\Lambda,p)\v$ and then we have:
$\C=\Nil\v$. Hence, $\Nil'$ is a regular exact category.\cqfd
\vskip 12pt
Let $H$ and $K$ be two $\F_n(\Lambda)\bv$-modules and $f:H\fl K$ be an $\F_n(\Lambda)\bv$-morphism. For any $i<n$, the cylinder of the inclusion $F_{i-1}(H)\subset 
F_i(H)$ will be denoted $E_i$ and the cylinder of $f:F_{n-1}(H)\fl F_n(K)$ will be denoted $E_n$. Thus we have an $\F_n(\Z)\bv$-module $\Gamma(f)$, where 
$F_i(\Gamma(f))$ is the module $\build\oplus_{j\leq i}^{} E_j$.

More precisely, we have formal elements $e_i$ of degree $0$ and $t_i$ of degree $1$, with $i=1,2,\dots,n$, such that any element $u$ in $F_i(\Gamma(f))$ can be written 
in a unique way:
$$u=\left\{\matrix{\build\sum_{j\leq i}^{}(e_j x_j+t_jy_{j-1})&\ \hbox{if}\ i<n\cr \build\sum_{j<n}^{}(e_jx_j+t_{j+1}y_j)+e_nz&\hbox{if}\ i=n\cr}\right.$$
with $x_j,y_j\in F_j(H)$ and $z\in F_n(K)$. The differential is defined as follows:
$$d(e_j x)=e_jd(x)\hskip 24pt d(e_n z)=e_n d(z)$$
$$d(t_{j+1}x)=\left\{\matrix{e_jx-e_{j+1} x-t_{j+1}d(x)&\hbox{if}\ j<n-1\cr e_{n-1}x-t_n d(x)-e_n f(x)&\hbox{if}\ j=n-1\cr}\right.$$
for any $j<n$, any $x\in F_j(H)$ and any $z\in F_n(K)$.

Moreover, we have a natural $\F_n(\Z)\bv$-morphism $\widehat f:\Gamma(f)\fl K$ such that, for any $i<n$, any $x\in F_i(H)$ and any $z\in F_n(K)$, we have:
$$\widehat f(e_i x)=f(x)\hskip 24pt \widehat f(t_{i+1}x)=0\hskip 24pt \widehat f(e_n z)=z$$

For any $i<n$, the submodule of $F_i(\Gamma(f))$ generated by $e_{i-1}x-e_i x$ and $t_i x$, with $x\in F_{i-1}(H)$ will be denoted $V_i$ and the submodule of 
$F_n(\Gamma(f))$ generated by $e_{n-1}x-e_n f(x)$ and $t_n x$ with $x\in F_{n-1}(H)$ will be denoted $V_n$. Then we get a $\F_n(\Lambda)\bv$-module $\widehat\Gamma(f)$ 
with: $F_i(\widehat\Gamma(f))=\build\oplus_{j\leq i}^{} V_j$ where the $\Lambda$-action is given by:
$$\theta(e_i x-e_{i+1}x)=e_{i-p}\theta(x)-e_{i+1-p}\theta(x)\hskip 24pt\theta(t_{i+1}x)=t_{i+1-p}\overline\theta(x)$$
$$\theta(e_{n-1}x-e_nf(x))=e_{n-1-p}\theta(x)-e_nf(\theta(x))$$
for any $i<n$, any $x\in F_i(H)$ and any $\theta\in\Lambda$ of weight $p>0$.

Each module $F_i(\widehat\Gamma(f))$ is acyclic and we have a sequence in $\F_n(\Z)\bv$:
$$0\fl \widehat\Gamma(f)\build\fl_{}^g\Gamma(f)\build\fl_{}^{\widehat f} K\fl 0$$
with $\widehat f\circ g=0$ and this sequence induces an exact sequence:
$$0\fl F_n(\widehat\Gamma(f))\build\fl_{}^gF_n(\Gamma(f))\build\fl_{}^{\widehat f} F_n(K)\fl 0$$

Let $\Lambda$ be an algebra in $\R_0$, $n>0$ be an integer and $X$ be an $\F_n(\Lambda)\bv$-module. We denote $\E(X)$ the class of triples $(H,K,f)$ where $K$ is a 
$\F_n(\Lambda)$-module, $f:K\fl X$ is an $\F_n(\Lambda)\bv$-morphism and $H$ is a sub-$\F_n(\Lambda)$-module of $K$ such that:

1) $K^\circ=\Gamma(f:H\fl X)$ 

2) for any $\theta\in\Lambda$ of weight $p>0$, we have:
$$\theta(e_n F_n(X))\subset \build\oplus_{i\leq n-p}^{}\Bigl(e_iF_i(H)\oplus t_i F_{i-1}(H)\Bigr)\otimes S^{\otimes p}$$

3) for any $i<n$, any $x\in F_i(H)$ and any $\theta\in\Lambda$ of weight $p$, we have:
$$\theta(e_i x)=e_{i-p}\theta(x)\hskip 24pt \theta(t_{i+1} x)=t_{i+1-p}\overline\theta(x)$$

4) the morphism $\widehat f:K\fl X$ is a $\F_n(\Lambda)\bv$-morphism.
\vskip 12pt
\noi{\bf 4.6 Lemma:} {\sl  Let $\Lambda$ be an algebra in $\R_0$ and $n>0$ be an integer. Let $X$ be an $\F_n(\Lambda)\bv$-module such that $F_n(X)$ is a 
$\Nil(\Lambda,n)'$-module. Then $\E(X)$ is non-empty.}

\vskip 12pt
\noi{\bf Proof:} Since $F_n(X)$ is admissible, we have a filtration:
$$0=C_0\subset C_1\subset C_2\subset\dots\subset C_m=F_n(X)$$
such that each quotient $C_i/C_{i-1}$ is in $\Nil(\Lambda,1)$. 

By setting: $F_i(X_k)=C_k\cap F_i(X)$, we get a filtration of $\F_n(\Lambda)\bv$-modules:
$$0=X_0\subset X_1\subset X_2\subset\dots\subset X_m=X$$

It is easy to see that $\E(X_0)$ is non-empty.

Let $p<m$ be an integer such that $\E(X_p)$ is non-empty and $(H,K,f)$ be an object of $\E(X_p)$. 

We have an $\F_n(\Lambda)\bv$-morphism $\widehat f:K\fl X_p$ and an exact sequence:
$$0\fl F_n(\widehat\Gamma(f))\fl F_n(\Gamma(f))\fl F_n(X_p)\fl 0$$
Therefore, $\widehat f:F_n(\Gamma(f))\fl F_n(X_p)$ is a homology equivalence.

We can construct by induction an $\F_n(\Lambda)\bv$-module $Y$, morphisms $g:K\fl Y$ and $h:Y\fl X_p$ with the following properties:

1) for any $i\leq n$ the morphism $F_i(K)\fl F_i(Y)$ is monic and morphisms $F_i(Y)/F_i(K)$ and $F_i(Y)/(F_i(K)+F_{i-1}(Y))$ are free in each degree and bounded from 
below

2) for any $i<n$, $F_i(Y)\fl F_i(X_p)$ is a surjective homology equivalence

By applying Lemma 4.3 to the following diagram:
$$\diagram{\Phi_{i-1}(Y)^\circ\build\oplus_{\Phi_{i-1}(K)}^{} \Phi_i(K))^\circ&\hfl{}{}&\Phi_i(Y)\cr\vfl{h}{}&&\vfl{h}{}\cr 
\Phi_{i-1}(X_p))^\circ\build\oplus_{\Phi_{i-1}(K)}^{}\Phi_i(K))^\circ&\hfl{}{}&\Phi_i(X_p))^\circ\cr}$$
we can construct by induction a $\Lambda$-action on $Y$ which is compatible with morphisms $g$ and $h$ and we get an exact sequence in $\F_n(\Lambda)\bv$:
$$0\fl Z\fl Y\build\fl_{}^h X_p\fl 0$$

We set: $K'=\Gamma(h)$ and $H'=\widehat\Gamma(h)$. We have an exact sequence in $\F_n(\Z)\bv$:
$$0\fl H'\fl K'\build\fl_{}^{\widehat h} X_p\fl 0$$
Since each $F_i(Z)$ is acyclic, the morphism $\widehat h:K'\fl X_p$ induces surjective homology equivalences $F_i(K')\fl F_i(X_p)$ and, because of Lemma 4.3, the 
$\Lambda$-action on $H'$ can be extended to a $\Lambda$-action on $K'$.

Since each module $F_i(H')$ is acyclic, the inclusion $X_p\subset X_{p+1}$ can be lifted to an inclusion $Y\subset Y'$ and we have a morphism $h':Y'\fl X_{p+1}$ and a 
commutative diagram in $\F_n(\Z)\bv$ with exact lines:
$$\diagram{0&\hfl{}{}&H'&\hfl{}{}&K'&\hfl{\widehat h}{}&X_p&\hfl{}{}&0\cr &&\vfl{=}{}&&\vfl{}{}&&\vfl{}{}&&\cr
0&\hfl{}{}&H'&\hfl{}{}&K''&\hfl{\widehat h'}{}&X_{p+1}&\hfl{}{}&0\cr}$$
Then, because of Lemma 4.3, there is a compatible $\Lambda$-action on $K''$ such that this diagram can be lifted to a commutative diagram in $\F_n(\Lambda)\bv$.

Since $F_n(X_{p+1})$ is a finite $\A$-complex and $F_n(H')$ is acyclic, there is in $F_n(K'')$ a finite $\A$-complex $C$ such that $\widehat h'$ induces an
isomorphism $C\build\fl_{}^\sim F_n(X_{p+1})$.

Notice that $F_{n-1}(K''/K)$ is free in each degree. Then there is a set $\widehat B$ of homogeneous elements in $Y$ such that the set $\{e_i x,t_{i+1}x\}$ with $i<n$ 
and $x\in B\cap F_i(Y)$ induces a basis in $F_{n-1}(K'')/F_{n-1}(K)$.

If $B$ is a subset of $\widehat B$ the submodule generated by $B$ and $F_n(K)$ will be denoted $A[B]$. Then we can find a finite subset $B$ of $\widehat B$ with the 
following properties:

$\bullet$ for any $x\in B$, $d(x)$ belongs to $A[B]$

$\bullet$ for any $\theta\in\Lambda$ of weight $p>0$, we have: 
$$\theta(C)\subset\build\oplus_{i\leq n-p}^{}\Bigl(e_i(A[B]\cap F_i(Y))+t_i(A[B]\cap F_{i-1}(Y))\Bigr)\otimes S^{\otimes p}$$
$$\theta(B)\subset A[B]\otimes S^{\otimes p}$$

Since we have: $\theta(e_ix)=e_i\theta(x)$ and $\theta(t_{i+1}x)=t_{i+1}\overline\theta(x)$, for any $i<n$ and any $x\in F_i(H')$, we get an $\F_n(\Lambda)$-module 
$K_1\subset K''$ defined by:
$$F_i(K_1)=\build\oplus_{j\leq i}^{}\Bigl(e_j(A[B]\cap F_j(Y))+t_i(A[B]\cap F_{j-1}(Y))\Bigr)$$
for $i<n$ and:
$$F_n(K_1)=F_{n-1}(K_1)\oplus t_n(A[B]\cap F_{n-1}(Y))\oplus e_n C$$
We have also an $\F_n(\Lambda)$-module $H_1$ with: 
$F_i(H_1)=F_i(K_1)$ for $i<n$ and $F_n(H_1)=F_{n-1}(K_1)\oplus t_n (A[B]\cap F_{n-1}(Y))$. These two modules induce an object of $\E(X_{p+1})$ and $\E(X_{p+1})$ is 
non-empty.
\cqfd   
\vskip 12pt
\noi{\bf 4.7 Lemma:} {\sl  Let $\Lambda$ be an algebra in $\R_0$ and $n>0$ be an integer. Let $X$ be an $\F_n(\Lambda)\bv$-module such that $\E(X)$ is non-empty. Then
there is an object $(H,K,f)$ in $\E(X)$ and a morphism $\varphi:J_1(\Lambda)\otimes F_n(X)\fl F_n(H)\otimes A[S]$ of $\A\v$-complexes with the following properties:

1) for any $\theta\in\Lambda$ of weight $p>0$, we have: $\varphi(\theta\otimes F_n(X))\subset F_{n-p}(H)\otimes S^{\otimes p}$

2) for any $\theta\in J_1(\Lambda)$ and any $x\in F_n(X)$, we have: $\widehat f(\varphi(\theta\otimes x))=\theta f(x)$

3) for any $\theta\in\Lambda$ of weight $p>0$ and any $x\in F_n(X)$, we have: $\theta e_n x=e_{n-p}\varphi(\theta\otimes x)$.}
\vskip 12pt
\noi{\bf Proof:} Let $(H,K,f)$ be an object in $\E(X)$. We have an exact sequence in $\F_n(\Lambda)$:
$$0\fl F_n(H)\fl F_n(K)\fl F_n(X)\fl 0$$
The $\Lambda$-action on $K$ is induced by the $\Lambda$-action on $H$ and the action on $e_n F_n(X)$. Actually, for any $\theta\in\Lambda$ of weight $p>0$, we have
maps $\lambda_i(\theta)$ and $\mu_{i+1}(\theta)$ in $\Map(F_n(X),F_i(H)\otimes S^{\otimes p})$ such that:
$$\theta e_n x=\build\sum_{i\leq n-p}^{}(e_i\lambda_i(\theta) x+t_i\mu_i(\theta)x)$$
for any $x\in\F_n(X)$.

The fact that we have a $\Lambda$-action on $K$ implies the following:

$\bullet$ for any $i$ and any homogeneous $\theta\in\Lambda$ and $\theta'\in J_1(\Lambda)$, we have: 
$$\lambda_i(\theta\theta')=\theta\lambda_i(\theta')\hskip 24pt \mu_i(\theta\theta')=\theta\mu_i(\theta')$$

$\bullet$ for any $i$ and homogeneous $\theta\in\Lambda$, we have:
$$\sum e_i d(\lambda_i(\theta))+\sum(e_{i-1}-e_i)\mu_i(\theta)-\sum t_i d(\mu_i(\theta))=\sum e_i\lambda_i(d(\theta))+\sum t_i\mu_i(d(\theta))$$

and this last condition is equivalent to:
$$\mu_i(d(\theta))=-d(\mu_i(\theta))\hskip 24pt \lambda_i(d(\theta))=d(\lambda_i(\theta))-\mu_i(\theta)+\mu_{i+1}(\theta)$$

As maps $\lambda_i$ and $\mu_{i+1}$ in $\Map(J_1(\Lambda)\otimes F_n(X\oplus Y),F_*(H)\otimes S^{\otimes *})$, these conditions become:
$$d(\mu_i)=0\hskip 24pt d(\lambda_i)=\mu_i-\mu_{i+1}$$

We have new maps $\widehat\lambda_i$ and $\widehat\mu_{i+1}$ defined as follows:
$$\build\sum_{i\leq n-p}^{}(e_i\widehat\lambda_i(\theta)+t_i\widehat\mu_i(\theta))=
\build\sum_{i\leq n-p}^{}(e_i\lambda_i(\theta)+t_i\mu_i(\theta))-\build\sum_{j\leq i<n-p}^{}((e_i-e_{i+1})\lambda_j(\theta)-t_{i+1}d(\lambda_j)(\theta))$$
where $p$ is the weight of $\theta$.

After computation, we find:
$$\widehat\lambda_i(\theta)=\left\{\matrix{\build\sum_{i\leq n-p}^{}\lambda_i(\theta)&\hbox{if}\ i=n-p\cr 0&\hbox{otherwise}\cr}\right.
\hskip 24pt\hbox{and}\hskip 24pt \widehat\mu_i(\theta)=0$$

With these new maps, we get a new $\Lambda$-action on $K$ and, with this new $K$, the lemma is easy to prove.\cqfd

\vskip 12pt
\noi{\bf 4.8 Lemma:} {\sl  Let $\Lambda$ be an algebra in $\R_0$ and $n>0$ be an integer. Let $X$ be an $\F_n(\Lambda)\bv$-module such that $F_n(X)$ is a direct summand of
a $\Nil(\Lambda,n)'$-module. Then $\E(X)$ is non-empty.}
\vskip 12pt
\noi{\bf Proof:} Since $F_n(X)$ is a direct summand of an admissible module, there is an $\F_n(\Lambda)\bv$-module $Y$ such that $F_n(X\oplus Y)$ is admissible. Because
of lemmas 4.6 and 4.7, there is a object $(H,K,f)$ in $\E(X\oplus Y)$ and we have an exact sequence in $\F_n(\Lambda)$:
$$0\fl F_n(H)\fl F_n(K)\fl F_n(X\oplus Y)\fl 0$$
The $\Lambda$-action on $K$ is induced by the $\Lambda$-action on $H$ and the action on $e_n F_n(X\oplus Y)$ is defined by:
$$\theta e_n z=e_{n-p}\varphi(\theta\otimes z)$$
for any $\theta\in\Lambda$ of weight $p>0$ and any $z\in F_n(X\oplus Y)$.

The morphism $f:H\fl X\oplus Y$ is defined by two morphisms $f_1:H\fl X$ and $f_2:H\fl Y$.

Let $H_1$ be the kernel of $f_2:H\fl Y$. There is a $\F_n(\Z)\bv$-module $H'$ and a morphism $g:H'\fl H_1$ such that:

$\bullet$ for any $i\leq n$, $F_i(H')$ is bounded from below and free over $A$

$\bullet$ for any $i\leq n$, the morphism $g:F_i(H')\fl F_i(H_1)$ is a surjective homology equivalence.

Because of Lemma 4.3 applied to the diagram:
$$\diagram{\Phi_{i-1}(H')&\hfl{}{}&\Phi_i(H')\cr\vfl{}{}&&\vfl{}{}\cr\Phi_{i-1}(H_1)&\hfl{}{}&\Phi_i(H_1)\cr}$$
we can construct, by induction, a $\Lambda$-action on $H'$ compatible with the morphism $g:H'\fl H_1$.

Let $x$ be an element of $F_n(X)$ and $\theta\in\Lambda$ of weight $p>0$. We have: $\widehat f(\varphi(\theta\otimes x)=\theta f(x)\in F_{n-p}(X)$. Then 
$\varphi(\theta\otimes x)$ is contained in $F_{n-p}(H_1)\otimes S^{\otimes p}$ and $\varphi$ is a morphism from $J_1(\Lambda)\otimes F_n(X)$ to $F_n(H_1)\otimes A[S]$.
But $H'\fl H_1$ is a surjective homology equivalence. Then $\varphi$ can be lifted to a morphism:
$$\psi: J_1(\Lambda)\otimes F_n(X)\fl F_n(H')\otimes A[S]$$

Morphisms $H'\fl H$ and $X\subset X\oplus Y$ induce a morphism $\Gamma(g:H'\fl X)\fl \Gamma(H\fl X\oplus Y$. Moreover $\psi$ induce a $\Lambda$-action on 
$K'=\Gamma(H'\fl X)$.

Let $\widehat B$ be a homogeneous basis of $H'$ compatible with the inclusions $F_i(H')\subset F_{i+1}(H')$. As before, if $B$ is a subset of $\widehat B$, the submodule
of $F_n(H')$ generated by $B$ will be denoted $A[B]$.

Hence, we can find a finite subset $B$ of $\widehat B$ such that:

$\bullet$ $d(B)\subset A[B]$

$\bullet$ for any $\theta\in\Lambda$ of weight $p>0$, $\theta(B)\subset A[B]\otimes S^{\otimes p}$

$\bullet$ for any $\theta\in\Lambda$ of weight $p>0$, $\theta(e_n F_n(X))\subset e_{n-p} A[B]\otimes S^{\otimes p}$

Let $H''$ be the module $A[B]$ and $K''$ be the submodule of $K'$ generated by $H''$ and $e_nF_n(X)$. The data $(H'',K'',\widehat g)$ is an object of $\E(X)$ and this  
class is non-empty.\cqfd
   
\vskip 12pt
\noi{\bf 4.9 Lemma:} {\sl  Let $\Lambda$ be an algebra in $\R_0$ and $n>0$ be an integer. Let $X$ be an $\F_n(\Lambda)\bv$-module such that $F_n(X)$ is a 
$\Nil(\Lambda,n)$-module. Then there is a $\F_n(\Lambda)$-module $Y$ and an $\F_n(\Lambda)\bv$-morphism $Y\fl X$ such that $F_n(Y)\fl F_n(X)$ is a surjective homology
equivalence.}
\vskip 12pt
\noi{\bf Proof:} Because of Lemma 4.4, we can find a $\F_n(\Lambda)\bv$-complex:
$$\dots\fl C_p\fl C_{p-1}\fl\dots\fl C_1\fl C_0\fl X\fl 0$$
inducing a $\Nil(\Lambda,n)'$-resolution of $F_n(X)$:
$$\dots\fl F_n(C_p)\fl F_n(C_{p-1})\fl\dots\fl F_n(C_1)\fl F_n(C_0)\fl F_n(X)\fl 0$$
Let $Z_p$ be the kernel of $C_p\fl C_{p-1}$. Since $\Nil(\Lambda,n)'$ is regular, Theorem 2.1 implies that $Z_p$ is, for $p$ large enough, a direct summand of an  
admissible module. Because of Lemma 4.8, there is a $\F_n(\Lambda)$-morphism $g:K\fl Z_p$ inducing a surjective homology equivalence $g:F_n(K)\fl F_n(Z_p)$. Thus we have a 
complex:
$$0\fl K\fl C_p\fl C_{p-1}\fl\dots\fl C_0\fl X$$
inducing an exact sequence:
$$F_n(K)\fl F_n(C_p)\fl F_n(C_{p-1})\fl\dots\fl F_n(C_0)\fl F_n(X)\fl 0$$
Let $\widehat C$ be the total complex of: $0\fl K\fl C_p\fl C_{p-1}\fl\dots\fl C_0\fl 0$. This complex is an $\F_n(\Lambda)$-module and the morphism $C_0\fl X$ induces 
an $\F_n(\Lambda)\v$-morphism $\widehat C\fl X$. Since the kernel of $F_n(K)\fl F_n(C_p)$ is acyclic, the morphism $\widehat C\fl X$ induces a surjective homology 
equivalence $F_n(\widehat C)\fl F_n(X)$.\cqfd
\vskip 12pt
\noi{\bf 4.10 Corollary:} {\sl Let $\Lambda$ be an algebra in $\R_0$ and $n>0$ be an integer. Then we have: $Nil_0(\Lambda,n)=0$ and $Nil_0(A,S,n)=0$.}
\vskip 12pt
\noi{\bf Proof:} If $X$ is a $\Nil(\Lambda,n)$-module, its class in $Nil(\Lambda,n)$ will be denoted $[X]$. Let $X$ be a $\Nil(\Lambda,n)$-module. Because of Lemma 4.9,
there is an $\F_n(\Lambda)$-module $Y$ and a surjective homology equivalence $f:F_n(Y)\fl X$. Thus we have:
$$[X]=[F_n(Y)]=0$$
and we have: $Nil_0(\Lambda,n)=0$. 

If $\Lambda$ is the algebra $\Z[\theta]$ with $\partial^\circ(\theta)=(0,1)$, we get (because of Theorem 1.17 and Lemma 4.1):
$$Nil_0(A,S,n)=Nil_0(\Lambda,n)=0$$
and the result follows.\cqfd
\vskip 12pt
\noi{\bf 4.11 Lemma:} {\sl Let $\Lambda$ be an algebra in $\R_0$, $n>0$ be an integer, $X$ be a $\F_n(\Lambda)$-module, $Y$ be a $\Nil(\Lambda,n)$-module and $f:F_n(X)\fl
Y$ be a $\Nil(\Lambda,n)$-morphism. Then $f$ factors through a $\F_n(\Lambda)$-module $Z$, via a $\F_n(\Lambda)$-morphism $X\fl Z$ and a surjective homology equivalence 
$F_n(Z) \build\fl_{}^\sim Y$.}
\vskip 12pt
\noi{\bf Proof:} The morphism $f:F_n(X)\fl Y$ induces a morphism $\widehat f:X\fl T(Y)$. Let $E$ be $(-1)$-cone of $\widehat f$. Then we have an exact sequence:
$$0\fl T(Y)\fl E\fl X\fl 0$$
in $\F_n(\Lambda)\bv$.

Let $H$ be the $0$-cone of $E\fl X$. This module is the cylinder of $\widehat f$ and we have a morphism $H\fl T(Y)$ inducing a surjective homology equivalence
$F_n(H)\fl F_n(T(Y))=Y$.

Because of Lemma 4.9, there is a $\F_n(\Lambda)$-module $K$ and a morphism $g:K\fl E$ inducing a surjective homology equivalence $F_n(K)\fl F_n(E)$. Let $Z$ be
the $0$-cone of the composite morphism $K\fl E\fl X$. We have a morphism $Z\fl H$ inducing a surjective homology equivalence $F_n(Z))\fl F_n(H)$. Hence, the 
morphism $\widehat f:X\fl T(Y)$ factors through the $\F_n(\Lambda)$-module $Z$ via a morphism $Z\fl T(Y)$ inducing a surjective homology equivalence $F_n(Z)\fl 
F_n(T(Y))=Y$.\cqfd
\vskip 12pt
If $\Lambda$ is an algebra in $\R_0$ and $n$ is an integer, the category of $\F_n(\Lambda)$-modules $X$ such that $F_n(X)$ is acyclic will be denoted $\F_n(\Lambda)^0$.
This category is an essentially small exact category and belongs to {\bf WE}$^h$.
\vskip 12pt
\noi{\bf 4.12 Lemma:} {\sl Let $\Lambda$ be an algebra in $\R_0$ and $n>0$ be an integer. Then the inclusion $\F_n(\Lambda)^0\subset \F_n(\Lambda)$ and the functor 
$F_n$ induce a homotopy fibration in K-theory:
$$K(\F_n(\Lambda)^0)\fl K(\F_n(\Lambda))\fl K(\Nil(\Lambda,n))$$}
\vskip 12pt
\noi{\bf Proof:} Let $w$ be the class of morphisms $f$ in $\F_n(\Lambda)$ such that $F_n(f)$ is a homology equivalence. Equipped with this class $w$, we get a
new category $w\F_n(\Lambda)$ in $\WE^h$ and, because $\F_n(\Lambda)$ is in $\WE^h$, Waldhausen's Localization Theorem implies that we have a homotopy fibration:
$$K(\F_n\Lambda)^0)\fl K(\F_n(\Lambda))\fl K(w\F_n(\Lambda))$$
On the other hand, the functor $F_n:w\F_n(\Lambda)\fl\Nil(\Lambda,n)'$ has the property:

$\bullet$ a morphism in $w\F_n(\Lambda)$ is an equivalence if and only if $F_n(f)$ is an equivalence.

Then, because of Lemma 4.11, the functor $F_n:w\F_n(\Lambda)\fl\Nil(\Lambda,n)$ has the approximation property and, because of Waldhausen's Approximation Theorem, 
$F_n$ induces a homotopy equivalence in K-theory. The result follows.\cqfd
\vskip 12pt
By shifting filtrations, we get an exact functor $\Sigma:\F_{n-1}(\Lambda)\fl \F_n(\Lambda)$ such that: $\Sigma(F_i(X))=F_{i+1}(X)$. This functor induces an exact 
functor $\Sigma:\F_{n-1}(\Lambda)^0\fl\F_n(\Lambda)^0$.

Let $w$ the class of $\F_n(\Lambda)^0$-morphisms $f:X\fl Y$ inducing an homology equivalence $F_1(X)\build\fl_{}^\sim F_1(Y)$. With these equivalences, we get a new 
Waldhausen exact category $w\F_n(\Lambda)^0\in\WE^h$.
\vskip 12pt
\noi{\bf 4.13 Lemma:} {\sl Let $\Lambda$ be an algebra in $\R_0$ and $n>1$ be an integer. Then we have a commutative diagram:
$$\diagram{K(\F_{n-1}(\Lambda)^0)&\hfl{}{}&K(\F_{n-1}(\Lambda))&\hfl{}{}&K(\Nil(\Lambda,n-1))\cr\vfl{\Sigma}{}&&\vfl{\Sigma}{}&&\vfl{}{}\cr
K(\F_n(\Lambda)^0)&\hfl{}{}&K(\F_n(\Lambda))&\hfl{}{}&K(\Nil(\Lambda,n))\cr\vfl{}{}&&\vfl{F_1}{}&&\cr K(w\F_n(\Lambda)^0)&\hfl{F_1}{}&K(\A_*)&&\cr}$$
where the two top lines and the two left columns are fibrations.}
\vskip 12pt
\noi{\bf Proof:} Let $\F_n(\Lambda)^{00}$ be the category of $\F_n(\Lambda)^0$-modules $X$ such that $F_1(X)$ is acyclic. Waldhausen's Localization Theorem 
implies we have a homotopy fibration:
$$K(\F_n(\Lambda)^{00})\fl K(\F_n(\Lambda)^0)\fl K(w\F_n(\Lambda)^0)$$
Let $X$ be a $\F_n(\Lambda)^{00}$-module. We have an exact sequence in $\F_n(\Lambda)^{00}$:
$$0\fl \Gamma_1(X)\fl X\fl \Sigma(\Gamma_2(X))\fl 0$$
where $\Gamma_1(X)$ is the complex $F_1(X)=X_1$ filtered by $0\subset X_1\subset\dots\subset X_1$ and $\Gamma_2(X)$ is the complex $X/X_1$ filtered by 
$0\subset X_2/X_1\subset\dots \subset X_n/X_1$. Thus $\Gamma_2$ is an exact functor from  $\F_n(\Lambda)^{00}$ to $\F_{n-1}(\Lambda)^0$.

Since $\Gamma_1(X)$ is acyclic in $\F_n(\Lambda)^{00}$, the Additivity Theorem implies that the functor $\Gamma_2$ induces a homotopy equivalence 
$K(\F_n(\Lambda)^{00})\build\fl_{}^\sim K(\F_{n-1}(\Lambda)^0)$. Moreover the functor $\Sigma:\F_{n-1}(\Lambda)^0\fl \F_n(\Lambda)^{00}$ induces a homotopy inverse of 
it: $K(\F_{n-1}(\Lambda)^0)\build\fl_{}^\sim K(\F_n(\Lambda)^{00})$. The result follows.\cqfd
\vskip 12pt
Let $\F_n(\Lambda)^1$ be the category of pairs $(X,t)$ where $X$ is a $\F_n(\Lambda)$-module and $t:F_n(X)^\circ\fl F_n(X)^\circ$ is a $1$-map from the underlying 
$\A$-complexes of $F_n(X)$ to itself such that: $t^2=0$ and $d(t)=$Id. A morphism $\varphi:(X,t)\fl (X',t)$ is this category is a $\F_n(\Lambda)$-morphism 
$\varphi:X\fl X'$ such that: $\varphi t=t\varphi$. It is easy to verify that $\F_n(\Lambda)^1$ is a category in $\WE^h$ and that the correspondence $(X,t)\mapsto X$ is 
a forgetful exact functor $\F_n(\Lambda)^1\fl \F_n(\Lambda)^0$ which respects the class $w$.
\vskip 12pt
\noi{\bf 4.14 Lemma:} {\sl The forgetful functor $w\F_n(\Lambda)^1\fl w\F_n(\Lambda)^0$ induces a homotopy equivalence in K-theory: $K(w\F_n(\Lambda)^1)\build\fl_{}^\sim 
K(w\F_n(\Lambda)^0)$.}
\vskip 12pt
\noi{\bf Proof:} 
Let $(X,t)$ be a $\F_n(\Lambda)^1$-module and $f:X\fl Y$ be a $\F_n(\Lambda)^0$-morphism. Let $Z$ be the cylinder of $f$. We have a $\F_n(\Lambda)^0$-inflation 
$g:X\fl Z$ and a homology equivalence $h:Z\build\fl_{}^\sim Y$.

In particular we have an inflation of acyclic finite $\A$-complexes $g:X^n\fl Z^n$ with cokernel $K$. Since these complexes are acyclic, we have a decomposition:
$$Z^n=X^n\oplus K$$
Since $K$ is acyclic, there is an $1$-map $t:K\fl K$ such that: $t^2=0$ and $d(t)$ is the identity of $K$. Maps $t:X^n\fl X^n$ and $t:K\fl K$ induce a $1$-map 
$t:Z^n\fl Z^n$ compatible with the inclusion $X^n\subset Z^n$.

Hence, the approximation property is satisfied and Waldhausen 's Approximation Theorem implies the result.\cqfd
\vskip 12pt

Let $\Lambda$ be an algebra in $\R_0$ and $n>0$ be an integer.
Consider $\Z[t]/(t^2)$ as a bigraded differential algebra with $\partial^\circ(t)=(1,0)$ and $d(t)=1$. Let $\Lambda^+$ be the free product $\Lambda*\Z[t]/(t^2)$ and
$\Lambda_n^+$ be the quotient of $\Lambda^+$ by the two-sided ideal generated by $J_n(\Lambda)$. This algebra is a bigraded differential algebra and belongs to $\R$.

Let $X$ be a $\F_n(\Lambda)^1$-module. Then the $\Lambda$-action and the action of $t$ induce a $\Lambda_n^+$-action on $F_n(X)=X_n$.

The algebra $\Lambda_n^+$ is generated by elements of the form: $u=u_0tu_1t\dots tu_p$ with $p\geq 0$, where each $u_i$ is a homogeneous element in $\Lambda$
such that:
$$\forall i=0,1,\dots,p,\ \ \partial_2^\circ(u_i)<n\ \ \ \hbox{and}\ \ \ \forall i\ \hbox{with}\ 0<i<p,\ \ \partial_2^\circ(u_i)>0$$

Moreover if $\partial_2^\circ(u_p)>0$ then this element acts trivially on $F_1(X)=X_1$ and, in order to be sure that $uX_1$ is contained in $X_1\otimes S^{\otimes *}$, 
we need to have: $p=0$ and $\partial_2^\circ(u_0)=0$ or $p>0$ and $\partial_2^\circ(u_0)=n-1$.

Let $\Lambda^{(n)}$ be the subalgebra of $\Lambda_n^+$ generated by elements $u=u_0tu_1t\dots u_{p-1}t$ such that: $p>0$, $\partial_2^\circ(u_0)=n-1$ and 
$0<\partial_2^\circ(u_i)<n$ for all $i>0$. 

Then the $\Lambda_n^+$-action on $X_n$ induces a $\Lambda^{(n)}$-action on $X_1$.
\vskip 12pt 
\noi{\bf 4.15 Lemma:} {\sl Let $n>0$ be an integer and $\Lambda$ be an algebra in $\R_0$. Then the algebra $\Lambda^{(n)}$ belongs to the class 
$\R_0$ and the functor $X\mapsto X_1$ induces an exact functor $\Psi_n:\F_n(\Lambda)^1\fl \Nil(\Lambda^{(n)})$.} 
\vskip 12pt 
\noi{\bf Proof:} Let $U$ be a basis of $\Lambda$ and $\overline U$ be the set of words in $U$. Then $\Lambda^{(n)}$ is freely generated by elements 
$u=u_0tu_1t\dots u_{p-1}t$ in $\Lambda_n^+$ such that: $p>0$, $u_i\in\overline U$ for all $i$, $\partial_2^\circ(u_0)=n-1$ and $0<\partial_2^\circ(u_i)<n-1$ for all
$i=1,2,\dots,p-1$. Let $\widehat U$ be the set of these elements. Then $\widehat U$ is a basis of $\Lambda^{(n)}$ and it is easy to see that $\Lambda^{(n)}$ is an 
algebra in $\R_0$.

Let $X$ be a $\F_n(\Lambda)^1$-module. Let $u=u_0tu_1t\dots u_{p-1}tu_p$ be an element in $\Lambda_n^+$ where each $u_i$ is a homogeneous element in $\Lambda$. Suppose 
$u$ acts non trivially on $X$. Then its degree is bounded. Therefore, $p$ is bounded and the weight is also bounded. Thus the action of $\Lambda^{(n)}$ on $X$ is
nilpotent and $X_1$ is a $\Nil(\Lambda^{(n)})$-module. The result follows.\cqfd

\vskip 12pt
\noi{\bf 4.16 Lemma:} {\sl Let $n>0$ be an integer and $\Lambda$ be an algebra in $\R_0$. Then the functor $\Psi_n:\F_n(\Lambda)^1\fl \Nil(\Lambda^{(n)})$ induces a
homotopy equivalence:
$$K(w\F_n(\Lambda)^1)\build\fl_{}^\sim K(\Nil(\Lambda^{(n)}))$$}
\vskip 12pt  
\noi{\bf Proof:} In order to prove this lemma, it will be enough to check the approximation property.

Let $X$ be a $\F_n(\Lambda)$-module, $Y$ be a $\Nil(\Lambda^{(n)})$-module and $f:\Psi_n(X)\fl Y$ be a morphism. We have to factor $f$ through a 
$\Nil(\Lambda^{(n)})$-module $\Psi_n(X')$ via a homology equivalence $\Psi_n(X')\build\fl_{}^\sim Y$ and a $\F_n(\Lambda)$-morphism $X\fl X'$.

Let $Z$ be the $-1$-cone of $f$. Then the projection $g:Z\fl \Psi_n(X)$ is a $\Nil(\Lambda^{(n)})$-morphism and its $0$-cone is homologically equivalent to $Y$.

There is an integer $m$ such that each element in $\Lambda_n^+$ of weight $\geq m$ acts trivially on $Z$ and $X$. Then $Z$ is a $\Nil(\Lambda^{(n)},m)$-module and,
 because of Lemma 4.11, there is an $\F_m(\Lambda^{(n)})$-module $K$ and a homology equivalence $F_m(K)\build\fl_{}^\sim Z$.

We will construct a $\F_n(\Lambda)$-module $H$ and a morphism $\pi:H\fl X$ such that we have: $\Psi_n(H)=F_m(K)$ and $\Psi_n(\pi)$ is the composite morphism 
$F_m(K)\build\fl_{}^\sim Z\build\fl_{}^g\Psi_n(X)$. Once this is done, we will denote $Y'$ the $0$-cone of $h$ and $\widehat f:X\fl Y'$ the induced morphism. We will 
then have a commutative diagram:
$$\diagram{\Psi_n(X)&\hfl{\widehat f}{}&\Psi_n(Y')\cr\vfl{=}{}&&\vfl{\sim}{}\cr \Psi_n(X)&\hfl{f}{}&Y\cr}$$
where $\Psi(Y')\build\fl_{}^\sim Y$ is a homology equivalence. Waldhausen's Approximation Theorem will then imply the result.

In order to be able to construct this module $H$ we will need to have a better understanding of the algebra $\Lambda_n^+$.
\vskip 12pt

Let $U$ be a basis of $\Lambda$, $\widehat U$ be the set of words in $U$ and $U^+$ be the disjoint union of $U$ and $\{t\}$. Let $V$ be the set of words in $U^+$
that are not zero in $\Lambda_n^+$ and $\overline V$ be the complement of $\widehat U$ in $V$.

Let $u=u_0tu_1t\dots tu_p$ be a word in $U^+$ with $p\geq0$ and $u_i\in\widehat U$ for all $i$. We set $a_i=\partial_2^\circ(u_i)$ and such a word is in $\overline V$ 
if and only if we have: $p>0$, $0\leq a_i<n$ for all $i$ and $a_i>0$ for all $i$ with $0<i<p$. If $u$ is in $\overline{V}$, we set: 
$\lambda(u)=\partial_2^\circ(u_0)=a_0$. We set also: $\lambda(1)=n-1$ and $\lambda$ is defined on $\overline{V}\cup\{1\}$. If $u$ is in $ V$, we set: 
$\mu(u)=\partial_2^\circ(u_p)=a_p$.

Let $\overline W$ be the subset of $\overline V$ defined by:
$$u=u_0tu_1t\dots tu_p\in \overline W\ \ \Longleftrightarrow\ \ 0<\partial_2^\circ(u_i)<n-1\ \hbox{for all}\ i\ \hbox{with}\ 0<i<p$$

We have subsets $\overline V_1$ and ${}_1\overline V_1$ of $\overline V$ defined by:
$$u\in \overline V_1\ \ \Longleftrightarrow\ \ \lambda(u)<n-1,\ \mu(u)=0$$
$$u\in {}_1\overline V_1\ \ \Longleftrightarrow\ \ \lambda(u)=n-1,\ \mu(u)=0$$
and we set: 
$$V_1=\overline V_1\cup\{1\}\hskip 24pt W_1=(\overline V_1\cap\overline W)\cup\{1\}\hskip 24pt {}_1W_1={}_1\overline V_1\cap\overline W$$
Let $a,b$ be two integers with: $0\leq a\leq b<n$. We set:
$$V(a,b)=\{v\in \overline V,\ \hbox{such that}\ \lambda(v)=n-1, \mu(v)\geq a\}\cup\{v\in\widehat U,\ \hbox{such that}\ b\leq\mu(v)<n\}$$
$$W(a,b)=V(a,b)\cap(\widehat U\cup\overline W)$$
and we denote $E(a,b)$ the $\Z$-module generated by $V(a,b)$ and $E_1$ the $\Z$-module generated by $V_1$.
\vskip 12pt
\noi{\bf 4.17 Lemma:} {\sl The module $E_1$ is a right $\Lambda^{(n)}$-module freely generated by $W_1$ and, for any integers $a,b$ with $0\leq a\leq b<n$, $E(a,b)$ is 
a left $\Lambda^{(n)}$-module freely generated by $W(a,b)$. Moreover these modules are stable under the differential.

There is a unique family $\gamma(u,v)\in\Lambda^{(n)}$, for $u,v\in W(0,0)$ such that:
$$\forall u\in W(0,0),\ d(u)=\build\sum_{v\in W(0,0)}^{} \gamma(u,v)v$$
$$\forall u,v\in W(0,0),\ \ d(\gamma(u,v))=-\build\sum_{w\in W(0,0)}^{} \overline{\gamma(u,w)}\gamma(w,v)$$

For any $u\in W_1$, the set $W(u)$ of $v\in\overline{W}$ such that $vu\in\Lambda^{(n)}$ is the set $W(0,n-\lambda(u)-1)\setminus W(n-\lambda(u),n-\lambda(u))$ and the set 
of $v\in\overline{W}$ such that $vu=0$ is the set $W(n-\lambda(u),n-\lambda(u))$.}
\vskip 12pt
\noi{\bf Proof:} We have a unique decomposition: $V_1=W_1{}_1V_1$. Hence, $E_1$ is a right $\Lambda^{(n)}$-module freely generated by $W_1$. For $0\leq a\leq b<n$ we have 
also a unique decomposition: $V(a,b)={}_1 V_1W(a,b)$ and $E(a,b)$ is a left $\Lambda^{(n)}$-module freely generated by $W(a,b)$.

For $u=u_0tu_1t\dots tu_p$ in $V$, we set: $\tau(u)=(a_0,\dots,a_p)$ with $a_i=\partial_2^\circ(u_i)$. For such an element $u$, $d(u)$ is a linear combination of
elements $v$ in some subset $E$ of $V$ and, for each $v\in E$, we have: $\tau(v)=\tau(u)$ or $\tau(v)$ is the sequence obtained from $\tau(u)$ by replacing $a_i,a_{i+1}$
by $a_i+a_{i+1}$ in the sequence (for some $i$ with $0\leq i<p$). Thus we have, for any $v\in E$: $\lambda(v)\geq\lambda(u)$ and $\mu(v)\geq\mu(u)$. 

Therefore, $E_1$ and modules $E(a,b)$ are stable under the differential.

Thus we have unique families $\gamma(u,v)\in\Lambda^{(n)}$, for $u,v\in W(0,0)$ such that we have:

$$\forall u\in W(0,0),\ d(u)=\build\sum_{v\in W(0,0)}^{} \gamma(u,v)v$$
with $\gamma(u,v)$ in $\Lambda^{(n)}$.

Let $u$ be an element in $W(0,0)$. we have:
$$d^2(u)=0=d(\build\sum_{v\in W(0,0)}^{} \gamma(u,v)v)=\sum_{v\in W(0,0)}^{} \Bigl(d(\gamma(u,v))v+\overline{\gamma(u,v)}d(v))$$
$$=\sum_{v\in W(0,0)}^{}d(\gamma(u,v))v+\sum_{v,w\in W(0,0)}^{}\overline{\gamma(u,v)}\gamma(v,w)w$$
and we get the formula:
$d(\gamma(u,v))=-\build\sum_{w\in W(0,0)}^{} \overline{\gamma(u,w)}\gamma(w,v)$ for all $u,v\in W_1$.

Moreover, since $E_1$ is stable under the differential, we have: $\gamma(u,v)\not=0$ implies $\lambda(u)\leq\lambda(v)$.

The last property of the lemma is easy to check.\cqfd
\vskip 12pt
Let $\F^\circ$ be the category of finite graded $\A$-modules. An object of $\F^\circ$ can be seen as a finite $\A$-complex without any differential.
Let $\F$ be the category of $\F^\circ$-modules $X$ equipped with a filtration $$\dots X^{i-1}\subset X^i\subset X^{i+1}\subset\dots$$
such that:

1) for all integer $i$, $X^i/X^{i-1}$ belongs to $\F^\circ$

2) we have: $X^i=0$ for $i\leq 0$ and $X^i=X$ for $i$ large enough

The set of elements of $X^i$ of degree $j$ will be denoted $X^i_j$.

If $k$ is an integer we denote $\F_k$ the subcategory of objects $X\in\F$ such that: $X^k=X$.

If we consider $\A\v$-modules instead of $\A$-modules, we get a category $\F\v$. More precisely, an object of $\F\v$ is a graded
$\A\v$-module $X$ equipped with a filtration
$$\dots X^{i-1}\subset X^i\subset X^{i+1}\subset\dots$$
such that: 

1) each quotient $X^i/X^{i-1}$ is a graded $\A\v$-module and $X^i_j/X^{i-1}_j=0$ for $|j|$ large enough

2) we have: $X^i=0$ for $i\leq 0$ and $X$ is the union of the $X^i$'s.
\vskip 12pt
Let $X$ and $Y$ be two objects of $\F$. If $a\geq 0$ is an integer, we denote $\Map_a(X,Y)$ the set of maps $f\in\Map(X,Y\otimes S^{\otimes a})$ sending each 
$X^i$ to $Y^{i-a}\otimes S^{\otimes a}$. This set is a graded $\Z$-module. An element $f\in\Map_a(X,Y)$ of degree $b$ will be called a $(b,a)$-map from $X$ 
to $Y$. Thus a $(b,a)$-map from $X$ to $Y$ if a map sending each $X^i_j$ to $Y^{i-a}_{j+b}\otimes S^{\otimes a}$.
\vskip 12pt
Let $\H$ be the class of data $(H,\delta,\pi,\varphi,\psi)$ where:

$\bullet$ $H$ associates to each $u\in W_1$ of weight $b$ an $\F_{m-b}$-module $H(u)$ such that $H(1)$ is the underlying $\F_m$-module of $K$. 
The direct sum of all $H(u)$ will also be denoted $H$.

$\bullet$ $\delta $ is a differential of degree $-1$ on $H$ and the components $\delta(u,v):H(u)\fl H(v)$ of $\delta$ respect the filtrations

$\bullet$ $\pi$ is a morphism of degree $0$ from $(H,\delta)$ to $X$ extending the morphism $g_1:K\fl F_1(X)$ and sends each $H(u)$ to $F_{n-\lambda(u)}(X)$

$\bullet$ $\varphi$ associates to each $\theta\in U^+$ and each $u\in W_1$, with $\theta u\in W_1$, a $\partial^\circ(\theta)$-map $\varphi(\theta)$ from $H(u)$ to 
$H(\theta u)$ and $\varphi(\theta)$ is bijective if $\theta=t$

$\bullet$ $\psi$ associates to each $u\in W_1$ and each $v\in W(u)$ a $\partial^\circ(v)$-map $\psi(v)$ from $H(u)$ to $H(1)$
\vskip 12pt
Consider such a data $(H,\delta,\pi,\varphi,\psi)$ in $\H$. Let $u$ be an element of $W_1$. We set: $k=n-\lambda(u)$ and: $E(k)=E(0,k-1)/E(k,k)$. If $a,b$ are two 
integers with $a\geq0$, we denote $T(k)^b_a$ the set of $\Z$-linear maps $f: E(k)\fl H(1)\otimes A[S]$ such that:

$\bullet$ for any homogeneous $u\in\Lambda^{(n)}$ with degree $i$ and any $x\in E(k)$, we have: $f(ux)=(-1)^{ia}uf(x)$. 

$\bullet$ for any $x\in E(k)$ of bidegree $(i,j)$, $f(x)$ belongs to $H(1)^{b-j}_{a+i}\otimes S^{\otimes j}$

Thus we get a graded filtered differential $\Z$-module $T(k)$.

On the other hand, $A$ acts on the left on $H(1)$ and $S$ and this action induces an action on $T(k)$. Therefore, $T(k)$ is a graded filtered differential right $A$-module.

Suppose $T(k)^b_a\not=0$. Then there exist a non-zero morphism $f:E(k)\fl H(1)\otimes A[S]$ sending each element of $E(k)$ of bidegree $(i,j)$ to 
$H(1)^{b-j}_{a+i}\otimes S^{\otimes j}$. Since $H(1)$ is a finite $\A$-complex, there are two integers $a_0\leq b_0$ such that $H(1)$ is null is degree outside of
$[a_0,b_0]$.  Hence, there is an non-zero element $x\in E(k)$ of bidegree $(i,j)$ such that: $a_0\leq a+i\leq b_0$ and $j\leq b$.

Since $\Lambda$ belongs to $\R_0$, the condition $j\leq b$ implies the existence of an integer $c$ such that $i\leq c$ and we get: $a_0-c\leq a\leq b_0$. Therefore,
for each $b$, we have: $T(k)^b_a=0$ for $|a|$ large enough and $T(k)$ belongs to $\F\v$.

Because of Lemma 4.17, $E(k)$ is a free left $\Lambda^{(n)}$-module and the family of maps $\psi(v)$ can be seen as a $A$-linear map $\widehat\psi:H(u)\fl 
T(n-\lambda(u))$ respecting degrees and filtrations. More precisely, the correspondence $\psi\leftrightarrow\widehat\psi$ is given by:
$$\psi(v)(x)=(-1)^{ij}\widehat\psi(x)(v)$$
for any $x\in H(u)$ of degree $i$ and any $v\in W(u)$ of degree $j$.

Therefore, such a data $(H,\delta,\pi,\varphi,\psi)$ in $\H$ is equivalent to the data $(H,\delta,\pi,\varphi,\widehat\psi)$.

For each integer $k\leq n$, we denote $H_k$ the direct sum of the modules $H(u)$ for $u\in W_1$ with $n-\lambda(u)\leq k$.

For each $\theta\in U^+$, we have a linear map $\varphi(\theta):H\fl H\otimes A[S]$ defined as follows:

On $H(u)$, with $u\in W_1$, $\varphi(\theta)$ is the map $\varphi(\theta)$ if $\theta u\in W_1$, the map $\psi(\theta)$ if 
$\theta u\in {}_1W_1\subset\Lambda^{(n)}$ and $\varphi(\theta)=0$ otherwise. 

This correspondence can be extended  to a ring homomorphism $\varphi()$ from $\Lambda_n^+$ to End$(H\build\otimes_A^{} A[S])$.

A data $(H,\delta,\pi,\varphi,\psi)$ in $\H$ will be said to be admissible if the following conditions are verified:

\noi 1) for any $u,v\in W_1$, we have: 
$$\delta(u,v)\not=0\ \Longrightarrow\ u=v\ \ \hbox{or}\ \ \Bigl(\partial_1^\circ(v)<\partial_1^\circ(u)\ \ \ \partial_2^\circ(v)\leq\partial_2^\circ(u)\ \ \ 
\lambda(v)\geq\lambda(v)\Bigr)$$

\noi 2) for all $\theta\in U^+$ with $\partial_1^\circ(\theta)=i$, we have on $H$:
$$\delta\circ\varphi(\theta)=\varphi(d(\theta))+(-1)^i\varphi(\theta)\circ \delta$$

\noi 3) for all $u\in W_1$ and all $v\in W(u)$ with $\partial_1^\circ(v)=i$, we have on $H(u)$:
$$d\circ\psi(v)=\build\sum_{w\in W(u)}^{}\gamma(v,w)\psi(w)+(-1)^i\psi(v)\circ\delta$$

\noi 4) for all $\theta\in U^+$ we have: $\pi\circ\varphi(\theta)=\theta\circ\pi$

\noi 5) for all $u\in W_1$, all $\theta\in U^+$ with $\theta u\in W_1$ and all $v\in W(\theta u)$ we have on $H(u)$: 
$$\psi(v\theta)=\psi(v)\circ\varphi(\theta)$$

It is easy to see that condition 5) is equivalent to the commutativity of the diagram:
$$\diagram{H(v)&\hfl{\varphi(\theta)}{}&H(u)\otimes S^{\otimes j}\cr\vfl{\widehat\psi}{}&&\vfl{\widehat\psi}{}\cr 
T(n-\lambda(v))&\hfl{\theta^*}{}&T(n-\lambda(u))\otimes S^{\otimes j}\cr}$$
where $j$ is the weight of $\theta$ and $\theta^*$ the morphism defined by:
$$\theta^*(f)(v)=(-1)^{ij}f(v\theta)$$
for $f\in T(n-\lambda(v))$, $v\in E(n-\lambda(u))$, $i=\partial_1^\circ(\theta)$ and $j=\partial^\circ(f)+\partial_1^\circ(v)$.

We can see also that condition 3) is equivalent to the fact that $\widehat\psi$ induces, for each integer $k<n$, a morphism:
$$\widehat\psi: H_k\fl T(k)$$
respecting degrees, filtrations and differentials.

\vskip 12pt
We will construct an admissible data in $\H$ by induction. We start by equipping $W_1$ with a total order such that we have for all $u$ and $v$ in $W_1$:
$$u\leq v\ \Longrightarrow \Big(\partial_2^\circ(u)<\partial_2^\circ(v)\Big)\ \hbox{or}\ \Big(\partial_2^\circ(u)=\partial_2^\circ(v)\ \hbox{and}\ \partial_1^\circ(u)
\leq\partial_1^\circ(v)\Big)$$

Suppose that, for a certain $u\not=1$ in $W_1$, the modules $H(v)$ are constructed for all $v<u$ as well as all the morphisms $\delta$, $\pi$,
$\varphi$ and $\psi$ occurring between modules already defined and that all conditions of admissibility are verified.

We have a decomposition: $u=\theta v$ with $\theta\in U^+$ and $v\in W_1$. Let $H(<\! u)$ be the direct sum of the modules $H(v)$ for $v< u$ and 
$\lambda(v)\geq\lambda(u)$. We denote also $H(\leq u)$ the module $H(u)\oplus H(<\! u)$.

We have to define $H(u)$, the differential $\delta:H(u)\fl H(\leq u)$, the projection $\pi:H(u)\fl X_{n-\lambda(u)}$, the map $\varphi(\theta)$ from $H(v)$ to
$H(u)\otimes S^{\otimes *}$ and $\partial^\circ(w)$-maps $\psi(w)$ from $H(u)$ to $H(1)$ for all $w\in W(u)$. After that, we have to check the admissibility conditions.

Suppose: $\theta=t$. In this case, we define $H(u)=H(tv)$ as the suspension of $H(v)$ (i.e. the $1$-cone of $0\fl H(v)$). Let $\varepsilon$ be the identity map $H(v)\fl 
H(tv)$. Then morphisms $\delta:H(u)\fl H(<\!u)$, $\pi:H(u)\fl X_n$, $\varphi(t):H(v)\fl H(u)$ and $\psi(w):H(u)\fl H(1)\otimes S^{\otimes*}$ are defined as 
follows:
$$\delta=\varepsilon^{-1}-\varphi(t)\delta\circ\varepsilon^{-1}\hskip 24pt\pi=t\pi\circ\varepsilon^{-1}\hskip 24pt\varphi(t)=\varepsilon$$
$$\forall w\in W(u),\ \ \psi(w)=\psi(wt)\circ\varepsilon^{-1}$$
and the admissibility conditions are easy to check.

Suppose $\theta$ is in $U$. We set: $k=n-\lambda(u)$, $a=\partial_2^\circ(u)$ and $(\alpha,\beta)=\partial^\circ(\theta)$. We have: $v\not=1$ and: $1<k<n$. We have to 
construct the $\F_{m-a}$-module $H(u)$ and morphisms:
$$\varphi(\theta):H(v)\fl H(u)\otimes S^{\otimes \beta}\hskip 24pt \widehat\psi:H(u)\fl T(k)$$
$$\delta: H(u)\fl H(\leq u)\hskip 24pt \pi:H(u)\fl F_k(X)$$

If $U$ is a graded filtered differential module, its graded filtered module of cycles will be denoted $U^\bullet$.

Thus we want to construct a commutative diagram:
$$\diagram{H(u)_i&\hfl{\delta}{}&H(u)^\bullet_{i-1}\oplus H(<\! u)^\bullet_{i-1}\cr \vfl{\widehat\psi\oplus\pi}{}&&\vfl{\widehat\psi\oplus\pi}{}\cr 
T(k)^{m-a}_i\oplus F_k(X)_i&\hfl{d}{}& {T(k)^\bullet}^{m-a}_{i-1}\oplus F_k(X)^\bullet_{i-1}\cr}\leqno{(D(u)_i)}$$
where $F_k(X)$ is considered as a $\F_1$-module.

Since $T(k)$ is in $\F\v$, there are two integers $a_0<b_0$ such that we have: $T(k)^{m-b}_i=F_k(X)_i={T(k)^\bullet}^{m-b}_{\ i-1}=H(<\! u)^\bullet_{i-1}=0$ for all 
$i<a_0$ and all $i>b_0$.

The module $H(u)$ will be constructed by induction on the degree. We start by setting: $H(u)_i=0$ for all $i<a_0$. Suppose $H(u)$ is constructed in all degree $<i$ as well
as the morphisms $\widehat\psi$ and $\pi$ on these modules. Then, in the diagram $(D(u)_i)$, everything is constructed except $H(u)_i$ and we can take the fiber product 
$M$ defined by the cartesian diagram:
$$\diagram{M&\hfl{h_0}{}&H(u)^\bullet_{i-1}\oplus H(<\! u)^\bullet_{i-1}\cr \vfl{h_1}{}&&\vfl{\widehat\psi\oplus\pi}{}\cr 
T(k)^{m-a}_i\oplus F_k(X)_i&\hfl{d}{}& {T(k)^\bullet}^{m-a}_{i-1}\oplus F_k(X^\bullet)_{i-1}\cr}\leqno{(D_i)}$$
and $M$ is a $\F\v$-module.

On the other hand, we have a commutative diagram:
$$\diagram{H(v)_j&\hfl{\delta}{}&H(v)^\bullet_{j-1}\oplus H(<\! v)^\bullet_{j-1}\cr \vfl{\widehat\psi\oplus\pi}{}&&\vfl{\widehat\psi\oplus\pi}{}\cr 
T(k+\beta)^{m-b}_j\oplus F_{k+\beta}(X)_j&\hfl{d}{}& {T(k+\beta)^\bullet}^{m-b}_{j-1}\oplus F_{k+\beta}(X^\bullet)_{j-1}\cr}\leqno{(D(v)_j)}$$
with $b=a-\beta$.

Since $S$ is flat on the left, the diagram $(D_i)\otimes S^{\otimes \beta}$ is cartesian. Then the action of $\theta$: $F_{k+\beta}(X)_{i-\alpha}\fl F_k(X)_i\otimes 
S^{\otimes\beta}$, the morphism $\varphi(\theta)$: $H(w)\fl H(w\theta)\otimes S^{\otimes\beta}$ and the morphism $\theta^*$: $T(k+\beta)^{m-b}_{i-\alpha}\fl 
T(k)^{m-a}_i\otimes S^{\otimes\beta}$ induce a morphism of diagrams from $(D(v)_{i-\alpha})$ to $(D_i\otimes S^{\otimes \beta})$ and then a $\F\v$-morphism 
$\varphi:H(v)_{i-\alpha}\fl M\otimes S^{\otimes \beta}$ sending, for each integer $j$, $H(v)^j_{i-\alpha}$ to $M^{j-\beta}\otimes S^{\otimes \beta}$.

Let $P$ be the image of $\varphi$. Since each $H(v)^j_{i-\alpha}$ is an $\A$-module, each $P^j$ is finitely generated and there exists a filtered module $P'\subset M$ 
such that, for each $j$, $P'^j$ is finitely generated and $\varphi(H(v)^j)$ is contained in $P'^{j-\beta}\otimes S^{\otimes\beta}$. 

Let $u: Q\fl P'$ be an $\A$-epimorphism  where $Q$ is a filtered $\A$-module such that, for all $j$, we have: $Q^j=0$ is $j\leq 0$, $Q^j=Q$ is $j\geq m-b$, $Q^j/Q^{j-1}$
belongs to $\A$ and $u:Q^j\fl P'^j$ is onto. If the morphism $\varphi$ is zero, the filtered module $Q$ will be chosen to be zero.

We set: $H(u)_i=Q$ and we have morphisms: 
$$\delta=h_0\circ u:H(u)_i\fl H(u)'_{i-1}\oplus H(<u)'_{i-1}$$
$$\widehat\psi\oplus \pi=h_1\circ u:H(u)_i\fl T(k)^{m-a}_i\oplus X^k_i$$

Moreover the morphism $\varphi:H(v)_{i-\alpha}\fl M\otimes S^{\otimes \beta}$ factors through $Q\otimes S^{\otimes \beta}$ via a morphism
$\varphi(\theta):H(v)_{i-\alpha}\fl H(u)_i\otimes S^{\otimes \beta}$. Therefore we get a morphism of diagrams $(D(v)_{i-\alpha})\fl (D(u)_i\otimes 
S^{\otimes\beta})$.

Thus we can construct $H(u)$ by induction on the degree and, because we have: $H(u)_i=0$ for $i$ large enough, $H(u)$ belongs to $\F_{m-a}$. The admissibility conditions
are easy to check. Hence, we can construct $H(u)$ by induction on $u$, for all $u\in W_1$ with $\partial_2^\circ(u)<m$. If $\partial_2^\circ(u\geq m$, we set: 
$H(u)=0$.

The module $H$, its filtration $\dots\subset H^k\subset H^{k+1}\subset\dots$ and the morphism $t$ define an $\F_n(\Lambda)^1$-module still denoted  $H$.

Thus we have an $\F_n(\Lambda)^1$-module $H$ and a morphism $\pi:H\fl X$ such that: $\Psi_n(H)=F_m(K)$ and $\Psi_n(h)$ is the composite morphism 
$F_m(K)\build\fl_{}^\sim Z\build\fl_{}^g\Psi_n(X)$. 

Let $Y'$ be the $0$-cone of $h$ and $\widehat f:X\fl Y'$ be the induced morphism. We have a commutative diagram:
$$\diagram{\Psi_n(X)&\hfl{\widehat f}{}&\Psi_n(Y')\cr\vfl{=}{}&&\vfl{\sim}{}\cr \Psi_n(X)&\hfl{f}{}&Y\cr}$$
where $\Psi_n(Y')\build\fl_{}^\sim Y$ is a homology equivalence. Waldhausen's Approximation Theorem will then implies Lemma 4.16.\cqfd
\vskip 12pt
\noi{\bf 4.18 Lemma:} {\sl Let $\Lambda$ be an algebra in $\R_0$ and $n>1$ be an integer. Then there is a homotopy fibration:
$$Nil(\Lambda^{(n)})\fl Nil(\Lambda,n-1)\build\fl_{}^f Nil(\Lambda,n)$$
where $f$ is induced by the inclusion $\Nil(\Lambda,n-1)\subset\Nil(\Lambda,n)$.}
\vskip 12pt
\noi{\bf Proof:} 
As before, we have a decomposition of infinite loop spaces:
$$K(\F_n(\Lambda)^0)\simeq K(\A)^{n-1}\times F_n^\circ(\Lambda)$$
Because of the Additivity Theorem, we have a homotopy equivalence: $K(\F_n(\Lambda))\simeq K(\A)^n$ and Lemma 4.12 implies a homotopy equivalence of infinite loop spaces:
$$\Omega(Nil(\Lambda,n))\simeq F_n^\circ(\Lambda)$$

Therefore, lemmas 4.13, 4.14 and 4.16 imply a homotopy fibration of infinite loop spaces:
$$\Omega(Nil(\Lambda,n-1))\fl \Omega(Nil(\Lambda,n))\fl Nil(\Lambda^{(n)})$$
and $Nil(\Lambda^{(n)})$ is the homotopy fiber of the map $Nil(\Lambda,n-1)\fl Nil(\Lambda,n)$ induced by the inclusion: $\Nil(\Lambda,n-1)\subset\Nil(\Lambda,n)$.\cqfd
\vskip 12pt
\noi{\bf 4.19 Lemma:} {\sl For any algebra $\Lambda$ in $\R_0$ and any integer $n>0$, the space $Nil(\Lambda,n)$ is contractible.}
\vskip 12pt
\noi{\bf Proof:} Let $p\geq 0$ be an integer and $\E_p$ the class of algebras $\Lambda$ in $\R_0$ such that, for any $n>0$, $Nil(\Lambda,n)$ is $p$-connected. Because of
Corollary 4.10, each algebra in $\R_0$ belongs to $\E_0$.

Let $p>0$ be an integer such that $\R_0\subset \E_{p-1}$. Let $\Lambda$ be an algebra in $\R_0$. Since any algebra $\Lambda^{(n)}$ are in $\E_{p-1}$, each space 
$Nil(\Lambda)$ is $(p-1)$-connected and we have, for any $n>1$, an exact sequence:
$$Nil_p(\Lambda^{(n)})\fl Nil_p(\Lambda,n-1)\fl Nil_p(\Lambda,n)\fl 0$$
Therefore, each morphism $Nil_p(\Lambda,n-1)\fl Nil_p(\Lambda,n)$ is onto. But $Nil(\Lambda,1)$ is contractible. Hence, we have: $Nil_p(\Lambda,n)=0$ for any $n>0$
and $\Lambda$ belongs to $\E_p$.

Thus we see, by induction, that each space $Nil(\Lambda,n)$ is $p$-connected for any $p$. The result follows.\cqfd
\vskip 12pt
\noi{\bf 4.20 Corollary:} {\sl For any integer $n>0$, the space $Nil(A,S,n)$ is contractible.}
\vskip 12pt
\noi{\bf Proof:} This is a straightforward consequence of lemmas 1.17, 4.1 and 4.19.\cqfd 
\vskip 12pt
\noi{\bf Proof of Theorem 2:} Because of Corollary 4.20, for any right-regular ring $A$ and any left-flat $A$-bimodule $S$, the space $Nil(A,S)$ is contractible.

In [V], an $\Omega$-spectrum $\underline Nil(A,S)$ is defined and the connective part of it is the space $Nil(A,S)$. For any $i\in\Z$ we set: 
$Nil_i(A,S)=\pi_i(\underline Nil(A,S))$ and we have a Bass-Quillen exact sequence (see [B], [Q], Lemma 2.13 in [V]):
$$0\rightarrow Nil_i(A,S)\rightarrow Nil_i(A_1,S_1)\oplus Nil_i(A_2,S_2)\rightarrow Nil_i(A_3,S_3)\rightarrow Nil_{i-1}(A,S)\rightarrow 0$$
with:
$$A_1=A[t]\hskip 24pt S_1=S[t]$$
$$A_2=A[t^{-1}]\hskip 24pt S_2=S[t^{-1}]$$
$$A_3=A[t^{\pm1}]\hskip 24pt S_3=S[t^{\pm1}]$$

Let $\B$ be the class of pairs $(A,S)$ where $A$ is a right-regular ring and $S$ is a left-flat $A$-bimodule. If $k$ is an integer, we denote $\B(k)$ the class of objects 
in $\B$ defined by:
$$(A,S) \in\B(k)\ \ \Longleftrightarrow\ \ \ \forall i\geq k,\  Nil_i(A,S)=0$$

We just have proven that $\B\subset \B(0)$.

Let $k<0$ be an integer. Suppose we have: $\B\subset \B(k+1)$.

Let $(A,S)$ be an object in $\B$. Consider the corresponding Bass-Quillen exact sequence. Because of Proposition 3.5 and Proposition 3.6, each $(A_i,S_i)$ belongs to $\B$
and the exact sequence in degree $i=k+1$ reduces to $0\fl Nil_k(A,S)\fl 0$. Hence, $(A,S)$ belongs to $\B(k)$ and we have: $\B\subset \B(k)$. Thus $\B$ is contained in
each $\B(k)$. Therefore, for any $(A,S)\in\B$, $\underline Nil(A,S)$ is contractible and Vanishing Theorem 2 is proven.\cqfd

\vskip 24pt
\noi{\bf References: }

\begin{list}{}{\leftmargin 24pt \labelsep 10pt \labelwidth 40pt \itemsep 0pt} 
\item[{[AMM]}] B. Antieau, A. Mathew, and M. Morrow, {\sl The K-theory of perfectoid rings}, Doc. Math. 27 (2022), 1923–1953.
\item[{[B]}] H. Bass -- {\sl Algebraic K-theory}, Benjamin (1968).
\item[{[Ga]}] P. Gabriel, {\sl Sur les catégories abéliennes}, Bull. Soc. Math. France, {\bf 90}, 1962, 323--448.
\item[{[Ge]}] S. M. Gersten, {\sl K-theory of free rings}, Comm. Algebra 1 (1974), 39–64.
\item[{[Gl]}] S. Glaz, {\sl Commutative coherent rings}, Lecture Notes in Mathematics, vol. 1371, Springer-Verlag, Berlin, 1989.
\item[{[GR]}] P. Gabriel and A. V. Roiter, {\sl Representations of Finite-Dimensional Algebras}, Encyclopaedia of Mathematical Sciences, {\bf 73}, 1997.
\item[{[Ka]}] M. Karoubi, -- {\sl Alg\`ebres de Clifford et K-th\'eorie}. Ann. Sci. Éc. Norm. Sup\'er. {\bf 1} (2), 1968, 161–-270.
\item[{[Ke]}] B. Keller, {\sl Chain complexes and stable categories}, Manuscripta Mathematica. {\bf 67}, 1990, 379–-417. doi:10.1007/BF02568439.
\item[{[Q]}] D. Quillen -- {\sl Higher algebraic K-theory I}, Proc. Conf. alg. K-theory, Lecture Notes in Math. {\bf 341} (1973), 85--147.
\item[{[Na]}] C. Dahlhausen -- {\sl Regularity of semi-valuation rings and homotopy invariance of algebraic K-theory}, 2024, Arxiv: math.KT 2403.02413.
\item[{[Ne]}] A. Neeman, {\sl The homotopy category of flat modules and Grothendieck duality}, Inv. Math.{\bf 174}, 2008, 255--308.
\item[{[TT]}] R.W. Thomason and T. Trobaugh, {\sl Higher algebraic K-theory of schemes and of derived categories}, The Grothendieck Festschrift III, Progress in Math., 
vol. 88, Birkh\"auser, (1990), 247–-435.
\item[{[V]}] P. Vogel, {\sl Algebraic K-theory of generalized free products and functors Nil}, (2019) Arxiv: math.KT 1909.02413. J. Pure Appl. Alg. 
{\bf 225} Issue 2 (2021) doi:10.1016/j.jpaa.2020.106488.
\item[{[W1]}] F. Waldhausen -- {\sl Algebraic K-theory of generalized free products, part 1 \& 2}, Annals of Math. {\bf 108} (1978), 135--256.
\item[{[W2]}] F. Waldhausen -- {\sl Algebraic K-theory of spaces}, Algebraic and geometric topology (New Brunswick, N.J., 1983), 318–419, 
  Lecture Notes in Math., 1126, Springer, Berlin, (1985).
\end{list}
\end{document}